%% file: thesis-arxiv.tex
\documentclass[10pt, a4paper, twoside, titlepage]{article}
\usepackage[a4paper, inner=1.15in, outer=1.35in, bottom=1.25in, marginparwidth=1in, marginparsep=0.125in]{geometry}
\pdfoutput=1

\usepackage[ngerman, english]{babel} 
\usepackage{amsmath}
\usepackage{amssymb}
\usepackage{amsthm}
\usepackage{thmtools}
\usepackage{algorithm, algpseudocode}
\usepackage{bbm}
\usepackage{blkarray}
\usepackage{pgfplots}
\usetikzlibrary{graphs,arrows.meta,matrix,fit,patterns,shapes}
\pgfplotsset{compat=1.8}
\usepackage{makecell}
\usepackage{svg}
\usepackage{colonequals}
\usepackage{soul}
\usepackage{color}
\usepackage{mathtools}
\usepackage{hyperref}
\usepackage[nameinlink,capitalize]{cleveref}
\usepackage[toc,acronym,symbols]{glossaries}
\usepackage{stmaryrd}
\usepackage{eso-pic}
\usepackage{transparent}
\usepackage{graphicx}
\usepackage[framemethod=TikZ]{mdframed}
\usepackage[shortlabels]{enumitem}
\usepackage[nottoc,numbib]{tocbibind}
\usepackage{ifoddpage}
\usepackage{lmodern} 
\usepackage[T1]{fontenc}
\usepackage{FiraSans}
\usepackage{merriweather} 
\usepackage{FiraMono}
\usepackage{mathpazo} 

\usepackage{fancyhdr}
\usepackage{titlesec} 
\numberwithin{equation}{section} 

\definecolor{codegreen}{RGB}{140,153,0}
\definecolor{codegray}{RGB}{128,128,128}
\definecolor{codepurple}{RGB}{173,20,87}
\definecolor{codeorange}{RGB}{240,107,0}
\definecolor{backcolour}{RGB}{237,237,237}
\definecolor{backcolour2}{RGB}{26,26,26}
\definecolor{solarizedbg}{RGB}{253,246,227}
\definecolor{solarizedtext}{RGB}{101,123,131}
\definecolor{codeblue}{RGB}{3,135,209}

\definecolor{tumBlue}{RGB}{0,101,189} 
\definecolor{tumDarkBlue}{RGB}{0,82,147} 
\definecolor{tumLightBlue}{RGB}{100,160,200} 
\definecolor{tumLighterBlue}{RGB}{152,198,234} 
\definecolor{tumOrange}{RGB}{227,114,34} 
\definecolor{tumGreen}{RGB}{162,173,0} 
\definecolor{tumGray}{RGB}{153,153,153} 
\definecolor{tumLight}{RGB}{218,215,203} 

\colorlet{sectionblue}{tumBlue}
\definecolor{linkred}{RGB}{127,0,0} 
\definecolor{darklinkred}{RGB}{50,0,0} 
\colorlet{headingcolor}{sectionblue}
\colorlet{headingcolormuted}[RGB]{headingcolor!20!white}
\colorlet{linkcolor}{linkred}

\makeatletter
\DeclareRobustCommand{\xhl}[1]{%
 \protected@edef\tmp{#1}%
 \expandafter\hl\expandafter{\tmp}%
}
\makeatother
\newcommand{\hlyellow}[1]{{\sethlcolor{solarizedbg}\protect\xhl{#1}}}
\newcommand*{\yellowemph}[1]{%
  \ifmmode\mathchoice{\tikz[baseline=(text.base)]\node(text)[rectangle, fill=solarizedbg, inner sep=0pt]{$\displaystyle #1$};}{\tikz[baseline=(text.base)]\node(text)[rectangle, fill=solarizedbg, inner sep=0pt]{$\textstyle #1$};}{\tikz[baseline=(text.base)]\node(text)[rectangle, fill=solarizedbg, inner sep=0pt]{$\scriptstyle #1$};}{\tikz[baseline=(text.base)]\node(text)[rectangle, fill=solarizedbg, inner sep=0pt]{$\scriptscriptstyle #1$};}\else\hlyellow{#1}\fi
}

\usepackage[labelfont={color=headingcolor,bf}]{caption}

\definecolor{darkblue}{RGB}{0,102,204}
\colorlet{colorlow}[hsb]{darkblue}
\definecolor{darkred}{RGB}{204,26,0}
\colorlet{colorhigh}[hsb]{darkred}
\colorlet{color1}[RGB]{colorlow!100!colorhigh}
\colorlet{color2}[RGB]{colorlow!85!colorhigh}
\colorlet{color3}[RGB]{colorlow!70!colorhigh}
\colorlet{color4}[RGB]{colorlow!50!colorhigh}
\colorlet{color5}[RGB]{colorlow!30!colorhigh}
\colorlet{color6}[RGB]{colorlow!15!colorhigh}
\colorlet{color7}[RGB]{colorlow!0!colorhigh}

\makeatletter
\def\moverlay{\mathpalette\mov@rlay}
\def\mov@rlay#1#2{\leavevmode\vtop{%
   \baselineskip\z@skip \lineskiplimit-\maxdimen
   \ialign{\hfil$\m@th#1##$\hfil\cr#2\crcr}}}
\newcommand{\charfusion}[3][\mathord]{
    #1{\ifx#1\mathop\vphantom{#2}\fi
        \mathpalette\mov@rlay{#2\cr#3}
      }
    \ifx#1\mathop\expandafter\displaylimits\fi}
\makeatother

\newcommand{\midd}[1]{\mathrel{}\middle#1\mathrel{}}

\mdfsetup{ skipbelow=-0.2em }

\newlength{\spaceblength}
\declaretheoremstyle[
    headfont=\bfseries,
    notefont=\bfseries,
    notebraces={}{\\[\parskip]}, 
    bodyfont=\normalfont\upshape,
    headpunct={},
    postheadspace=\spaceblength,
    spacebelow=\parskip,
    spaceabove=\parskip,
    headformat={%
        \checkoddpage\ifoddpage\rlap{\hskip\textwidth\hskip10pt\color{headingcolor}\ \ \ \NAME\ \NUMBER}\hskip-\spaceblength{\NOTE}%
        \else\makebox[0pt][r]{\color{headingcolor}\NAME\ \NUMBER\hskip10pt\ \ \ \ }\hskip-\spaceblength{\NOTE}\fi%
    },
    mdframed={
        nobreak=true,
        linecolor=headingcolor!20,
        innertopmargin=8pt,
        innerbottommargin=6pt,
        innerleftmargin=8pt,
        innerrightmargin=8pt,
        skipabove=0.7em,
        linewidth=2pt } 
]{boxstyle}
\declaretheoremstyle[
    headfont=\bfseries\itshape,
    notefont=\normalfont\bfseries,
    notebraces={}{\\[\parskip]}, 
    bodyfont=\normalfont,
    headpunct={\indent},
    postheadspace=\spaceblength,
    spacebelow=\parskip,
    spaceabove=\parskip,
    headformat={%
        \checkoddpage\ifoddpage\rlap{\hskip\textwidth\color{headingcolor}\ \ \ \NAME}\hskip-\spaceblength{\NOTE}%
        \else\makebox[0pt][r]{\color{headingcolor}\NAME\ \ \ \ }\hskip-\spaceblength{\NOTE}\fi%
    },
    qed=\qedsymbol
]{proofstyle}
\declaretheoremstyle[
    headfont=\bfseries,
    notefont=\bfseries,
    notebraces={}{\\[\parskip]}, 
    bodyfont=\normalfont,
    headpunct={\indent},
    postheadspace=\spaceblength,
    spacebelow=\parskip,
    spaceabove=\parskip,
    headformat={%
        \checkoddpage\ifoddpage\rlap{\hskip\textwidth\color{headingcolor}\ \ \ \NAME\ \NUMBER}\hskip-\spaceblength{\NOTE}%
        \else\makebox[0pt][r]{\color{headingcolor}\NAME\ \NUMBER\ \ \ \ }\hskip-\spaceblength{\NOTE}\fi%
    },
    qed={$\begingroup\color{headingcolormuted}\blacktriangleleft\endgroup$}
]{examplestyle}
\declaretheoremstyle[
    headfont=\normalfont,
    notefont=\bfseries,
    notebraces={}{\\[\parskip]}, 
    bodyfont=\normalfont,
    headpunct={\indent},
    postheadspace=\spaceblength,
    spacebelow=\parskip,
    spaceabove=\parskip,
    headformat={%
        \checkoddpage\ifoddpage\rlap{\hskip\textwidth\color{headingcolor}\ \ \ \NAME}\hskip-\spaceblength{\NOTE}%
        \else\makebox[0pt][r]{\color{headingcolor}\NAME\ \ \ \ }\hskip-\spaceblength{\NOTE}\fi%
    },
    qed={$\begingroup\color{headingcolormuted}\blacktriangleleft\endgroup$}
]{remarkstyle}

\declaretheorem[style=boxstyle,numberwithin=section]{definition}
\declaretheorem[style=boxstyle,sibling=definition]{theorem}

\declaretheorem[style=boxstyle,sibling=definition]{lemma}

\declaretheorem[style=examplestyle,sibling=definition]{example}
\declaretheorem[style=proofstyle,numbered=no,name=Proof]{tproof}

\declaretheorem[style=remarkstyle,numbered=no]{remark}

\addto\extrasenglish{

}

\newlength{\secskip}
\titleformat{\section}
  {\LARGE\bfseries\color{headingcolor}\sffamily}
  {\makebox[0pt][r]{\thesection\hspace*{\secskip}}}{0em}{}
\titleformat{\subsection}
  {\large\bfseries\color{headingcolor}\sffamily}
  {\makebox[0pt][r]{\thesubsection\hspace*{\secskip}}}{0em}{}
\titleformat{\subsubsection}
  {\bfseries\color{headingcolor}\sffamily}
  {\makebox[0pt][r]{\thesubsubsection\hspace*{\secskip}}}{0em}{}

\hypersetup{
    colorlinks=true,
    allcolors=.,
    urlcolor=linkcolor,
    citecolor=linkcolor,
    linkcolor=linkcolor,
}

\newcommand{\spectitle}[3]{
	\sffamily
	\def\svgwidth{60pt}
  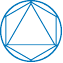\hfill\raisebox{5pt}[0pt][0pt]{\def\svgwidth{100pt}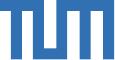}\\[8em]
	{\large\merriweathersanslight Technical University of Munich\\[2em]
	Department of Mathematics\\[2em]
	Bachelor's Thesis in Mathematics}\\[8em]
	{\color{headingcolor}\Huge\bfseries #1}\\[4em]
	{\LARGE #2}\\[10em]
	{\large\merriweathersanslight {\sffamily Supervisor:} Prof.~Dr.~rer.~nat.~habil.~Nina Gantert\\[1.5em]
	{\sffamily Submission Date:} #3}
}

\newtagform{roman}[]()

\DeclareMathOperator{\VAR}{Var}

\DeclareMathOperator{\BR}{br}

\DeclareMathOperator{\STRENGTH}{Strength}

\newcommand{\Strength}[1] {\STRENGTH\left(#1\right)}
\newcommand{\br}[1] {\BR\left(#1\right)}

\newcommand{\Prb}[1] {\mathbb{P}\left[#1\right]}
\newcommand{\Prbc}[2] {\mathbb{P}_{#1}\left[#2\right]}
\newcommand{\Prbcs}[2] {\mathbb{P}^{#1}\left[#2\right]}
\newcommand{\Prbcc}[3] {\mathbb{P}_{#1}^{#2}\left[#3\right]}
\newcommand{\bPrb}[1] {\mathbf{P}\left[#1\right]}

\newcommand{\Ex}[1] {\mathbb{E}\left[#1\right]}
\newcommand{\bEx}[1] {\mathbf{E}\left[#1\right]}
\newcommand{\Var}[1] {\VAR\left[#1\right]}

\newcommand{\dx}[1] {\;\mathrm{d}#1}
\newcommand{\mudx}[2] {\;#1\left(\mathrm{d}#2\right)}

\newcommand{\bbP} {\mathbb{P}}
\newcommand{\bbE} {\mathbb{E}}

\newcommand{\bbR} {\mathbb{R}}
\newcommand{\bbN} {\mathbb{N}}
\newcommand{\bbZ} {\mathbb{Z}}

\newcommand{\bP} {\mathbf{P}}

\DeclareMathOperator*{\argmin}{arg\,min}
\newcommand{\divides}{\,\big\rvert}

\tikzset{
  diagonal fill a/.style n args=2{path picture={%
  \draw[fill=#1, draw=none] (path picture bounding box.south west) --
              (path picture bounding box.north east) -- (path picture bounding box.south east) -- cycle; %
  \draw[fill=#2, draw=none] (path picture bounding box.south west) --
              (path picture bounding box.north east) -- (path picture bounding box.north west) -- cycle;}},
  diagonal fill b/.style n args=2{path picture={%
  \draw[fill=#1, draw=none] (path picture bounding box.south west) --
              (path picture bounding box.north east) -- (path picture bounding box.south east) -- cycle; %
  \draw[pattern=north west lines, pattern color=#2, draw=none] (path picture bounding box.south west) --
              (path picture bounding box.north east) -- (path picture bounding box.north west) -- cycle;}},
  diagonal fill c/.style n args=2{path picture={%
  \draw[pattern=north west lines, pattern color=#1, draw=none] (path picture bounding box.south west) --
              (path picture bounding box.north east) -- (path picture bounding box.south east) -- cycle; %
  \draw[fill=#2, draw=none] (path picture bounding box.south west) --
              (path picture bounding box.north east) -- (path picture bounding box.north west) -- cycle;}},
  diagonal fill d/.style n args=2{path picture={%
  \draw[pattern=north west lines, pattern color=#1, draw=none] (path picture bounding box.south west) --
              (path picture bounding box.north east) -- (path picture bounding box.south east) -- cycle; %
  \draw[pattern=north west lines, pattern color=#2, draw=none] (path picture bounding box.south west) --
              (path picture bounding box.north east) -- (path picture bounding box.north west) -- cycle;}},
  table nodes/.style={
    rectangle,
    draw=none,
    align=center,
    minimum height=7mm,
    text depth=0.5ex,
    text height=2ex,
    inner xsep=0pt,
    outer sep=0pt
  },      
  table/.style={
    matrix of nodes,
    row sep=-\pgflinewidth,
    column sep=-\pgflinewidth,
    nodes={
        table nodes
    }
  }
}

\usepackage{ifthen}
\newboolean{print}
\setboolean{print}{false}
\usepackage{import}

\title{Linearly Edge-Reinforced Random Walks}
\author{Fabian Michel}
\date{14/09/2020}

\newglossaryentry{root}{name=\ensuremath{\rho},description={root of a tree},type=symbols,sort={rho}}
\newglossaryentry{par}{name=\ensuremath{\protect\overleftarrow{v}},description={parent of a node $v$},type=symbols,sort={v}}
\newglossaryentry{children}{name=\ensuremath{C\left(v\right)},description={children of a node $v$},type=symbols,sort={cv}}
\newglossaryentry{child}{name=\ensuremath{v_i},description={$i$-th child of a node $v$},type=symbols,sort={vi}}
\newglossaryentry{numch}{name=\ensuremath{b_v},description={number of children of a node $v$},type=symbols,sort={bv}}
\newglossaryentry{paredge}{name=\ensuremath{e_v},description={edge from a node $v$ to its parent},type=symbols,sort={ev}}
\newglossaryentry{nodedist}{name=\ensuremath{\left|v\right|},description={distance of a node $v$ to the root},type=symbols,sort={v}}
\newglossaryentry{edgedist}{name=\ensuremath{\left|e\right|},description={distance of an edge $e$ to the root (equals the distance of the endpoint further away from the root)},type=symbols,sort={e}}
\newglossaryentry{distset}{name=\ensuremath{T_n},description={set of nodes at distance $n$ from the root},type=symbols,sort={tn}}
\newglossaryentry{ancestor}{name=\ensuremath{u < v},description={node $u$ is an ancestor of node $v$},type=symbols,sort={uv}}
\newglossaryentry{neighbor}{name=\ensuremath{u \sim v},description={node $u$ is a neighbor of node $v$},type=symbols,sort={uv}}
\newglossaryentry{lca}{name=\ensuremath{u \land v},description={lowest common ancestor of $u$ and $v$},type=symbols,sort={uv}}

\newglossaryentry{ntree}{name=\ensuremath{T^{\left(n\right)}},description={tree constructed from another tree $T$ as follows: $T^{\left(n\right)}$ contains all nodes $v$ of $T$ with $n \divides \left|v\right|$ where two nodes $u, v$ are connected by an edge if, and only if, $u \leq v$ and $\left|u\right| + n = \left|v\right|$ (in the tree $T$)},type=symbols,sort={tn}}

\newglossaryentry{omuv}{name=\ensuremath{\omega\left(u,v\right)},description={probability to go from node $u$ to node $v$ in the environment $\omega$},type=symbols,sort={omegauv}}
\newglossaryentry{pmeasom}{name=\ensuremath{\bbP_\omega^v},description={probability measure of the Markov chain starting at node $v$ in the environment $\omega$},type=symbols,sort={pomegav}}
\newglossaryentry{pmeas}{name=\ensuremath{\bbP^v},description={probability measure of the mixed Markov chains starting at node $v$},type=symbols,sort={pv}}
\newglossaryentry{pmeasomr}{name=\ensuremath{\bbP_\omega},description={probability measure of the Markov chain starting at the root in the environment $\omega$},type=symbols,sort={pomega}}
\newglossaryentry{pmeasr}{name=\ensuremath{\bbP},description={probability measure of the mixed Markov chains starting at the root},type=symbols,sort={p}}
\newglossaryentry{ompmeas}{name=\ensuremath{\bP},description={probability measure describing the distribution of the environment $\omega$},type=symbols,sort={p}}

\newglossaryentry{edgeweight}{name=\ensuremath{w_n\left(\left\{u,v\right\}\right)},description={weight of the edge $\left\{u,v\right\}$ at time $n$},type=symbols,sort={wnuv}}
\newglossaryentry{reinforcement}{name=\ensuremath{\Delta},description={reinforcement parameter of the edge-reinforced random walk},type=symbols,sort={delta}}

\newglossaryentry{br}{name=\ensuremath{\br{T}},description={branching number of the tree $T$},type=symbols,sort={brt}}

\newglossaryentry{beta}{name=\ensuremath{B\left(\alpha,\beta\right)},description={Beta distribution with parameters $\alpha > 0, \beta > 0$ with density $\frac{\Gamma\left(\alpha + \beta\right)}{\Gamma\left(\alpha\right) \Gamma\left(\beta\right)} \cdot x^{\alpha-1} \cdot \left(1 - x\right)^{\beta-1}$ on $\left[0, 1\right]$},type=symbols,sort={balphabeta}}

\newglossaryentry{statmeas}{name=\ensuremath{\mu_\omega},description={stationary measure on the Markov chain given by the environment $\omega$},type=symbols,sort={muomega}}

\newglossaryentry{av}{name=\ensuremath{A_v},description={random variable defining the random environment $\omega$; $\omega\left(v, v_i\right) = \frac{A_{v_i}}{1 + \sum_{j = 1}^{b_v} A_{v_j}}$},type=symbols,sort={av}}
\newglossaryentry{cv}{name=\ensuremath{C_v},description={conductance along the edge from $v$ to its parent},type=symbols,sort={cv}}
\newglossaryentry{cutset}{name=\ensuremath{\Pi},description={a cutset: a finite set of nodes not including the root such that every infinite path starting at the root intersects $\Pi$ and there is no pair of nodes in $\Pi$ with one node being the ancestor of the other},type=symbols,sort={pi}}

\newglossaryentry{probchild}{name=\ensuremath{p_k},description={probability to have $k$ children for every node in a Galton-Watson tree},type=symbols,sort={pk}}
\newglossaryentry{gwmean}{name=\ensuremath{m},description={expected number of children for every node in a Galton-Watson tree},type=symbols,sort={m}}

\begin{document}

\begin{titlepage}
	\makeatletter
	
	\sffamily
	\def\svgwidth{60pt}
  \import{svg-inkscape/}{tum_mathematik_svg-tex.pdf_tex}\hfill\raisebox{5pt}[0pt][0pt]{\def\svgwidth{100pt}\import{svg-inkscape/}{TUM_svg-tex.pdf_tex}}\\[8em]
	{\large\merriweathersanslight Technical University of Munich\\[2em]
	Department of Mathematics\\[2em]
	Bachelor's Thesis in Mathematics}\\[8em]
	{\color{headingcolor}\Huge\bfseries Linearly Edge-Reinforced Random Walks}\\[4em]
	{\LARGE \@author}\\[10em]
	{\large\merriweathersanslight {\sffamily Supervisor:} Prof.~Dr.~rer.~nat.~habil.~Nina Gantert\\[1.5em]
	{\sffamily Submission Date:} \@date}

	\AddToShipoutPicture*{\put(0,0){%
		\parbox[b][\paperheight]{\paperwidth}{%
		\vfill
		\centering
		{\transparent{0.15}\includegraphics[width=\paperwidth]{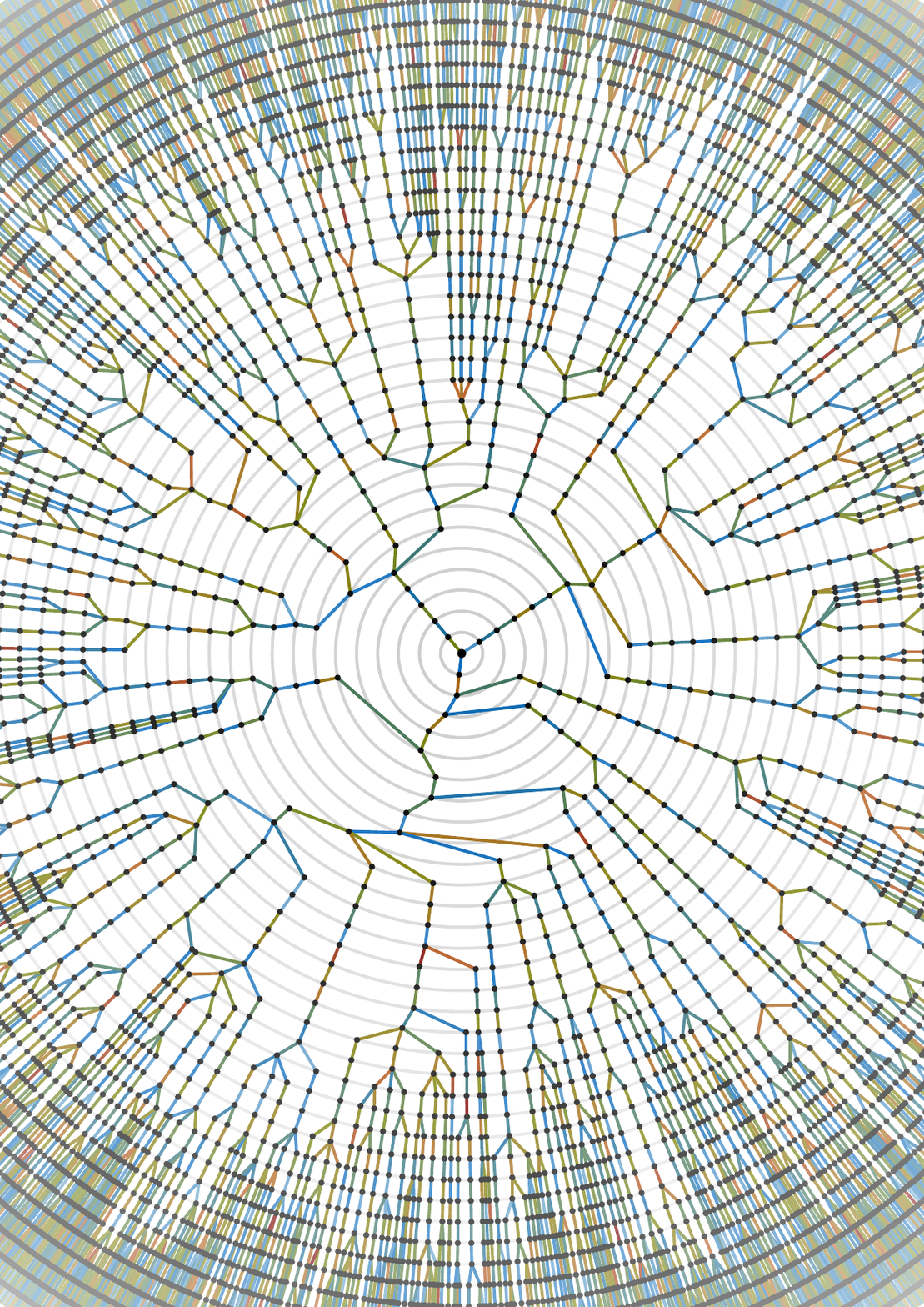}}%
		\vfill
	}}}
	\makeatother
\end{titlepage}

\ifthenelse{\boolean{print}}{\thispagestyle{empty} { ~ } \newpage}{}

\thispagestyle{empty}
{
	\makeatletter
	
	\sffamily
	\def\svgwidth{60pt}
  \import{svg-inkscape/}{tum_mathematik_svg-tex.pdf_tex}\hfill\raisebox{5pt}[0pt][0pt]{\def\svgwidth{100pt}\import{svg-inkscape/}{TUM_svg-tex.pdf_tex}}\\[6em]
	{\large\merriweathersanslight Technical University of Munich\\[1.5em]
	Department of Mathematics\\[1.5em]
	Bachelor's Thesis in Mathematics}\\[6em]
	{\color{headingcolor}\Huge\bfseries Linearly Edge-Reinforced Random Walks}\\[2em]
	{\color{headingcolor}\huge Linear selbstverst\"arkende Irrfahrten}\\[3em]
	{\LARGE \@author}\\[8em]
	{\large\merriweathersanslight {\sffamily Supervisor:} Prof.~Dr.~rer.~nat.~habil.~Nina Gantert\\[1.5em]
	{\sffamily Submission Date:} \@date}

	\makeatother
	\ifthenelse{\boolean{print}}{}{~\\[5em]\noindent {\large\merriweathersanslight final digitial version (page numbers do not agree with printed version)}}
}
\newpage

\ifthenelse{\boolean{print}}{\thispagestyle{empty} { ~ } \newpage}{}

\thispagestyle{empty}
{
	\makeatletter
	~\\[20em]
	\noindent I confirm that this bachelor's thesis is my own work and I have documented all sources and material used.\\[1em]
	\noindent Ich erkl\"are hiermit, dass ich diese Arbeit selbst\"andig und nur mit den angegebenen Hilfsmitteln angefertigt habe.
	
	\noindent ~\\[3em]
	\rule{8cm}{.4pt}\\
	\@author, M\"orlenbach, \@date
	\makeatother
}
\newpage
{
  \thispagestyle{fancy}
  ~\\[4em]
  \begin{center}
    \LARGE\bfseries\color{headingcolor}\sffamily Abstract
  \end{center}

  This thesis examines linearly edge-reinforced random walks on infinite trees. In particular,
  recurrence and transience of such random walks on general (fixed) trees as well as on Galton-Watson
  trees (i.e.~random trees) is characterized, and shown to be related to the branching number of these trees and a
  so-called reinforcement parameter. A phase transition from transience to recurrence takes place
  at a critical parameter value. As a tool, random walks in random environment are introduced and
  known results are repeated, together with detailed proofs. A result on quasi-independent
  percolation is proved as a by-product. Finally, for the edge-reinforced random walk on $\bbZ$,
  the existence of a kind of stationary / limiting distribution with finite moments is shown.

  ~\\[2em]
  \begin{center}
    \LARGE\bfseries\color{headingcolor}\sffamily Zusammenfassung
  \end{center}

  Diese Bachelorarbeit analysiert linear selbstverst\"arkende Irrfahrten auf unendlichen B\"aumen.
  Insbesondere wird die Rekurrenz und Transienz solcher Irrfahrten auf allgemeinen (festen) B\"aumen
  sowie auf Galton-Watson-B\"aumen (also zuf\"alligen B\"aumen) charakterisiert und es wird gezeigt,
  dass ein Zusammenhang zu der Verzweigungszahl dieser B\"aume und einem sogenannten Verst\"arkungs-Parameter
  existiert. F\"ur einen bestimmten Parameterwert findet ein Phasen\"ubergang von Transienz zu
  Rekurrenz statt. Als Werkzeug werden Irrfahrten in zuf\"alligen Umgebungen betrachtet und bereits
  bekannte Resultate wiederholt, komplett mit detaillierten Beweisen. Als Nebenprodukt wird eine
  Aussage zu quasi-unabh\"angiger Perkolation bewiesen. Zum Schluss wird die Existenz einer Art
  station\"arer / Grenz- Verteilung mit endlichen Momenten f\"ur die selbstverst\"arkende Irrfahrt
  auf $\bbZ$ gezeigt.
}
\newpage

{
\hypersetup{
  hidelinks
}
\tableofcontents
}

\newpage

{
\hypersetup{
  hidelinks
}
\listoffigures
\listoftables

~\\[2em]
\subsection*{A note on cross-references in this thesis}
Cross-references and citations are marked in \textcolor{linkred}{red}, and they are clickable and
directly link to the referenced object. Titles, definitions, theorems,
etc.~are colored in \textcolor{sectionblue}{blue}.

}

\newpage

\section{Introduction}
\label{sec:intro}

The central topic of this thesis is the linearly edge-reinforced random walk (ERRW), a special type of random
walk on a graph. The edges in the graph are weighted, and the probability to leave a node via one of
the incident edges is proportional to the respective edge weight compared to the weights of the other
incident edges. Initially, all edges have weight $1$, and each time an edge is crossed, its weight
is increased by a parameter $\Delta > 0$. In this thesis, this type of random walk is considered on
an infinite tree. Often, ERRW also refers to a more general concept of random walk
where different reinforcement schemes are possible, i.e.~the edge weight is not necessarily
increased by a constant value upon traversal but by values which change over time.

Intuitively, increasing $\Delta$ increases the probability that the random walk revisits parts of
the tree which it has visited before. The random walk becomes ``more'' recurrent. For small
values of $\Delta$, the ERRW is more similar to the simple random walk on the graph.
An interesting question is when the ERRW is recurrent (a.s.~returning to its
starting point) and when it is transient (never returning to its starting point with positive
probability). \autoref{sec:transrecerrw} will show that there is a phase transition at a critical
value $\Delta_0$: for $\Delta < \Delta_0$, the random walk is transient, and for $\Delta > \Delta_0$
it is recurrent.

Instead of fixing a tree on which the ERRW is run, it is also possible to consider
a random tree, and then the ERRW on this tree. For a special case of random trees,
Galton-Watson trees, there is a simple criterion for a.s.~recurrence and transience
of the ERRW on an instantiation of the random tree. This is formalized in
\autoref{thm:errwgw}.

In order to prove the central results on the ERRW, this thesis uses a detour via a
concept called random walk in random environment (RWRE). A RWRE is a Markov chain (MC)
where the transition probabilies between two states of the MC are random variables, defining the environment
(thus the name of the RWRE). First, the transition probabilies are chosen randomly, and
the random walk then runs on the resulting MC.

As the ERRW is not a Markov chain (MC), in contrast to many typical random walks,
it can be very hard to analyze. Pemantle showed in \cite{errwpemantle} that there is a connection
between the ERRW and RWRE. Indeed, the law of the ERRW is the same
as the law of a certain RWRE. Therefore, \autoref{sec:transrecrwre} first proves
results for RWRE, following \cite{rwrelyonspemantle}.

\subsection{Literature}
\label{ssec:lit}

The proofs are not self-contained but mostly rely on a single source only: \cite{probtreenet},
``Probability on Trees and Networks'' by Lyons and Peres. The book introduces electrical
network terminology which is helpful to analyze random walks. For this thesis, the relevant
parts are \cite{probtreenet}[Chapter 1, 2, 3, 5]. \cite{probtreenet} is an excellent book, and its
lecture is highly recommended. As a much shorter collection of notes, \cite{introrw} also gives a good
introduction of the key ideas used in \cite{probtreenet}, but without the electrical terminology
and by far not in the generality needed in this thesis.

The central part of this thesis involves \cite{rwrelyonspemantle}, which considers RWRE
on trees and also on Galton-Watson trees. The proofs originally relied on \cite{rwperctrees} for
random walks and percolation and \cite{ising} for quasi-independent percolation, and these sources
can still provide additional insights. As mentioned in the previous paragraph, this thesis
replaced these sources by the newer, all-encompassing \cite{probtreenet}, which was written by the
same author and uses a unified notation.

\cite{errwpemantle} is maybe the most influential work on this particular topic as it was the
first paper to demonstrate the connection between ERRW and RWRE. The proofs
in this paper were simplified and statements were generalized over time, which explains why this
is not the main source of this thesis. \cite{rrwrenlund} provides an overview over results on the
ERRW, and some parts are very close to the questions considered in this thesis. Instead of
Galton-Watson trees, \cite{rrwrenlund} considers regular trees but most calculations will turn
out to be the same for Galton-Watson trees. The main argument needed to extend the results to
Galton-Watson trees is actually due to \cite{rwrelyonspemantle} which considers RWRE
on Galton-Watson trees instead of just fixed trees as in \cite{rrwrenlund}.

Finally, further questions can be asked about RWRE and ERRW which go beyond
recurrence and transience. One such question is whether the speed of the random walk (defined as
the limit of the distance to the starting point at time $n$ divided by $n$) is positive in the
transient case. This question goes beyond the scope of this thesis, but \cite{transientaidekon}
is still an interesting source: it shows positive speed for regular trees, just as
\cite{errwpemantle} / \cite{rrwrenlund} considered recurrence and transience on regular trees.
Similarly to the results shown in this thesis, it could be possible to extend the positivity of
the speed to Galton-Watson trees.

\subsection{Main Results}
\label{ssec:results}

\autoref{thm:transrecrwre} and \autoref{thm:rwregw} characterize recurrence and transience of
RWRE for fixed trees and Galton-Watson trees, respectively. The results are due to
\cite{rwrelyonspemantle} and simply adopted (not extended in any way) in this thesis. They show
a phase transition between recurrence and transience which depends on the so-called branching
number of the tree which measures a kind of average degree of the nodes in the tree.
\autoref{lem:percolation} extends the results on percolation (random edge removal) of
\cite[Section 5.3]{probtreenet} to a case where the assumptions on the independence of the edge
removal is less strict.

\autoref{thm:errw} and \autoref{thm:errwgw} characterize recurrence and transience of the
ERRW for fixed trees and Galton-Watson trees, respectively. Again, there is a phase
transition which depends on the reinforcement parameter $\Delta$ as well as the branching number of the
tree and mean of the Galton-Watson tree, respectively. These results have not
been stated in the same generality before (as far as the author knows), but they follow from
\cite{rwrelyonspemantle} and the method from \cite{errwpemantle}.

\autoref{sec:errw_z} turns to a special case of the ERRW: the ERRW on $\bbZ$.
It is shown that a kind of stationary / limiting distribution exists, and \autoref{thm:zstat}
proves the new result that all moments of this distribution exist and are finite (again, the
result is new as far as the author knows).

\clearpage
\newpage

\section{Preliminaries}
\label{sec:prelim}

Throughout this thesis, random walks on infinite trees will be considered. All considered trees
will be locally finite, that is, the number of edges incident to a node is always finite.
The following notation is used in connection with trees:

\begin{definition}[Infinite Tree]
  \label{def:inftree}
  An \textbf{infinite (rooted) tree} $T$ is a tuple $\left(V, E\right)$ of vertices (or nodes) and
  edges, with $V$ an infinite set containing $\rho$ and
  $E \subseteq \left\{ \left\{u, v\right\} \midd| u,v \in V \land u \neq v \right\}$.
  As a tree, $T$ contains no cycles and is connected, and all nodes have finite degree.
  The following notation will be used:
  \begin{itemize}
    \item The vertex $\gls{root} \in V$ is called the \textbf{root} of the tree.
    \item $u, v \in V$ are \textbf{neighbors} if they are connected by an edge, written $\gls{neighbor}$.
    \item The \textbf{parent} of a vertex $v \in V$, denoted by $\gls{par}$, is the
    neighbor of $v$ which lies on the unique shortest path from $\rho$ to $v$. $\overleftarrow{\rho}$
    is undefined.
    \item The \textbf{children} of a vertex $v \in V$, a set denoted by $\gls{children}$, are
    all neighbors of $v$ which are not equal to $\overleftarrow{v}$. It is assumed that
    $C\left(v\right)$ is finite for every $v \in V$, and the children are denoted by
    $C\left(v\right) = \left\{\glslink{child}{v_1}, \ldots, v_{b_v}\right\}$ where $\gls{numch}$ is the number of children
    of vertex $v$.
    \item The \textbf{edge} leading from $\overleftarrow{v}$ to $v$ for a $v \in V$ is denoted by
    $\gls{paredge} := \left\{ \overleftarrow{v}, v \right\}$.
    \item The \textbf{distance} $\gls{nodedist}$ of a vertex $v \in V$ is the (edge) length of the
    unique shortest path from $\rho$ to $v$. If $E \ni e = e_v$, then set
    $\gls{edgedist} = \left|v\right|$. Let $\gls{distset} := \left\{v \in V : \left|v\right| = n\right\}$ be
    the set of vertices at distance $n$ from the root.
    \item $u \in V$ is an \textbf{ancestor} of $v \in V$, written $\gls{ancestor}$, if $u \neq v$ and $u$ lies on the
    unique shortest path from $\rho$ to $v$. If $u < v$ or $u = v$, the notation $u \leq v$ is used.
  \end{itemize}
\end{definition}

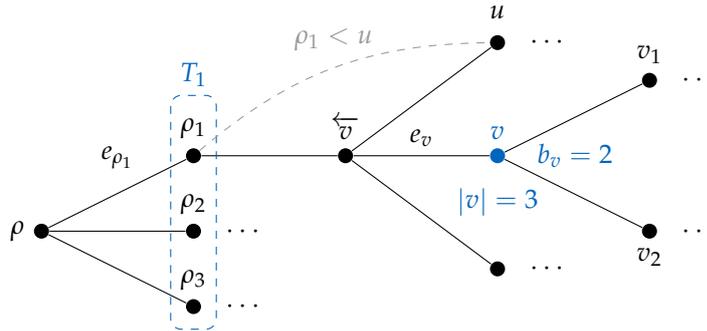
\begin{figure}[H]
  \begin{center}
    \input{figures/tree_notation.tex}
    \caption{Notation for trees}
    \label{fig:treenotation}
  \end{center}
\end{figure}

\subsection{Random Walks}
\label{ssec:rw}

One type of random walk on an infinite tree is the RWRE, which is used as a tool to
analyze the central topic of this thesis, the ERRW. The RWRE is basically
a Markov chain, but with random transition probabilities, which are also referred to as the
environment. The formal definition is as follows:

\begin{definition}[Random Walk in Random Environment]
  \label{def:rwre}
  Consider a tree $T = \left(V, E\right)$. Let $\left(\gls{av}\right)_{v \in V \setminus \left\{\gls{root}\right\}}$
  be a collection of random variables over some probability space with $A_v \geq 0$. The \textbf{environment}
  $\omega : V \times V \to \left[0, 1\right]$ on the tree $T$ is defined by
  \begin{align*}
    \omega\left(v, \gls{child}\right) := \frac{A_{v_i}}{\mathbbm{1}_{v \neq \rho} + \sum_{j = 1}^{b_v} A_{v_j}} \textrm{ for } 1 \leq i \leq \gls{numch}
    \qquad
    \omega\left(v, \gls{par}\right) := \frac{1}{1 + \sum_{j = 1}^{b_v} A_{v_j}} \;\; \left(v \neq \rho\right)
  \end{align*}
  The Markov chain (MC) $\left(X_n\right)_{n \geq 0}$, starting at $v \in V$, is now
  defined by
  \begin{align*}
    \glslink{pmeasom}{\bbP_\omega^v}\left[X_0 = v\right] = 1
    \qquad
    \Prbcc{\omega}{v}{X_{n+1} = w \midd| X_n = u} = \glslink{omuv}{\omega\left(u,w\right)}
  \end{align*}
  $X_n$ is a random walk taking values in $V$. If $X_n$ starts at the root $\rho$,
  $\glslink{pmeasomr}{\bbP_\omega}$ will be used instead of $\bbP_{\omega}^{\rho}$.

  The distribution of $\omega$ is called $\gls{ompmeas}$ and the probability measures $\glslink{pmeas}{\bbP^v}$, $\glslink{pmeasr}{\bbP}$ are defined
  by setting
  \begin{align*}
    \Prbcs{v}{\cdot} := \int \Prbcc{\omega}{v}{\cdot} \mudx{\bP}{\omega}
  \end{align*}
  $\left(X_n\right)_{n \geq 0}$, distributed according to $\bbP^v$, is called \textbf{random walk in random environment}
  (RWRE), or mixture of MCs.
\end{definition}

In the following, it will always be assumed that $A_v, v \in V \setminus \left(\left\{\rho\right\} \cup \glslink{distset}{T_1}\right)$
are identically distributed, and $A$ will denote a random variable with the same distribution. Furthermore,
$A_v$ and $A_u$ are assumed to be independent for $u, v \in V \setminus \left\{\rho\right\}$ with
$\overleftarrow{u} \neq \overleftarrow{v}$.

Note that the $A_v$ can be reconstructed from the environment $\omega$. For $v \neq \rho$, it holds that
$A_{v_i} = \omega\left(v, v_i\right) \cdot \omega\left(v, \overleftarrow{v}\right)^{-1}$. Hence, instead of the $A_v$,
$\omega$ is sometimes specified directly. If that is the case, the random variables
$\omega\left(v, v_i\right) \cdot \omega\left(v, \overleftarrow{v}\right)^{-1}$ are identically distributed
for $v \neq \rho$. Note that setting $A_{\rho_i} = \omega\left(\rho, \rho_i\right)$ will complete the
definition of the random variables $A_v$. In particular, the random variables $A_{\rho_i}$ can always
be chosen such that $0 \leq A_{\rho_i} \leq 1$ (this is not necessarily true for the $A_v$ in general).

Following \cite[Chapter 2]{probtreenet}, the following terms will be introduced for a
MC given by an environment $\omega$ (which, in turn, is defined via the random variables
$A_v$) and for trees in general:
\begin{enumerate}[(1)]
  \item \label{term:cutset} A \textbf{cutset} $\gls{cutset} \subseteq V$ is a finite set of nodes not
  including the root $\rho$ such that every infinite path starting at $\rho$ intersects $\Pi$ and
  there is no pair $u, v \in \Pi$ with $\gls{ancestor}$. Such a set cuts the connection between the
  root and $\infty$.
  \item \label{term:conductance} The \textbf{conductance} $\gls{cv}$ along the edge $\gls{paredge}$ for some
  $\rho \neq v \in V$ is defined as $C_v := \prod_{\rho \neq u \leq v} A_u$. Given these
  conductances, the tree $T$ can be seen as an electrical network. $C_v$ is a random variable
  depending on the value of the $A_v$, i.e.~depending on the environment $\omega$.
  \item \label{term:statmeas} $\gls{statmeas}: V \to \left[0, \infty\right)$ is called a \textbf{stationary
  measure} for the MC given by the environment $\omega$ if it satisfies
  $\mu_\omega\left(v\right) = \sum_{\gls{neighbor}} \mu_\omega\left(u\right) \omega\left(u, v\right)$
  for all $v \in V$.
  $\mu_\omega$ is further called a stationary distribution if it is a probability measure.
  
  Setting $\mu_\omega\left(v\right) := C_v + \sum_{i=1}^{b_v} C_{v_i}$
  (where $C_\rho$ should be $0$) defines a stationary measure. It holds that
  \begin{align*}
    \omega\left(v, v_i\right) = \frac{C_{v_i}}{\mu_\omega\left(v\right)}
    \qquad
    \omega\left(v, \overleftarrow{v}\right) = \frac{C_v}{\mu_\omega\left(v\right)} \;\; \left(v \neq \rho\right)
  \end{align*}
  \item \label{term:rectrans} The MC given by $\omega$ is called
  \textbf{recurrent} if $\Prbc{\omega}{X_n = \rho \textrm{ for infinitely many } n} = 1$, and
  \textbf{transient} otherwise. In the transient case, it holds that
  $\Prbc{\omega}{X_n = \rho \textrm{ for infinitely many } n} = 0$, so the root is only visited
  finitely often a.s., and $\Prbc{\omega}{\exists n > 0: X_n = \rho} < 1$. In the recurrent
  case, it hols that $\Prbc{\omega}{\exists n > 0: X_n = \rho} = 1$ (obviously). Since the considered
  MCs are a.s.~irreducible (every node can be reached from every other node
  with positive probability), this corresponds to the standard notion of recurrence and transience.

  The MC is further called \textbf{positive recurrent} if the expected return time to the
  root is finite, i.e.~$\bbE_\omega\left[\min\left\{ n > 0 : X_n = \rho \right\}\right] < \infty$ (this implies recurrence).
  Positive recurrence is equivalent to the existence of a (unique) stationary distribution. Note that the
  existence of a stationary measure does not imply the existence of a stationary distribution since
  the stationary measure may be infinite (if it is finite, it can always be normalized). Indeed,
  in the null recurrent case (recurrent, but not positive recurrent), there exists a unique
  stationary measure which is infinite.
\end{enumerate}

For more background on electrical networks (and conductances in particular) and their relation to
random walks, \cite{probtreenet} gives an excellent introduction (in Chapters 1 and 2).

\begin{figure}[H]
  \begin{center}
    \input{figures/rwre_instantiation.tex}
    \caption{An instantiation of the environment of a RWRE and related terms}
    \label{fig:rwre}
    Values of the random variables $A_v$ in \textcolor{tumOrange}{orange}, conductances in \textcolor{tumBlue}{blue}.
  \end{center}
\end{figure}
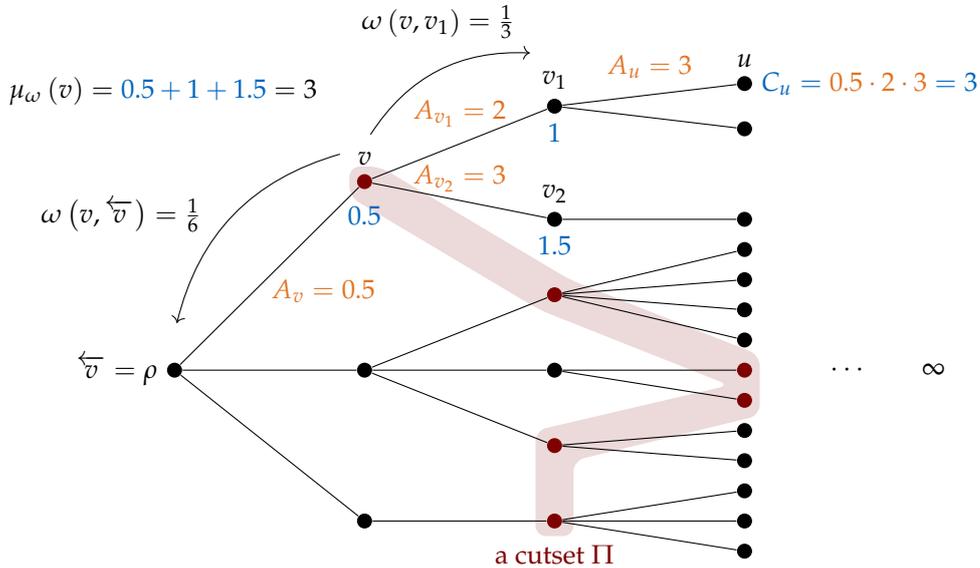

The actual, central topic of this thesis is another type of random walk, the ERRW. The
ERRW is not a MC because the transition probabilities change over time: the
probability to take an edge in the tree is roughly proportional to the number of traversals the
edge has seen so far.

\begin{definition}[Linearly Edge-Reinforced Random Walk]
  \label{def:errw}
  Consider a tree $T = \left(V, E\right)$ and weight functions $w_n : E \to \left[0, \infty\right)$
  for $n \geq 0$ with $w_0 \equiv 1$. Further, let $\gls{reinforcement} \geq 0$ be a \textbf{reinforcement parameter}.
  The \textbf{linearly edge-reinforced random walk} (ERRW) $\left(X_n\right)_{n \geq 0}$ with parameter $\Delta$ on $T$
  is defined as follows:
  \begin{align*}
    \textrm{for } n \geq 1: \textrm{ }
    \gls{edgeweight} = 1 + \Delta \sum_{i=1}^n \mathbbm{1}_{\left(X_{i-1} = u \land X_i = v\right) \lor \left(X_{i-1} = v \land X_i = u\right)}
  \end{align*}
  i.e.~the weight of each edge is increased by $\Delta$ upon every traversal, and
  \begin{align*}
    \Prb{X_0 = \gls{root}} = 1
    \qquad
    \Prb{X_{n+1} = v \midd| X_n = u} = \frac{w_n\left( \left\{u, v\right\} \right)}{w_n\left(\glslink{paredge}{e_u}\right) + \sum_{j=1}^{\glslink{numch}{b_u}} w_n\left(e_{\glslink{child}{u_j}}\right)}
  \end{align*}
  i.e.~the probability to transition from $u$ to $v$ at time $n$ is proportional to the weight $w_n$ associated
  to the edge $\left\{u, v\right\}$, compared to the sum of the weights $w_n$ of all incident edges.
\end{definition}

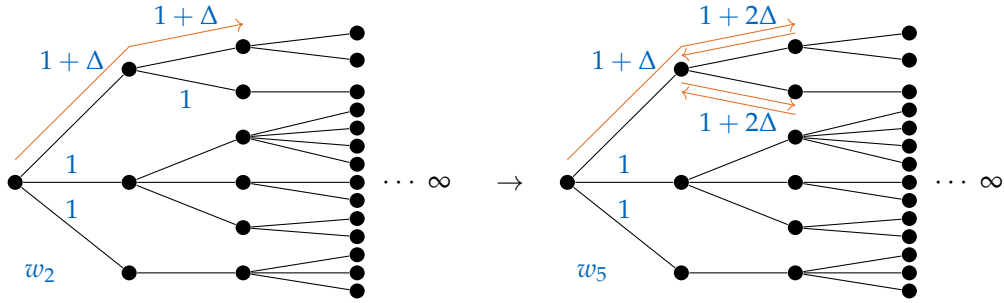
\begin{figure}[H]
  \begin{center}
    \input{figures/errw_example.tex}
    \caption{A path taken by the ERRW and the corresponding change in edge weights}
    \label{fig:errw}
  \end{center}
\end{figure}

\autoref{fig:errw_on_gw} on page \pageref{fig:errw_on_gw} shows an example for the ERRW
on a Galton-Watson tree (which will be defined below). The example shows that with a low
reinforcement parameter $\Delta$, the visible part of the tree is exited relatively fast, while the
random walk stays close to the root much longer for higher values of $\Delta$.

\subsection{Galton-Watson Trees}
\label{ssec:gwtree}

A random walk can also take place on a random tree. This thesis looks at Galton-Watson trees, which
are defined as follows:

\begin{definition}[Galton-Watson Tree]
  \label{def:gwtree}
  A \textbf{Galton-Watson tree} is constructed as follows: starting with the root $\rho$ and a
  probability distribution over the non-negative integers, denoted by $\gls{probchild}, k \in \bbN_{\geq 0}$
  with $\sum_{k=0}^\infty p_k = 1$, the number of children of each node is chosen independently
  according to the probability disitribution $p_k, k \in \bbN_{\geq 0}$. In terms of the notation
  of \autoref{def:inftree}, $\Prb{\gls{numch} = k} = p_k$ for all nodes $v$, and the events
  $\left\{b_v = k\right\}$ are independent for all nodes $v$.

  The event that the constructed tree is finite is called \textbf{extinction}.

  The expected number of children is denoted by $m := \sum_{k=0}^\infty k \cdot p_k$, also called
  \textbf{mean} of the Galton-Watson tree for short.
\end{definition}

\begin{figure}
  \begin{center}
    \def\svgwidth{11cm}
    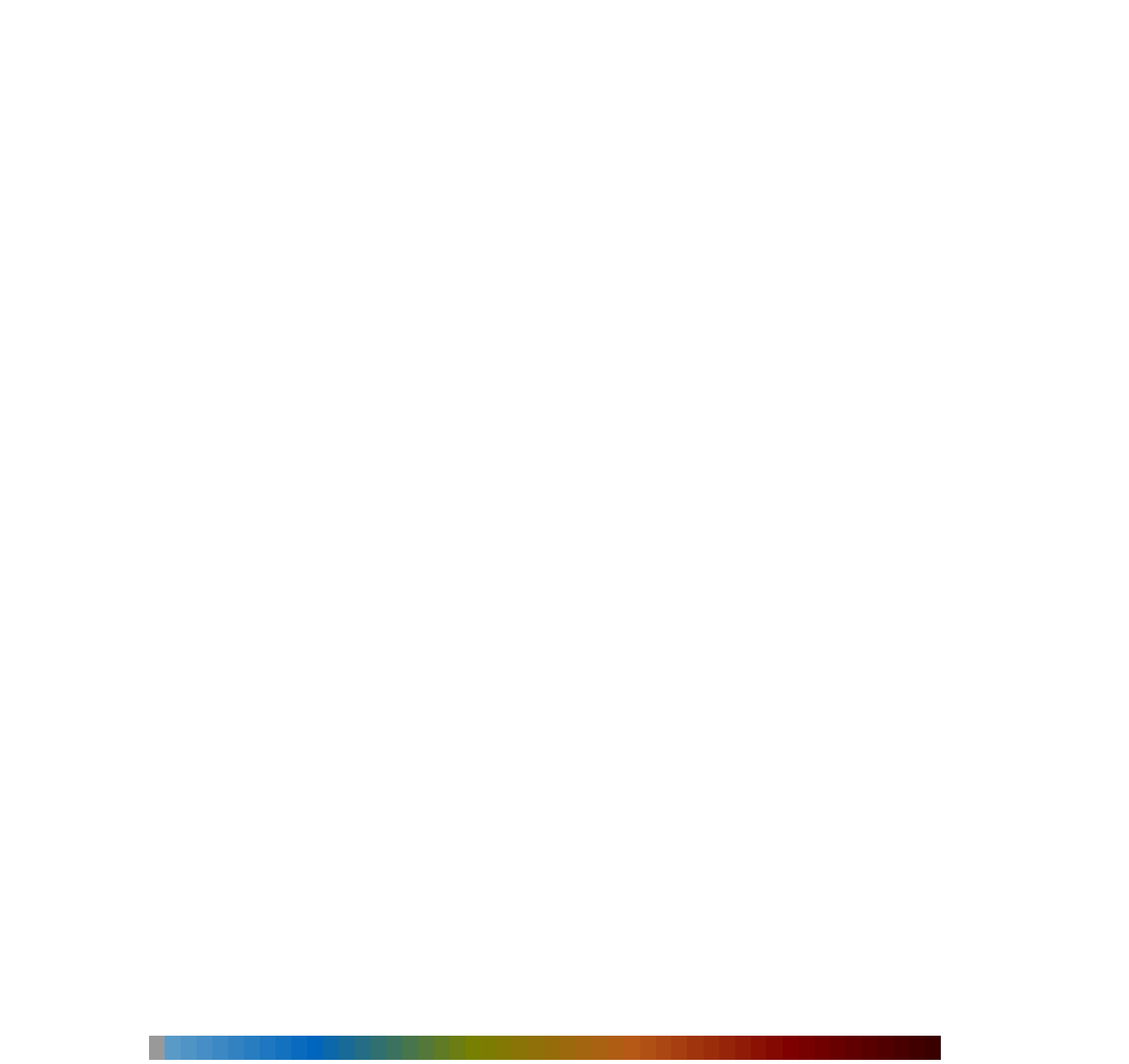 \\[2em]
    \def\svgwidth{11cm}
    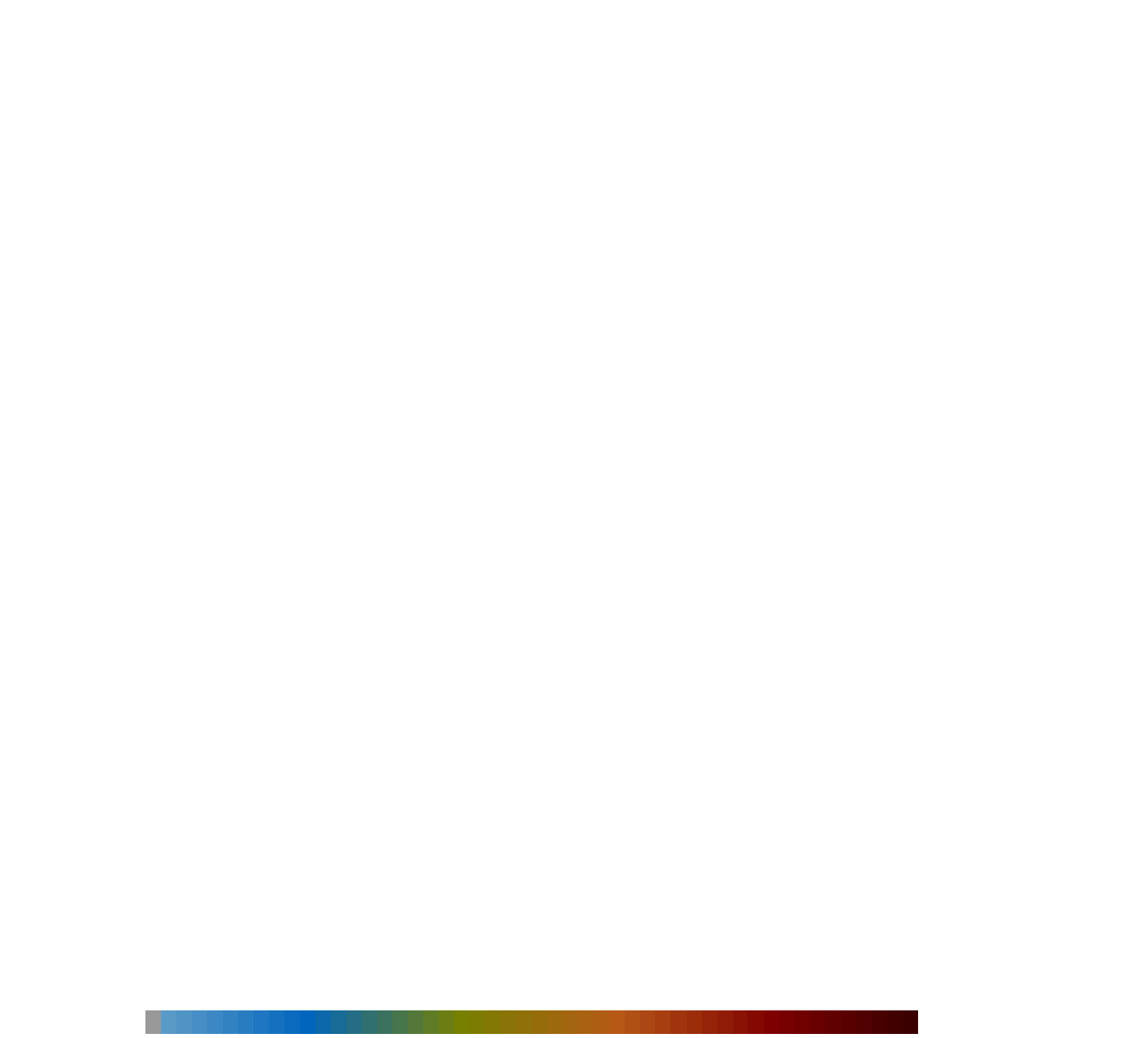
    \caption{ERRW on a Galton-Watson tree}
    \parbox{\textwidth}{Both trees are Galton-Watson trees with $p_1 = 0.85, p_2 = 0.1, p_3 = 0.05$.
    The ERRW was simulated on the trees for $100000$ steps or until the boundary of the
    visible part of the trees was reached. The weights $w_n$ of the edges at the end of the simulations
    are color coded.}
    \label{fig:errw_on_gw}
  \end{center}
\end{figure}

The nodes of a Galton-Watson tree could formally be defined as a subset of $\bbN_{\geq 0} \times \bbN_{\geq 0}$
where $\left(0, 0\right) = \gls{root}$, and a collection of iid random variables
$L_i^{\left(n\right)}, \left(i, n\right) \in \bbN_{\geq 0} \times \bbN_{\geq 0}$ taking values in
the non-negative integers, distributed as specified by $p_k, k \in \bbN_{\geq 0}$. If a node
$\left(n, i\right)$ is contained in the tree (the $i$-th node, starting with $0$, at distance $n$
from the root), then $\left(n, 0\right), \ldots, \left(n, i-1\right)$ also belong to the tree.
The number of children of $\left(n, i\right)$, i.e.~$\glslink{numch}{b_{\left(n, i\right)}}$,
is chosen according to $L_i^{\left(n\right)}$, and the children of $\left(n, i\right)$ are given by
\begin{align*}
  \glslink{children}{C\left(n, i\right)} := \left\{ \left(n+1, a\right), \ldots, \left(n+1, a+b_{\left(n, i\right)}-1\right) \right\}
  \quad \textrm{where} \quad
  a = \sum_{j=0}^{i-1} b_{\left(n, j\right)}
\end{align*}
\autoref{fig:errw_on_gw} on page \pageref{fig:errw_on_gw} shows two instantiations of Galton-Watson
trees where the number of children is identically distributed.

\subsection{Flows and Branching Number}
\label{ssec:flowbr}

In order to argue about transience and recurrence of random walks, it is useful to consider flows
on trees. Many results from \cite{probtreenet} are closely related to flows on trees and a concept
called branching number, which will be introduced below, and which is also defined on the basis of
flows.

\begin{definition}[Flow]
  \label{def:flow}
  Let $T = \left(V, E\right)$. A \textbf{flow} $f : E \to \bbR$ on $T$ from a set of source nodes
  $A \subseteq V$ to a set of sink nodes $B \subseteq V$ satisfies the following properties:
  \begin{enumerate}[(i)]
    \item \label{def:flowcons} $\forall v \in V \setminus \left(A \cup B\right): f\left(\gls{paredge}\right) = \sum_{i = 1}^{\gls{numch}} f\left(e_{\gls{child}}\right)$
    \item \label{def:flowsrc} $\forall v \in A: f\left(e_v\right) \leq \sum_{i = 1}^{b_v} f\left(e_{v_i}\right)$
    \item \label{def:flowsink} $\forall v \in B: f\left(e_v\right) \geq \sum_{i = 1}^{b_v} f\left(e_{v_i}\right)$
  \end{enumerate}
  Since there is no parent of the root $\rho$, $e_{\rho}$ is undefined, and $f\left(e_{\rho}\right)$
  should therefore be replaced by $0$ in the equations above.

  For $v \in V$, the flow $f\left(e_v\right)$ along the edge $e_v = \left\{ \gls{par}, v \right\}$
  is interpreted as a flow of $f\left(e_v\right)$ units from $\overleftarrow{v}$ to $v$ if
  $f\left(e_v\right) \geq 0$, otherwise as a flow of $\left|f\left(e_v\right)\right|$ units from
  $v$ to $\overleftarrow{v}$.

  The \textbf{strength} of a flow is defined as
  $\Strength{f} := \sum_{v \in A} \left(\sum_{i = 1}^{b_v} f\left(e_{v_i}\right)\right) - f\left(e_v\right)$.
  A unit flow is a flow with strength $1$.

  A flow from the root $\rho$ to infinity (short: \textbf{flow from $\rho$ to $\infty$}) is a flow
  $f$ on $T$ where $A = \left\{\rho\right\}$ and $B = \varnothing$.
\end{definition}

\ref{def:flowcons} is the flow conservation property: the incoming flow from the parent node
equals the sum of the outgoing flows to the children (if $f$ takes negative values, there can also
be an outgoing flow to the parent or an incoming flow from a child). At the source nodes, the
outgoing flow can be larger than the incoming flow, as specified by \ref{def:flowsrc}, while
the reverse is true for the sink nodes, as defined in \ref{def:flowsink}.

The concept of a branching number for infinite trees will be used to characterize transience
and recurrence of random walks. The branching number measures the average number of children per
node in the tree. It is intuitively clear that a higher number of children will make the random
walk more transient, explaining why the branching number is such a useful tool in this context.
A related concept, the growth rate, defined as $\lim_{n\to\infty} \left|\gls{distset}\right|^{\frac{1}{n}}$
(if the limit exists) turns out to be too imprecise to give accurate statements on transience and
recurrence. The branching number, in contrast, takes additional structural properties of the tree
into account, more than just the average growth per generation.

\begin{definition}[Branching Number]
  \label{def:br}
  Let $T = \left(V, E\right)$ be an infinite tree with root $\rho \in V$. Then, the \textbf{branching number}
  of $T$ is given by
  \begin{align*}
    \gls{br} := \sup \left\{ \lambda > 0 : \exists \textrm{ flow } f
    \textrm{ from } \rho \textrm{ to } \infty \textrm{ with } f \not\equiv 0
    \textrm{ such that } \forall e \in E: \left|f\left(e\right)\right| \leq \lambda^{-\gls{edgedist}} \right\}
  \end{align*}
  Since $f \not\equiv 0$ is required, the following textual definition can be given:
  the branching number is the supremum over all $\lambda$ such that restricting the flow on
  edges at distance $n$ from the root to at most $\lambda^{-n}$ still admits a non-zero flow from
  the root to infinity.
\end{definition}

Note that $\br{T} \geq 1$, since a flow of $1$ along an infinite path fulfils the constraint
$\left|f\left(e\right)\right| \leq \lambda^{-\left|e\right|}$ with $\lambda = 1$. Since $T$ is
always an infinite tree in this thesis, $T$ contains such an infinite path (because the number
of vertices is infinite and every tree is connected). In contrast, the branching number of a finite
tree is considered to be zero. Examples for the branching number, and its relation to the growth
rate can be found in \cite[Section 1.2]{probtreenet}.

Flows and cuts in a graph are closely related: the max-flow min-cut theorem states that the strength
of a maximal flow restricted by edge capacities in a finite graph between two sets of vertices
equals the capacity of the minimal cut -- details can be found in \cite[page 75]{probtreenet}.
This theorem can be extended to infinite graphs, as shown by \cite[Theorem 3.1]{probtreenet}. For
the infinite case, the strength of the maximal flow equals the infimum over all capacities of
cutsets. Whenever the max-flow min-cut theorem is used in this thesis, it will be used in the
generalized version. The max-flow min-cut theorem gives an alternative expression for the branching
number:
\begin{align}
  \label{equ:br_alt}
  \br{T} = \sup \left\{ \lambda > 0 : \inf_{\gls{cutset} \textrm{ cutset}} \sum_{v \in \Pi} \lambda^{-\gls{nodedist}} > 0 \right\}
  = \inf \left\{ \lambda > 0 : \inf_{\Pi \textrm{ cutset}} \sum_{v \in \Pi} \lambda^{-\left|v\right|} = 0 \right\}
\end{align}

Consider the tree $\gls{ntree}$ which contains all nodes $v \in V$ with $n \divides \left|v\right|$
($n$ divides $\left|v\right|$) where two nodes $u, v$ are connected by an edge if, and only if, $u \leq v$ and
$\left|u\right| + n = \left|v\right|$ (where all distances to the root are measured in the original
tree $T$). Later, it will be necessary to calculate the branching number of this tree, which is
directly related to the branching number of the original tree $T$.

\begin{lemma}[Branching Number of $T^{\left(n\right)}$]
  \label{lem:tnbr}
  Let $T = \left(V, E\right)$ be an infinite tree and $T^{\left(n\right)}$ as
  defined above. Then $\br{T^{\left(n\right)}} = \br{T}^n$.
\end{lemma}

\begin{tproof}
  Set $T^{\left(n\right)} = \left(V^{\left(n\right)}, E^{\left(n\right)}\right)$.
  In contrast to the definition of $T^{\left(n\right)}$, $\left|e\right|$ and $\left|v\right|$
  for an edge $e$ or a node $v$ are measured in the tree they belong to, so if $v \in V^{\left(n\right)}$
  with $\left|v\right| = k$ is considered, the same $v$, seen as a node $v \in V$ of the tree $T$
  will have distance $\left|v\right| = kn$ to the root in $T$. It will be implicitly clear if a
  node is considered as a node of $T$ or $T^{\left(n\right)}$.

  Let $f$ be a non-zero flow in $T$ from $\rho$ to $\infty$ which satisfies $\left|f\left(e\right)\right| \leq \lambda^{-\left|e\right|}$.
  Then it is easy to construct a flow $g$ in $T^{\left(n\right)}$ which satisfies
  $\left|g\left(e\right)\right| \leq \lambda^{-n \left|e\right|}$ by simply setting the flow
  on an edge in $T^{\left(n\right)}$ which ends in the node $v$ to $f\left(\gls{paredge}\right)$ where
  $e_v \in E$ is an edge in $T$.

  The reverse direction is equally easy to proof: if $g$ is a non-zero flow in $T^{\left(n\right)}$
  satisfying $\left|g\left(e\right)\right| \leq \left(\lambda^n\right)^{-\left|e\right|}$,
  define the flow $f$ in $T$ by setting $f\left(e_v\right)$ (with $v \in V$ such that $n \divides \left|v\right|$
  and $e_v \in E$) to $g\left(e_v\right) \cdot \lambda^{-n+1}$ (where $v$ is seen as $v \in V^{\left(n\right)}$ and
  $e_v \in E^{\left(n\right)}$ is a different edge than before). Then, there is a unique choice
  for $f$ on the remaining edges. For this direction, the strength of $f$ is different from the
  strength of $g$, but $f$ is still a non-zero flow satisfying $\left|f\left(e\right)\right| \leq \lambda^{-\left|e\right|}$.
\end{tproof}

Finally, for the case of the random Galton-Watson tree, the branching number is simply the expected
number of children per node:

\begin{lemma}[Branching Number of a Galton-Watson Tree]
  \label{lem:gwbr}
  Let $T$ be a Galton-Watson tree with mean $\gls{gwmean} > 1$. Then, conditional on non-extinction of $T$,
  $\br{T} = m$ a.s.
\end{lemma}

This lemma is proved in \cite[Corollary 5.10]{probtreenet}.
\clearpage
\newpage

\section{Transience and Recurrence of RWRE}
\label{sec:transrecrwre}

This section is devoted to the characterization of recurrent and transient RWREs. In
order to prove the main result, \autoref{thm:transrecrwre}, three lemmata are required when following
the proof from \cite{rwrelyonspemantle}; two will be proved in advance, one was
extracted from the proof and will be proved later.

\begin{lemma}[Chernoff-Cram\'er Theorem \textnormal{(\cite[page 129]{rwrelyonspemantle})}]
  \label{lem:cct}
  There are two versions of this lemma, an additive and a multiplicative one:
  \begin{enumerate}[(i)]
    \item \label{lem:cct_add} Let $X$ be a real-valued random variable, and set
    $\gamma\left(a\right) := \inf_{t \geq 0} \left(-at + \log \Ex{e^{tX}}\right)$. If $S_n$ denotes the sum
    of $n$ independent random variables with the same distribution as $X$, then, $\forall a \in \bbR$:
    \begin{align}
      \label{equ:cct_add}
      \frac{1}{n} \cdot \log \Prb{S_n \geq n \cdot a} \overset{n \to \infty}{\longrightarrow} \gamma\left(a\right)
      \quad \textrm{ and } \quad \frac{1}{n} \cdot \log \Prb{S_n \geq n \cdot a} \leq \gamma\left(a\right)
    \end{align}
    \item \label{lem:cct_mul} Let $A$ be a real-valued random variable with $0 < A < \infty$ a.s.,
    and let $U_n$ denote the product of $n$ independent random variables with the same distribution
    as $A$. Then, $\forall s > 0$:
    \begin{align*}
      \Prb{U_n \geq s^n}^{\frac{1}{n}} \overset{n \to \infty}{\longrightarrow} \inf_{t \geq 0} \left(s^{-t} \cdot \Ex{A^t}\right)
      \quad \textrm{ and } \quad \Prb{U_n \geq s^n}^{\frac{1}{n}} \leq \inf_{t \geq 0} \left(s^{-t} \cdot \Ex{A^t}\right)
    \end{align*}
  \end{enumerate}
\end{lemma}

\begin{tproof}
  Since \ref{lem:cct_mul} is a corollary of \ref{lem:cct_add}, first consider \ref{lem:cct_add}.
  Denote by $Y_1, Y_2, \ldots$ independent random variables with the same distribution as $X$, and
  set $S_n = \sum_{i = 1}^n Y_i$. Consider for any $t > 0$
  \begin{align*}
    \Prb{S_n \geq n \cdot a} &= \Prb{e^{tS_n} \geq e^{t \cdot n \cdot a}} \overset{\textrm{Markov inequ.}}{\leq} \Ex{e^{tS_n}} e^{-t \cdot n \cdot a} \\
    \implies \log \Prb{S_n \geq n \cdot a} &\leq -t \cdot n \cdot a + \log \Ex{e^{tS_n}} \\
    \implies \frac{1}{n} \cdot \log \Prb{S_n \geq n \cdot a} &\leq -at + \log \sqrt[n]{\prod_{i=1}^n \Ex{e^{tY_i}}} = -at + \log \Ex{e^{tX}}
  \end{align*}
  Note that the inequality also holds for $\Ex{e^{tX}} = \infty$, and trivially for $t = 0$. This shows that
  \begin{align*}
    \frac{1}{n} \cdot \log \Prb{S_n \geq n \cdot a} \leq \gamma\left(a\right)
  \end{align*}
  The convergence property is harder to prove (and actually the property which will be used later).
  Refer to \cite{cramerthm} for a proof.

  Now, to proof \ref{lem:cct_mul}, set $X := \log A$. With $Y_1, Y_2, \ldots$ iid with
  the same law as $X$, as above, $Z_1, Z_2, \ldots$ iid with the same law as $A$,
  $S_n = \sum_{i = 1}^n Y_i$, and $U_n = \prod_{i=1}^n Z_i$, it holds that $\log U_n = S_n$
  (in distribution). In addition, setting $s = e^a$,
  \begin{align*}
    \frac{1}{n} \cdot \log \Prb{S_n \geq n \cdot a} = \frac{1}{n} \cdot \log \Prb{\log U_n \geq n \cdot a}
    = \frac{1}{n} \cdot \log \Prb{U_n \geq \left(e^a\right)^n} = \frac{1}{n} \cdot \log \Prb{U_n \geq s^n}
  \end{align*}
  On the other hand,
  \begin{align*}
    \gamma\left(a\right) = \gamma\left(\log s\right) = \inf_{t \geq 0} \left(-t \log s + \log \Ex{e^{t\log A}}\right) =
    \inf_{t \geq 0} \left(\log s^{-t} + \log \Ex{A^t}\right)
  \end{align*}
  Now, applying \ref{lem:cct_add} and taking $\exp\left(\cdot\right)$ on both sides of \Cref{equ:cct_add},
  the following statement can be derived for $U_n$:
  \begin{align*}
    \Prb{U_n \geq s^n}^{\frac{1}{n}} \overset{n \to \infty}{\longrightarrow} \inf_{t \geq 0} \left(s^{-t} \cdot \Ex{A^t}\right)
    \quad \textrm{ and } \quad \Prb{U_n \geq s^n}^{\frac{1}{n}} \leq \inf_{t \geq 0} \left(s^{-t} \cdot \Ex{A^t}\right) & \qedhere
  \end{align*}
\end{tproof}

\begin{figure}
  \begin{center}
    \begin{tabular}{cc}
      ~~~~~~~~~~~$\Delta = 1$ & ~~~~~~~~~~~~$\Delta = 0.2$ \\
      \input{simulations/plot_exat-d1.00.tex} & \input{simulations/plot_exat-d0.20.tex} \\
      \input{simulations/plot_exatmod-d1.00.tex} & \input{simulations/plot_exatmod-d0.20.tex} \\[0.7em]
      \input{simulations/plot_maxs-d1.00.tex} & \input{simulations/plot_maxs-d0.20.tex}
    \end{tabular}
    \caption{Alternative expression for $p$}
    \parbox{\textwidth}{The expressions from \autoref{lem:altp} for the random variable $A$ which
    is later relevant for the ERRW (see \autoref{sec:errw_z}). This $A$ depends on the
    reinforcement parameter $\Delta$ of the ERRW, explaining the different plots for the
    different values of $\Delta$. In the middle row, the different curves show the values for
    different $s$, from $s = 1$ (\textcolor{tumBlue}{dark blue}) to $s = 0.1$ (\textcolor{tumBlue!40}{light blue}).}
    \label{fig:altp}
  \end{center}
\end{figure}
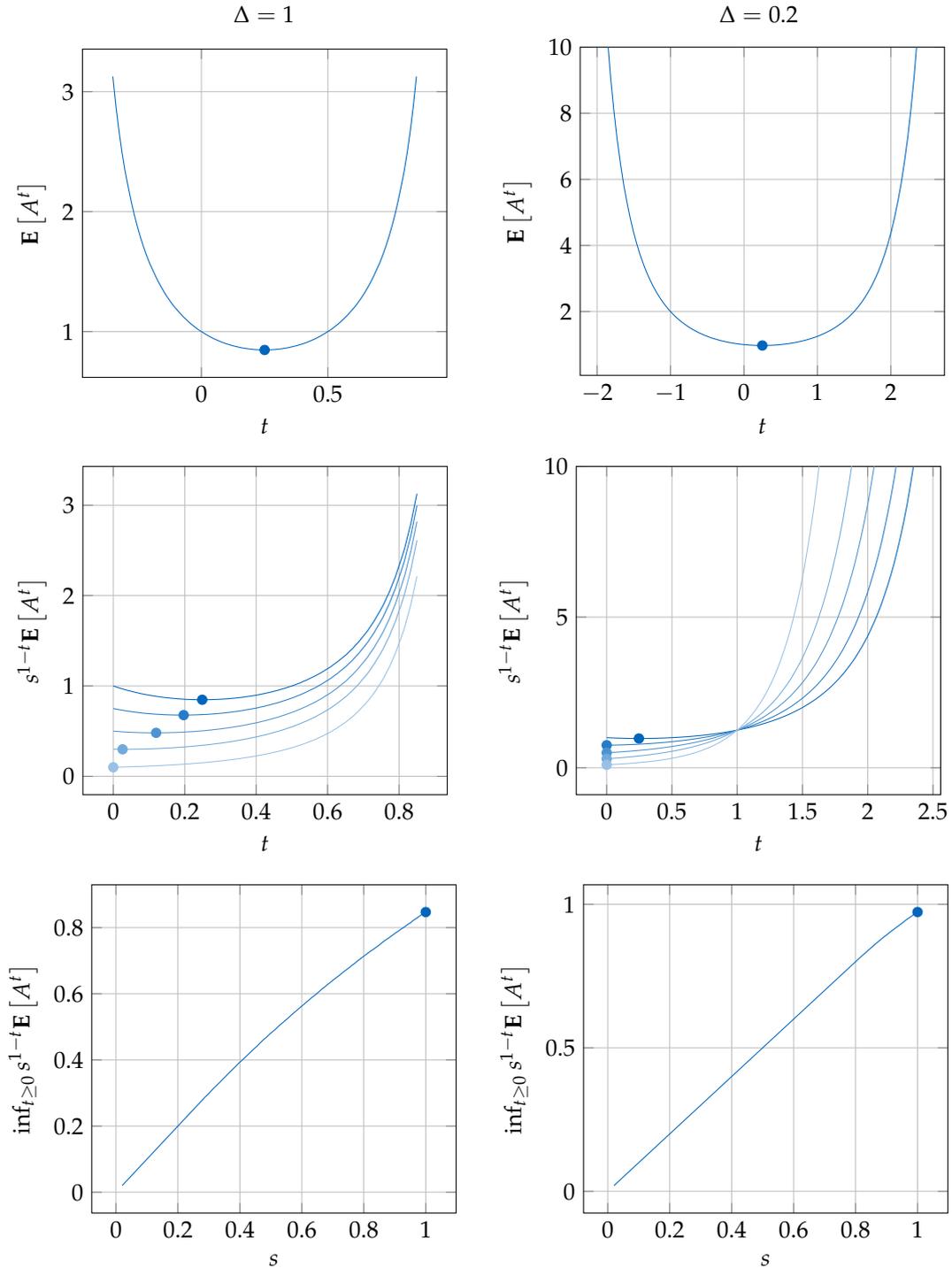

\clearpage

The characterization of recurrence and transience in \autoref{thm:transrecrwre} depends on the
quantity $p$, which is defined as $p := \min_{0 \leq t \leq 1} \bEx{A^t}$ for some random variable
$A$ such that $0 < A < \infty$ a.s. The following lemma reformulates the expression for
$p$, which will be useful later.

\begin{lemma}[Alternative Expression for $p$ \textnormal{(\cite[page 129]{rwrelyonspemantle})}]
  \label{lem:altp}
  Let $A$ be a real-valued random variable such that $0 < A < \infty$ a.s. Then
  \begin{align}
    \label{equ:altp}
    \min\limits_{0 \leq t \leq 1} \bEx{A^t} = \max\limits_{0 < s \leq 1} \inf\limits_{t \geq 0} s^{1-t} \bEx{A^t}
  \end{align}
\end{lemma}

\begin{tproof}
  The direction $\geq$ of \Cref{equ:altp} is easy to prove:
  \begin{align*}
    \textrm{for } 0 < s \leq 1 \textrm{: } \min\limits_{0 \leq t \leq 1} \bEx{A^t}
    \;\;\overset{s^{1-t} \leq 1 \textrm{ for } t \leq 1}{\geq}\;\; \min\limits_{0 \leq t \leq 1} s^{1-t} \bEx{A^t}
    \geq \inf_{t \geq 0} s^{1-t} \bEx{A^t}
  \end{align*}
  Since this holds for all $0 < s \leq 1$, it follow that $\geq$ holds in \Cref{equ:altp}.

  For the reverse direction, let $t_0 := \argmin_{0 \leq t \leq 1} \bEx{A^t}$. Since $A > 0$,
  $\bEx{A^t}$ is strictly convex in $t$ if not $A \equiv 1$ a.s.~(and it is trivial to
  verify \Cref{equ:altp} in the case $A \equiv 1$ a.s.). Consider the following three cases:
  \begin{itemize}
    \item $0 < t_0 < 1$. Then, by the strict convexity, $t_0$ will be the unique global minimum of
    $\bEx{A^t}$ and since $s^{1-t} \geq 1$ for $0 < s \leq 1$ and $t > 1$, it follows that
    $\inf_{t \geq 0} s^{1-t} \bEx{A^t} = \min\limits_{0 \leq t \leq 1} s^{1-t} \bEx{A^t}$.
    Choosing $s = 1$ shows $\leq$ in \Cref{equ:altp}.
    \item $t_0 = 0$. Then $\bEx{A^t} > 1$ for $t > 0$. It follows that $\inf_{t \geq 0} s^{1-t} \bEx{A^t} = s$
    for $0 < s \leq 1$ and again choosing $s = 1$ shows $\leq$ in \Cref{equ:altp}.
    \item $t_0 = 1$. This case is not as simple as the other two. For the $A$ used later in
    connection with the ERRW, the discussion in \autoref{sec:errw_z} will show that
    $t_0 = \frac{1}{4}$, hence this case is not relevant for the central statements on the
    ERRW. Therefore, this part is not proved here, but a proof can be found in
    \cite{rwrelyonspemantle}. \qedhere
  \end{itemize}
\end{tproof}

The plots in \autoref{fig:altp} on page \pageref{fig:altp} show the involved expressions
for the $A$ which will be used later in connection with the ERRW. See \autoref{sec:errw_z}
for more details on the distribution of $A$.

\begin{theorem}[Transience and Recurrence of RWRE \textnormal{(\cite[Theorem 1]{rwrelyonspemantle})}]
  \label{thm:transrecrwre}
  Consider a RWRE on a tree $T = \left(V, E\right)$ with the notation from
  \autoref{def:rwre}. In particular, the following assumptions are made:
  $\gls{av}, v \in V \setminus \left(\left\{\gls{root}\right\} \cup \glslink{distset}{T_1}\right)$
  are identically distributed, and independent for nodes with different parent nodes. $A$ is a
  random variable with the same distribution. Furthermore, w.l.o.g.~$0 < A_{\glslink{child}{\rho_i}} \leq 1$
  a.s. for $1 \leq i \leq \glslink{numch}{b_{\rho}}$. Finally, assume $0 < A < \infty$
  a.s.~(every node can a.s.~be reached with positive probability).
  
  Let $p := \min_{0 \leq t \leq 1} \bEx{A^t}$. Then
  \begin{enumerate}[(i)]
    \item \label{thm:rwre_trans} if $p \cdot \gls{br} > 1$, the RWRE is a.s.~transient.
    \item \label{thm:rwre_rec} if $p \cdot \br{T} < 1$ or (more generally)
    $\inf_{\gls{cutset} \textrm{ cutset}} \sum_{v \in \Pi} p^{\gls{nodedist}} = 0$,
    the RWRE is a.s.~recurrent.
    \item \label{thm:rwre_posrec} if $p \cdot \varlimsup_{n \to \infty} \left|\gls{distset}\right|^{\frac{1}{n}} < 1$
    or (more generally) $\sum_{v \in V} p^{\left|v\right|} < \infty$,
    the RWRE is a.s.~positive recurrent.
  \end{enumerate}
\end{theorem}

\begin{remark}
  The stated assumptions on the independence and identical distribution of the $A_v$ are stronger
  than necessary. Since transience and recurrence are maintained when the tree is modified in
  finitely many places, it is intuitively clear that it suffices to assume identical distribution
  up to a finite number of exceptions and independence for pairs of nodes with different parents
  up to a finite number of exceptional pairs. However, the proof becomes more cumbersome in this
  case, and since it is already long enough in its present, detailed version, it is left to the
  reader to check this.
\end{remark}

\begin{tproof}
  The three statements are proved separately, following \cite{rwrelyonspemantle} but with additional
  details and the proof now only relies on \cite{probtreenet} as external source.
  \begin{enumerate}[(i)]
    \item The idea of the proof is to a.s.~find a subtree on which the random walk is
    transient, which implies transience of the random walk on the original tree. For the proof of
    this part, assume that all $A_v$ are identically distributed. At the end of the proof of
    \ref{thm:rwre_trans}, it will be evident that the proof also holds for the stated assumptions. By
    \autoref{lem:altp}, it is possible to choose $s \in \left(0, 1\right]$ such that
    $p = \inf\limits_{t \geq 0} s^{1-t} \bEx{A^t}$. Furthermore, by \autoref{lem:cct}
    \ref{lem:cct_mul}, for $v \in V$ with $\left|v\right| = n \geq 1$,
    \begin{align}
      \label{equ:cvprob}
      \bPrb{\gls{cv} \geq s^n}^{\frac{1}{n}} = \bPrb{\prod_{\rho \neq \glslink{ancestor}{u \leq v}} A_v \geq s^n}^{\frac{1}{n}} \overset{n \to \infty}{\longrightarrow} \inf_{t \geq 0} s^{-t} \cdot \bEx{A^t}
    \end{align}
    Note that
    \begin{align*}
      \inf_{t \geq 0} s^{-t} \cdot \bEx{A^t} = s^{-1} \inf_{t \geq 0} s^{1-t} \cdot \bEx{A^t} = s^{-1} p
    \end{align*}
    and that $s^{-1} p > \left(s \br{T}\right)^{-1}$ since $p \cdot \br{T} > 1$ by assumption. The
    limit on the right hand side of \Cref{equ:cvprob} is $s^{-1} p$, and it is strictly larger
    than $\left(s \br{T}\right)^{-1}$. This implies in particular that there must be some $n \geq 1$
    such that
    \begin{align*}
      \bPrb{C_v \geq s^n}^{\frac{1}{n}} > \left(s \br{T}\right)^{-1}
      \quad \implies \quad \bPrb{C_v \geq s^n} > \left(s \br{T}\right)^{-n}
    \end{align*}
    for all nodes $v \in V$ with $\left|v\right| = n$. As the events
    $\left\{C_v \geq s^n \textrm{ and } \forall \rho \neq u \leq v : A_u \geq \varepsilon\right\}$
    converge to the event $\left\{C_v \geq s^n\right\}$ for $\varepsilon \searrow 0$ (recall
    that $A > 0$ a.s.), an $\varepsilon > 0$ which is small enough can be chosen such that
    \begin{align}
      \label{equ:q}
      q := \bPrb{C_v \geq s^n \textrm{ and } \forall \rho \neq u \leq v : A_u \geq \varepsilon} > \left(s \br{T}\right)^{-n}
    \end{align}
    for all $v \in T_n$.

    Now, the transient subgraph of $T$ can be constructed, and $q$ is used to show that such an
    infinite subgraph exists a.s. Consider the tree $\gls{ntree}$ which contains all nodes
    $v \in V$ with $n \divides \left|v\right|$ where two nodes $u, v$ are connected by an edge if,
    and only if, $u \leq v$ and $\left|u\right| + n = \left|v\right|$ (where all distances to the
    root are measured in the original tree $T$). By \autoref{lem:tnbr}, $\br{T^{\left(n\right)}} = \br{T}^n$.

    Let $\psi$ be the random subgraph of $T^{\left(n\right)} = \left(V^{\left(n\right)}, E^{\left(n\right)}\right)$
    where an edge $\left\{u, v\right\} \in E^{\left(n\right)}$ with $u < v$ is deleted if
    \begin{align*}
      \prod_{u < w \leq v, w \in V} A_w < s^n \textrm{ or there is } w \in V \textrm{ with } u < w \leq v \textrm{ and } A_w < \varepsilon
    \end{align*}
    The events that two edges are deleted are independent whenever the starting nodes of the edges
    (the nodes closer to the root) are different, and the probability that an edge is deleted is
    $1 - q$ by the assumptions on independence and distribution of the $A_v$ (still assumed to be
    all identically distributed).
    
    By \Cref{equ:q}, $q \cdot \left(s \br{T}\right)^n > 1$, hence $q^{\frac{1}{n}}s\br{T} > 1$.
    It is therefore possible to choose $r \in \left(\left(q^{\frac{1}{n}}s\br{T}\right)^{-1}, 1\right)$.
    Since $r < 1$ and $s \leq 1$, it holds that $\left(rs\right)^{-n} > 1$. On the other hand,
    since $r > \left(q^{\frac{1}{n}}s\br{T}\right)^{-1}$, it also holds that
    $\left(rs\right)^{-1} < q^{\frac{1}{n}}\br{T}$, so $q \br{T^{\left(n\right)}} = q \br{T}^n > \left(rs\right)^{-n} > 1$.

    By \autoref{lem:percolation} (proved below on page \pageref{lem:percolation}), there is
    a.s.~a subtree $\psi^\ast$ of $\psi$ with $\br{\psi^\ast} > \left(rs\right)^{-n}$. Now,
    $\psi^\ast$ is a subtree of $T^{\left(n\right)}$. Consider the corresponding subtree
    $T^\ast = \left(V^\ast, E^\ast\right)$ of $T$ (i.e.~the graph with vertices and edges
    which lie between two nodes contained in $\psi^\ast$). By construction,
    \begin{enumerate}[(a)]
      \item \label{prf:asexist} $T^\ast$ with the following three properties exists
      a.s.~(for almost all values of the random variables $A_v$):
      \item \label{prf:brt} $\br{T^\ast} > \left(rs\right)^{-1}$ (because of \autoref{lem:tnbr} and
      $\br{\psi^\ast} > \left(rs\right)^{-n}$)
      \item \label{prf:av} $A_v \geq \varepsilon$ for all $v \in V^\ast$
      \item \label{prf:prodav} $\prod_{u < w \leq v, w \in V^\ast} A_w \geq s^n$ for $u,v \in V^\ast$ with $n \divides \left|u\right| = \left|v\right| - n$
    \end{enumerate}
    
    \begin{figure}[H]
      \begin{center}
        \def\svgwidth{5.4cm}
        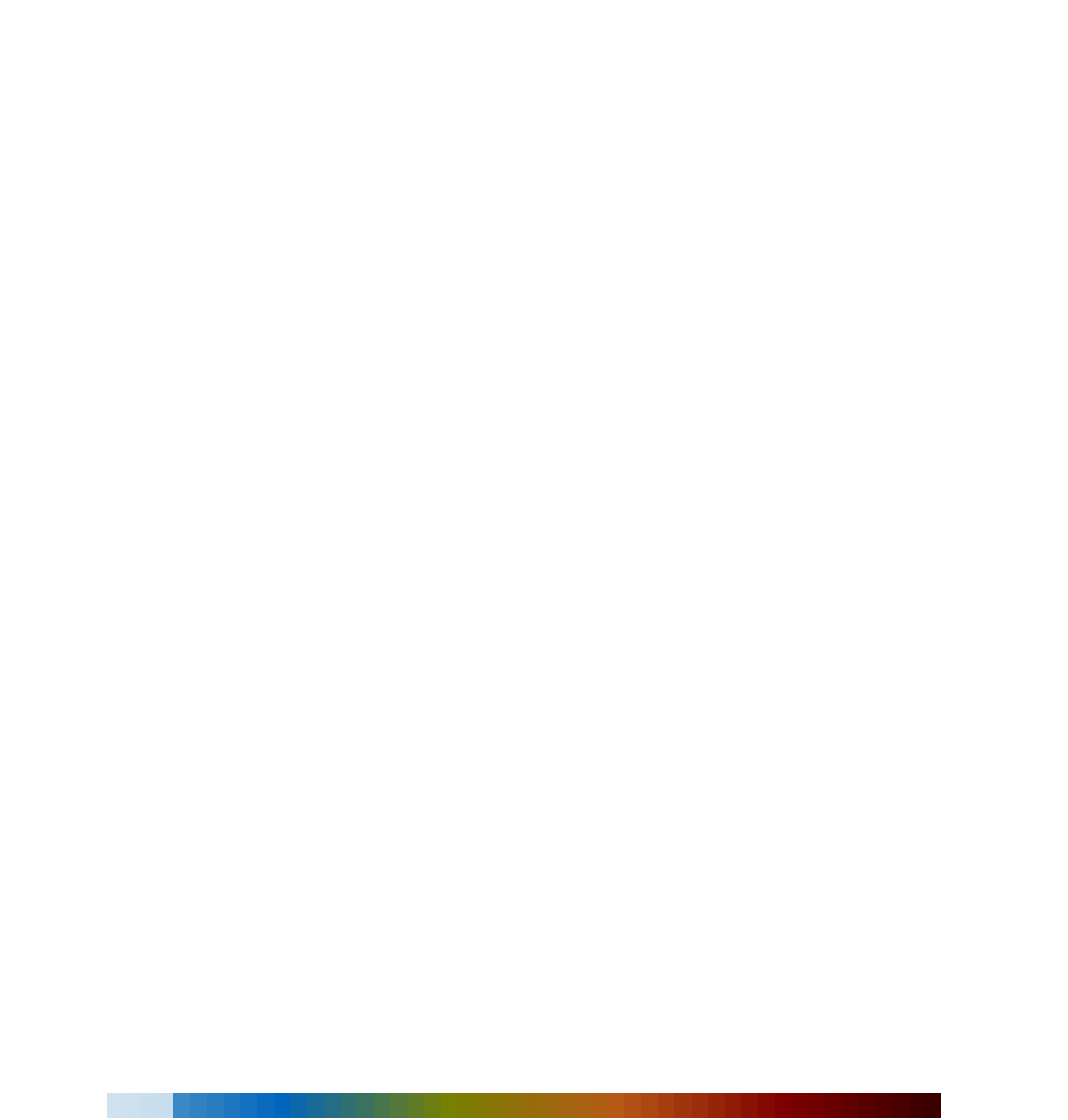
        \caption{Construction of the transient subtree}
        \parbox{\textwidth}{An example where $T$ is a Galton-Watson tree with
        $p_1 = 0.7, p_2 = 0.15, p_3 = 0.1, p_4 = 0.05$, hence $\br{T} = m = 1.5$. The values of the
        variables $A_v$ were randomly selected according to the common distribution which will later
        be relevant for the ERRW, and which depends on the reinforcement parameter $\Delta$
        (the condition from \autoref{thm:transrecrwre} \ref{thm:rwre_trans} holds).
        $\log\left(1 + 100A_v\right)$ is color coded. The bold nodes belong to $T^{\left(n\right)}$
        and the fully colored parts of the tree fulfil conditions \ref{prf:av} and \ref{prf:prodav}.}
        \label{fig:transsub}
      \end{center}
    \end{figure}

    Therefore, $C_v \geq s^{\left|v\right|} \varepsilon^{n-1}$ by \ref{prf:av}, \ref{prf:prodav} and
    the fact that $s \leq 1$. Hence
    \begin{align*}
      \inf_{\Pi \textrm{ cutset in } T^\ast} \sum_{v \in \Pi} r^{\left|v\right|} C_v
      \geq \varepsilon^{n-1} \inf_{\Pi \textrm{ cutset in } T^\ast} \sum_{v \in \Pi} \left(rs\right)^{\left|v\right|}
      > 0
    \end{align*}
    where the last expression is $> 0$ by \ref{prf:brt} and the alternative expression for the
    branching number in \Cref{equ:br_alt}. By the max-flow min-cut theorem, there is a non-zero
    flow $f$ from the root of $T^\ast$ to $\infty$ obeying the capacity constraints
    $\left|f\left(\gls{paredge}\right)\right| \leq r^{\left|v\right|}C_v$.

    As $r < 1$, $\sum_{n \geq 1} r^n < \infty$. Applying \cite[Proposition 3.4]{probtreenet} to the
    non-zero flow $f$ yields the existence of a finite energy flow (definition on
    \cite[page 38]{probtreenet}, but not relevant here) from the root of $T^\ast$ to $\infty$,
    and therefore also the existence of a unit flow to $\infty$ with finite energy. By
    \cite[Theorem 2.11]{probtreenet}, this implies that the random walk on $T^\ast$ is transient.
    As $T^\ast$ is a subtree of $T$, the RWRE on $T$ is transient as well
    (a.s.~by \ref{prf:asexist}).

    Finally, note that due to the assumption that the $A_v$ might be distributed differently
    for $v \in T_1$, $T^\ast$ might not have the stated properties \ref{prf:av} and \ref{prf:prodav}
    if $T^\ast$ contains the original root $\rho$ of $T$. But in this case, it is possible to simply
    remove $\rho$ from $\psi^\ast$, which removes $\rho$ and all nodes with distance $< n$ from
    $T^\ast$. $T^\ast$ might decompose into multiple trees when deleting these nodes, but each of
    them will have properties \ref{prf:av}, \ref{prf:prodav} and \ref{prf:asexist} is also maintained.
    Finally, at least one of them will have property \ref{prf:brt} (this is easy to prove, and also
    treated in \cite[Example 1.3]{probtreenet}). Hence, the subtree with the required properties
    still exists.

    \item The proof relies on \cite[Theorem 2.3]{probtreenet}: the random walk is transient if,
    and only if, the effective conductance from $\rho$ to $\infty$ is positive, which is equivalent
    to the effective resistance from $\rho$ to $\infty$ being finite (the definitions of effective conductance
    and resistance can be found on \cite[page 27]{probtreenet}, they are not important here). Hence, the random walk is
    recurrent if the effective resistance is infinite. By the Nash-Williams criterion from
    \cite[page 37]{probtreenet}, the effective resistance is larger than or equal to
    $\left(\sum_{v \in \Pi} C_v\right)^{-1}$ for any cutset $\Pi$. If
    $\inf_{\Pi \textrm{ cutset}} \sum_{v \in \Pi} C_v = 0$, then
    $\sup_{\Pi \textrm{ cutset}} \left(\sum_{v \in \Pi} C_v\right)^{-1} = \infty$, implying
    infinite effective resistance, and therefore recurrence.

    First consider the case where $p \cdot \br{T} < 1$. Then $\br{T} < p^{-1}$. Consequently, by
    the alternative expression for the branching number from \Cref{equ:br_alt},
    $\inf_{\Pi \textrm{ cutset}} \sum_{v \in \Pi} p^{\left|v\right|} = 0$.
    Thus, to proof \ref{thm:rwre_rec}, it is sufficient to show that $\inf_{\Pi \textrm{ cutset}} \sum_{v \in \Pi} p^{\left|v\right|} = 0$
    implies a.s.~recurrence.
    For a particular $0 < t \leq 1$ and $p = \bEx{A^t}$ with $\inf_{\Pi \textrm{ cutset}} \sum_{v \in \Pi} p^{\left|v\right|} = 0$,
    it holds that
    \begin{align*}
      \bEx{\sum_{v \in \Pi} C_v^t} &= \sum_{v \in \Pi} \bEx{\prod_{\rho \neq u \leq v} A_u^t} \overset{\textrm{independence}}{=}
      \sum_{v \in \Pi} \prod_{\rho \neq u \leq v} \bEx{A_u^t} \\
      &\overset{0 < A_{\glslink{child}{\rho_i}} \leq 1}{\leq}\;\;
      \sum_{v \in \Pi} \prod_{u \leq v, \left|u\right| > 1} \bEx{A^t} = \sum_{v \in \Pi} p^{\left|v\right| - 1} = p^{-1} \sum_{v \in \Pi} p^{\left|v\right|}
    \end{align*}
    By Fatou's lemma, it holds that
    \begin{align*}
      \bEx{\inf_{\Pi \textrm{ cutset}} \sum_{v \in \Pi} C_v^t} \leq \inf_{\Pi \textrm{ cutset}} \bEx{\sum_{v \in \Pi} C_v^t}
      \leq p^{-1} \cdot \inf_{\Pi \textrm{ cutset}} \sum_{v \in \Pi} p^{\left|v\right|} = 0
    \end{align*}
    Choosing $t = \argmin_{0 \leq t \leq 1} \bEx{A^t} > 0$ (where $> 0$ follows since
    $p = \min_{0 \leq t \leq 1} \bEx{A^t}$ cannot be $1$ by the assumption in \ref{thm:rwre_rec}),
    and noting that for $C_v < 1$, $C_v \leq C_v^t$, as well as that there must be a sequence of
    cutsets $\Pi_n$ where $C_v^t < 1$ for all $v \in \Pi_n$ eventually, one gets
    \begin{align*}
      \inf_{\Pi \textrm{ cutset}} \sum_{v \in \Pi} C_v \leq \inf_{\Pi \textrm{ cutset}} \sum_{v \in \Pi} C_v^t
      = 0 \textrm{ a.s.}
    \end{align*}
    This concludes the proof, see first paragraph.

    \item The following criterion is given in \cite[Exercise 2.1 (f)]{probtreenet}: a random
    walk is positive recurrent if, and only if, the sum of the conductances over all eges is finite
    (which is immediately clear since a finite sum of conductances is equivalent to the existence of
    a stationary measure, also see remarks \ref{term:statmeas} and \ref{term:rectrans} after \autoref{def:rwre}).

    First consider the case where $p \cdot \varlimsup_{n \to \infty} \left|T_n\right|^{\frac{1}{n}} < 1$.
    Then
    \begin{align*}
      \sum_{v \in V} p^{\left|v\right|} = \sum_{n \geq 0} \sum_{v \in T_n} p^n = \sum_{n \geq 0} \left(\left|T_n\right|^{\frac{1}{n}} \cdot p\right)^n < \infty
    \end{align*}
    where finiteness follows with the root test by the assumption.

    Thus, to proof \ref{thm:rwre_posrec}, it is sufficient to show that $\sum_{v \in V} p^{\left|v\right|} < \infty$
    implies a.s.~positive recurrence. As in the proof of \ref{thm:rwre_rec}, it holds
    that for a particular $0 < t \leq 1$ and $p = \bEx{A^t}$ with $\sum_{v \in V} p^{\left|v\right|} < \infty$ that
    $\bEx{\sum_{\rho \neq v \in V} C_v^t} \leq p^{-1} \sum_{\rho \neq v \in V} p^{\left|v\right|} < \infty$.
    Hence, $\sum_{\rho \neq v \in V} C_v^t < \infty$ a.s.~and therefore $C_v < 1$ for
    all but finitely many $v$ a.s. For $C_v < 1$, it holds that $C_v \leq C_v^t$, so one
    may conclude with $t = \argmin_{0 \leq t \leq 1} \bEx{A^t} > 0$ (again $> 0$ by the assumption
    in \ref{thm:rwre_posrec}) that
    \begin{align*}
      \sum_{\rho \neq v \in V} C_v = \underbrace{\sum_{\rho \neq v \in V, C_v < 1} C_v}_{\leq \sum_{\rho \neq v \in V} C_v^t < \infty}
      + \underbrace{\sum_{\rho \neq v \in V, C_v \geq 1} C_v}_{\textrm{finitely many terms, hence finite sum}} < \infty \textrm{ a.s.}
    \end{align*}
    The a.s.~finiteness of the sum of the conductances now implies a.s.~positive
    recurrence, see above. \qedhere
  \end{enumerate}
\end{tproof}

To complete the proof of the previous theorem, the following lemma is still necessary:

\begin{lemma}[Special Quasi-Independent Percolation and Branching Number]
  \label{lem:percolation}
  Let $T = \left(V, E\right)$ be an infinite tree with root $\rho$. Consider the following
  special case of quasi-independent percolation on $T$ (see \cite[Section 5.4]{probtreenet} for the general
  definition of quasi-independent percolation): an edge $\left\{u, v\right\} \in E$ is kept with
  probability $q > 0$ and deleted with probability $1 - q$. The events that an edge is deleted are
  further assumed to be independent for two edges whenever the starting nodes (the nodes closer to
  $\gls{root}$) of the edges are different. Denote the resulting random subgraph of $T$ by $\psi$.

  If $q \cdot \gls{br} > 1$, then for any $s < q \cdot \br{T}$, there exists a.s.~a
  subtree $\psi^\ast$ (possibly with root $\neq \rho$) in $\psi$ with $\br{\psi^\ast} > s$.
  More generally, there even exists a.s.~a subtree $\psi^\dagger$ (possibly with root
  $\neq \rho$) in $\psi$ with $\br{\psi^\dagger} = q \cdot \br{T}$ (but this statement is not
  further needed in this thesis).
\end{lemma}

\begin{tproof}
  Denote by $\rho \leftrightarrow v$ for $v \in V$ the event that $\rho$ is connected to $v$ in
  $\psi$ and by $K_v \subseteq V$ the connected component of $v \in V$ in $\psi$. The goal is to
  show $\Prb{\exists v \in V: \br{K_v} > s} = 1$, where $\br{K_v} = 0$ if the connected component
  is finite.
  
  Later, it will be necessary to demonstrate the a.s.~existence of a $v \in V$ with
  $\left|K_v\right| = \infty$. By Kolmogorov's zero-one law, it holds that
  $\Prb{\exists v \in V: \left|K_v\right| = \infty} \in \left\{0, 1\right\}$ because the event that
  an infinite subtree exists is independent of how any finite subgraph of $\psi$ looks
  (and by the independence assumption on the edge removal). Hence, it suffices to show that
  $\Prb{\exists v \in V: \left|K_v\right| = \infty} > 0$. But then, it is also sufficient to show that
  $\Prb{\left|K_\rho\right| = \infty} > 0$, because $\Prb{\left|K_\rho\right| = \infty} = 0$ would imply
  $\Prb{\left|K_v\right| = \infty} = 0$ for all $v \in V$.
  
  For $u, v \in V$, let $\gls{lca}$ be the \textbf{lowest common ancestor} of $u$ and $v$ ($u \land v$ is the node
  maximizing the distance to $\rho$ with the property that $u \land v$ lies both on the unique path from
  $\rho$ to $u$ and from $\rho$ to $v$). Let $w := u \land v$ and assume $u \neq w \neq v$. Then, there are distinct nodes
  $\overrightarrow{w^u}$ and $\overrightarrow{w^v}$ which are the next nodes after $w$ on the paths
  from $w$ to $u$ and $w$ to $v$, respectively, and it holds that (by the independence assumption
  on the edge removal)
  \begin{align*}
    \Prb{\rho \leftrightarrow u, \rho \leftrightarrow v}
    &= \underbrace{\Prb{\rho \leftrightarrow w}}_{q^{\left|w\right|}}
    \cdot \underbrace{\Prb{w \leftrightarrow \overrightarrow{w^u}, w \leftrightarrow \overrightarrow{w^v}}}_{\leq q}
    \cdot \underbrace{\Prb{\overrightarrow{w^u} \leftrightarrow u}}_{q^{\left|u\right| - \left|w\right| - 1}}
    \cdot \underbrace{\Prb{\overrightarrow{w^v} \leftrightarrow v}}_{q^{\left|v\right| - \left|w\right| - 1}}\\
    &\leq q^{\left|u\right| + \left|v\right| - \left|w\right| - 1}\\
    \implies \Prb{\rho \leftrightarrow u, \rho \leftrightarrow v} &\cdot \Prb{\rho \leftrightarrow w}
    \leq q^{\left|u\right| + \left|v\right| - 1}
    = q^{-1} q^{\left|u\right|} q^{\left|v\right|}
    = q^{-1} \Prb{\rho \leftrightarrow u} \Prb{\rho \leftrightarrow v}
  \end{align*}
  It is easy to check that the final inequality also holds in the case that $u = w$ or $v = w$ or
  $u = v = w$. Therefore, with $M := q^{-1}$, the considered percolation is a quasi-independent percolation by
  the definition in \cite[Section 5.4]{probtreenet}.

  Consider, for a moment, any edge-keeping probability $1 \geq q > 0$ (without the assumption $q \cdot \br{T} > 1$).
  Now, applying \cite[Theorem 5.19]{probtreenet}, and, in the words of the proof of said theorem,
  ``the proof of Theorem 5.[15] can be followed to the desired conclusion'' which is that
  for all $q > \frac{1}{\br{T}}$, it holds that $\Prb{\left|K_\rho\right| = \infty} > 0$
  (only using one direction of \cite[Theorem 5.15]{probtreenet}, slightly reformulated).
  \cite[Proposition 5.8]{probtreenet}, the fact that $\Prb{\rho \leftrightarrow v} = q^{\left|v\right|}$
  and the alternative expression of the branching number in \Cref{equ:br_alt} yield that
  for all $q < \frac{1}{\br{T}}$, it holds that $\Prb{\left|K_\rho\right| = \infty} = 0$.

  Finally, consider the case where, after the first edge removal with parameter $q$, a second,
  independent percolation where each edge is kept with probability $p$ is carried out. Assume
  that the events that two edges are deleted are independent for any two different edges this time,
  and also independent of the previous edge removal. Then, the resulting percolation is equivalent
  to a percolation with parameter $qp$ fulfilling the same independence assumptions as stated in
  \autoref{lem:percolation}. By the arguments above, it thus holds that
  $qp > \frac{1}{\br{T}} \implies \Prb{\left|K_\rho^{\left(qp\right)}\right| = \infty} > 0$ and
  $qp < \frac{1}{\br{T}} \implies \Prb{\left|K_\rho^{\left(qp\right)}\right| = \infty} = 0$ where $K_\rho^{\left(qp\right)}$ is the
  connected component of $\rho$ in the random graph $\psi_{qp}$, which is defined as $\psi$, but with
  removal parameter $qp$ instead of $q$.

  A Bernoulli percolation (see \cite[Section 5.2]{probtreenet}) with parameter $p$ on the graph $\psi$
  yields $\psi_{qp}$. For $q > \frac{1}{\br{T}}$ (as assumed in \autoref{lem:percolation}), it holds that $\Prb{\left|K_\rho\right| = \infty} > 0$,
  so there exists a.s.~a $v \in V$ with $\left|K_v\right| = \infty$ (check again the second paragraph
  of this proof). $K_v$ is a tree. But now, \cite[Theorem 5.15]{probtreenet} implies that
  $\frac{1}{\br{K_v}} = \sup\left\{ p : \Prb{\left|K_v^{\left(qp\right)}\right| = \infty} = 0 \right\}$.
  If $qp > \frac{1}{\br{T}}$, then $\Prb{\left|K_v^{\left(qp\right)}\right| = \infty} > 0$
  by the arguments above, and $qp < \frac{1}{\br{T}}$ implies $\Prb{\left|K_v^{\left(qp\right)}\right| = \infty} = 0$.
  Thus
  \begin{align*}
    \frac{1}{\br{K_v}} = \sup\left\{ p : \Prb{\left|K_v^{\left(qp\right)}\right| = \infty} = 0 \right\}
    = \frac{1}{q \br{T}} \implies
    \br{K_v} = q \br{T}
  \end{align*}
  So there exists a.s.~a subtree of $\psi$ with branching number $q \br{T}$. This concludes
  the proof.
\end{tproof}

\subsection{RWRE on Galton-Watson trees}
\label{ssec:rwregw}

\autoref{thm:transrecrwre} is already quite general. However, for the special case of Galton-Watson
trees, a little more is known: the critical case $p \cdot \br{T} = 1$ can be shown to imply recurrence
as well, while this is not true for general trees (see \cite[page 131]{rwrelyonspemantle}).

\begin{theorem}[RWRE on Galton-Watson Trees \textnormal{(\cite[Theorem 3]{rwrelyonspemantle})}]
  \label{thm:rwregw}
  Let $T$ be a Galton-Watson tree with mean $\gls{gwmean} > 1$. Consider a random environment on
  $T$ as in \autoref{thm:transrecrwre}, and let again $p := \min_{0 \leq t \leq 1} \bEx{A^t}$. Then
  \begin{enumerate}[(i)]
    \item \label{thm:rwregw_trans} If $pm > 1$, then, conditional on non-extinction of $T$, the
    RWRE is a.s.~transient.
    \item \label{thm:rwregw_rec} If $pm \leq 1$, then the RWRE is a.s.~recurrent.
    \item \label{thm:rwregw_posrec} If $pm < 1$, then the RWRE is a.s.~positive
    recurrent.
  \end{enumerate}
\end{theorem}

\begin{remark}
  In the case $A \equiv \lambda^{-1}$ where $\lambda > 0$ is constant, \autoref{thm:rwregw}
  fully characterizes recurrence and transience of the $\lambda$-biased random walk (see
  \cite[Theorem 1.7]{probtreenet}) on the Galton-Watson tree $T$.
\end{remark}

\begin{tproof}
  Again, the three parts are proved separately:
  \begin{enumerate}[(i)]
    \item By \autoref{lem:gwbr}, it holds that $\gls{br} = m$ a.s.~conditional on
    non-extinction. If $p \cdot \br{T} > 1$, then by \autoref{thm:transrecrwre} \ref{thm:rwre_trans}
    the RWRE is a.s.~transient.

    \item \autoref{thm:transrecrwre} \ref{thm:rwre_trans} only shows a.s.~recurrence if
    $pm < 1$. In order to extend this result to the case $pm = 1$, more work is required.

    As in the proof \autoref{thm:transrecrwre} \ref{thm:rwre_rec}, the Nash-Williams criterion is used, but in a
    slightly different way: by \cite[page 37]{probtreenet} and \cite[Theorem 2.3]{probtreenet},
    the random walk is recurrent if there is a series of disjoint cutsets $\glslink{cutset}{\Pi_n}, n \geq 1$
    such that $\sum_{n=1}^\infty \left(\sum_{v \in \Pi_n} C_v\right)^{-1} = \infty$. The following
    proof is basically a transcription of the proof given in \cite{rwrelyonspemantlecorr}, and the
    series of cutsets that is used will be $\Pi_n := \gls{distset}$. The expected value $\Ex{\cdot}$
    should be understood as taking the expectation with respect to the combined probability space of
    possible Galton-Watson trees $T$ and the possible values of the random variables $A_v$, where
    $T$ should now be seen as (tree-valued) random variable as well.
    
    Consider $t = \argmin_{0 \leq t \leq 1} \bEx{A^t} > 0$ ($> 0$ by the same argument as in
    the proof of \autoref{thm:transrecrwre} \ref{thm:rwre_rec}). As in the proof of
    \autoref{thm:transrecrwre} \ref{thm:rwre_rec}, it holds that
    \begin{align*}
      \Ex{\sum_{\gls{nodedist} = n} C_v^t} &= \Ex{\Ex{\sum_{v \in T_n} C_v^t \midd| T}} \\
      & \leq \Ex{p^{-1} \sum_{v \in T_n} p^{\left|v\right|}}
      = p^{-1} \Ex{\left|T_n\right|} p^n = p^{-1} \left(mp\right)^n \leq p^{-1}
    \end{align*}
    where the last equality follows from \cite[Proposition 5.5]{probtreenet}, and the final inequality
    follows from the assumption.

    In order to show $\sum_{n=1}^\infty \left(\sum_{v \in T_n} C_v\right)^{-1} = \infty$ a.s.,
    it is sufficient to prove that the probability of this expression being smaller than any constant $L > 0$ is zero.
    For this purpose, the following consequence of the power mean
    inequality\footnote{see e.g.~\url{https://en.wikipedia.org/wiki/Generalized_mean\#Inequality_between_any_two_power_means}}
    will be helpful: for $s \geq -1$ (so in particular for $s = t$), and $a_1, \ldots, a_N > 0$,
    \begin{align*}
      \left(\frac{1}{N}\sum_{n=1}^N a_n^{-1}\right)^{-1} \leq \left(\frac{1}{N}\sum_{n=1}^N a_n^s\right)^\frac{1}{s}
    \end{align*}
    Setting $s := t$, taking the $t$-th power and setting $a_n := \sum_{v \in T_n} C_v$ yields
    \begin{align*}
      \left(\frac{1}{N} \sum_{n=1}^N \left(\sum_{v \in T_n} C_v\right)^{-1}\right)^{-t}
      \leq \frac{1}{N}\sum_{n=1}^N \left(\sum_{v \in T_n} C_v\right)^t
      \leq \frac{1}{N}\sum_{n=1}^N \sum_{v \in T_n} C_v^t
    \end{align*}
    where the last inequality follows, for example, from the fact that the $\left\lVert \cdot \right\rVert_t$-norm
    of a fixed vector is non-increasing in $t$ (hence $\left\lVert \cdot \right\rVert_1 \leq \left\lVert \cdot \right\rVert_t$
    for $t \leq 1$). For any $N \in \bbN$, it now holds that
    \begin{align*}
      \Prb{\sum_{n=1}^\infty \left(\sum_{v \in T_n} C_v\right)^{-1} \leq L}
      &\leq \Prb{\sum_{n=1}^N \left(\sum_{v \in T_n} C_v\right)^{-1} \leq L}\\
      &= \Prb{\left(\frac{1}{N} \sum_{n=1}^N \left(\sum_{v \in T_n} C_v\right)^{-1}\right)^{-t} \geq \left(\frac{N}{L}\right)^t} \\
      &\overset{\textrm{Markov inequ.}}{\leq}
      \Ex{\left(\frac{1}{N} \sum_{n=1}^N \left(\sum_{v \in T_n} C_v\right)^{-1}\right)^{-t}} \cdot \left(\frac{L}{N}\right)^t \\
      &\overset{\textrm{PM inequ.}}{\leq} \Ex{\frac{1}{N} \sum_{n=1}^N \sum_{v \in T_n} C_v^t} \cdot \left(\frac{L}{N}\right)^t \\
      &\leq \left(\frac{1}{N} \sum_{n=1}^N p^{-1}\right) \cdot \left(\frac{L}{N}\right)^t
      = p^{-1}\left(\frac{L}{N}\right)^t
    \end{align*}
    Since $N \in \bbN$ was arbitrary, one may conclude that $\Prb{\sum_{n=1}^\infty \left(\sum_{v \in T_n} C_v\right)^{-1} \leq L} = 0$
    for all $L > 0$, which concludes the proof.

    \item This part is not proved here, because it is not used anymore in this thesis. Check
    \cite[Theorem 3 (iii)]{rwrelyonspemantle} for details (note however that in \cite{rwrelyonspemantle},
    it is only said that the proof is similar to the proof of \autoref{thm:transrecrwre} \ref{thm:rwre_posrec}). \qedhere
  \end{enumerate}
\end{tproof}

\clearpage
\newpage

\section{Transience and Recurrence of ERRW}
\label{sec:transrecerrw}

The results for transience and recurrence of RWRE help to analyze ERRWs.
\cite[Lemma 1 and the following paragraphs]{errwpemantle} provide the necessary connection between
ERRWs and RWRE. Consider an infinite tree $T$ and the ERRW with
reinforcement parameter $\gls{reinforcement}$ on $T$. For any node $v$ other than the root, the
edge weights $\glslink{edgeweight}{w_n}$ will be $1 + \Delta$ for the edge $\gls{paredge}$ to
$\gls{par}$, and $1$ for every edge to a child $\gls{child}$ with $1 \leq i \leq \gls{numch}$ when
the node is visited for the first time.

\begin{figure}[H]
  \begin{center}
    \begin{tikzpicture}
      \node[circle,fill=black,inner sep=0.7mm] (N0) at (0, 0) {};
      \node[circle,fill=black,inner sep=0.7mm] (N1) at (2, 0) {};
      \node[circle,fill=black,inner sep=0.7mm] (N2) at (4, 1) {};
      \node (Nb) at (4, 0) {$\vdots$};
      \node[circle,fill=black,inner sep=0.7mm] (N3) at (4, -1) {};
      \node[below=2mm] at (N0) {$\overleftarrow{v}$};
      \node[below=2mm] at (N1) {$v$};
      \node[above=2mm] at (N2) {$v_1$};
      \node[below=2mm] at (N3) {$v_{b_v}$};
      \draw (N0) -- node[above] {$1 + \Delta$} (N1) -- node[above] {$1$} (N2);
      \draw (N1) -- node[above] {$1$} (N3);
    \end{tikzpicture}
    \caption{First visit to a node $v$ on the ERRW}
    \label{fig:firstvis}
  \end{center}
\end{figure}
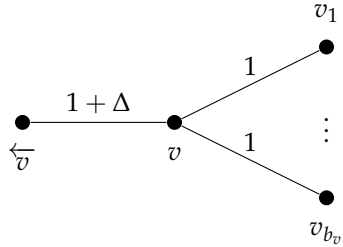

The node will then be left along one of the edges. If the node is visited again, then the ERRW
will take that same edge again when it returns to $v$ because $T$ is a tree. The edge weight of
the edge will have increased by $2\Delta$. This process is closely related to drawing balls from an
urn. For $\Delta = 1$, imagine an urn which initially contains $2$ black balls as well as $1$ ball
for every edge to a child, all with different colors. A ball is then drawn randomly from the urn,
and replaced by $3$ balls of the same color as the ball which was drawn (the urn now contains $2$
more balls). This process is repeated. The probability to take the edge to the parent node at each
visit to the node $v$ corresponds exactly to the probability to draw a black ball from the urn, if
a ball is drawn whenever the node $v$ is visited, and the edge corresponding to the drawn ball is
taken. The same holds for the other colors and the edges they correspond to.

In general, $\Delta$ does not need to be an integer, but the idea is the same. The urn model was
extensively studied and the sequence of draws was shown to be distributed just as a mixture of
iid draws where each color $i$ is chosen with probability $q_i$ in every draw, and
$q_i$ is a random variable. $q$, seen as a vector with components $q_i$, has a Dirichlet distribution
whose parameters depend on $\Delta$; the exact distribution can be found in \cite[Lemma 1]{errwpemantle}.

The ERRW can thus be simulated by placing independent urns at every node, and choosing
the next edge according to the ball which was drawn from the urn. The urns, in turn, can be replaced
by a mixture of iid draws. This makes it possible to define a RWRE which has
the exact same distribution as the ERRW: the RWRE is defined on the same
graph, and the random environment $\omega$ is given by the random vector $q$ described above.

\begin{figure}
  \begin{center}
    \def\svgwidth{5.6cm}
    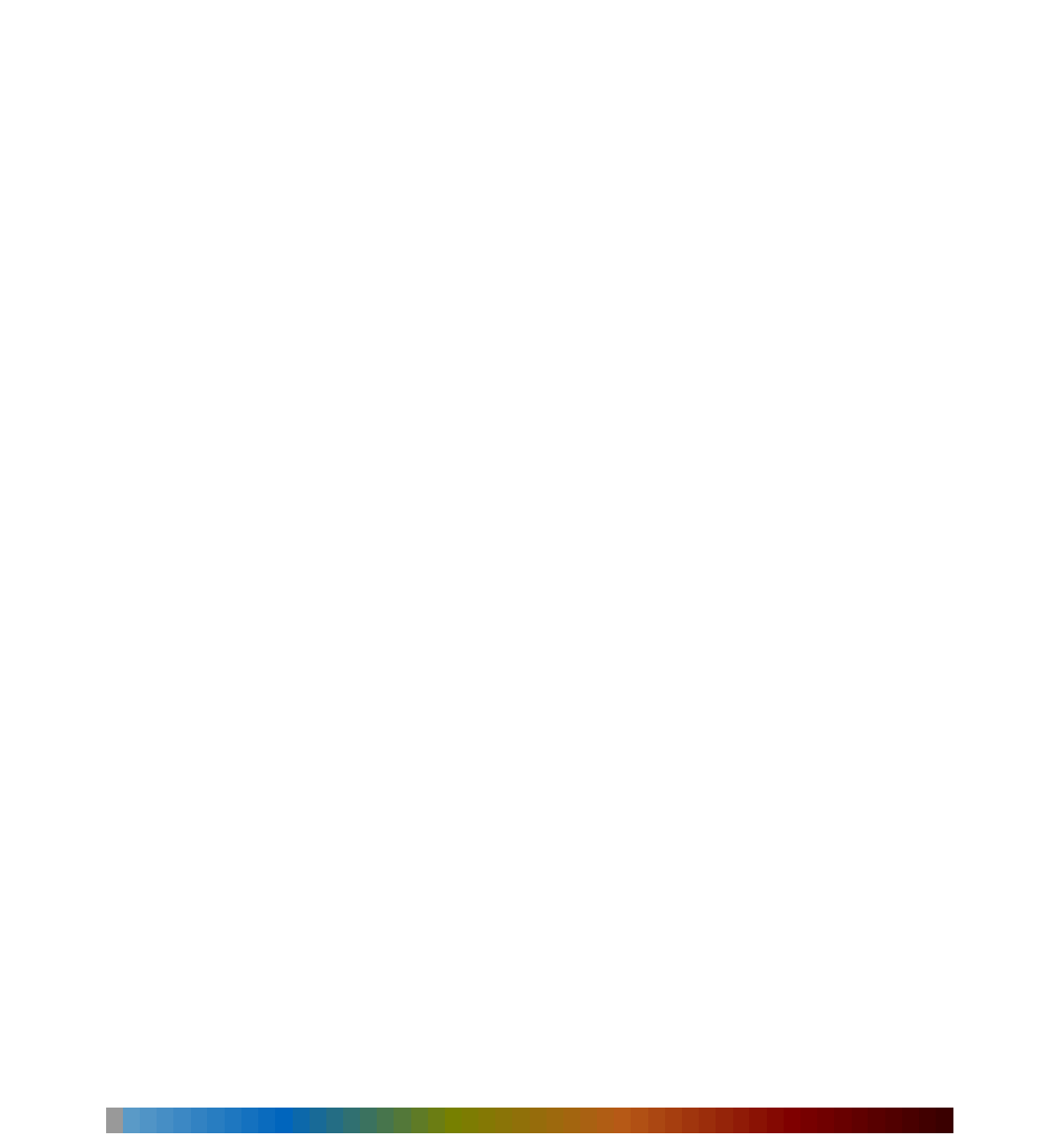 ~~~~~~~~~~~
    \def\svgwidth{5.6cm}
    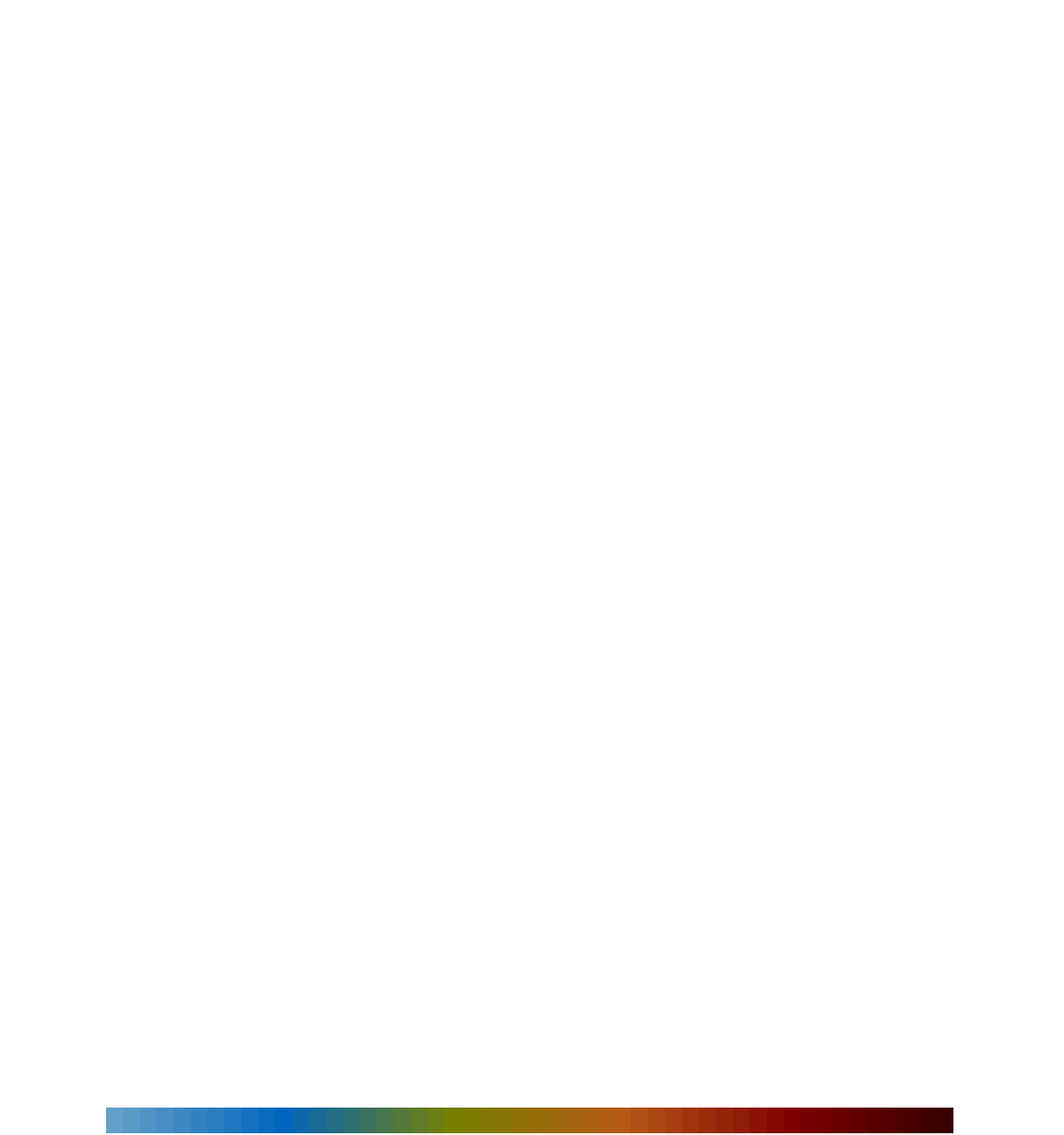 \\[2em]
    \def\svgwidth{5.8cm}
    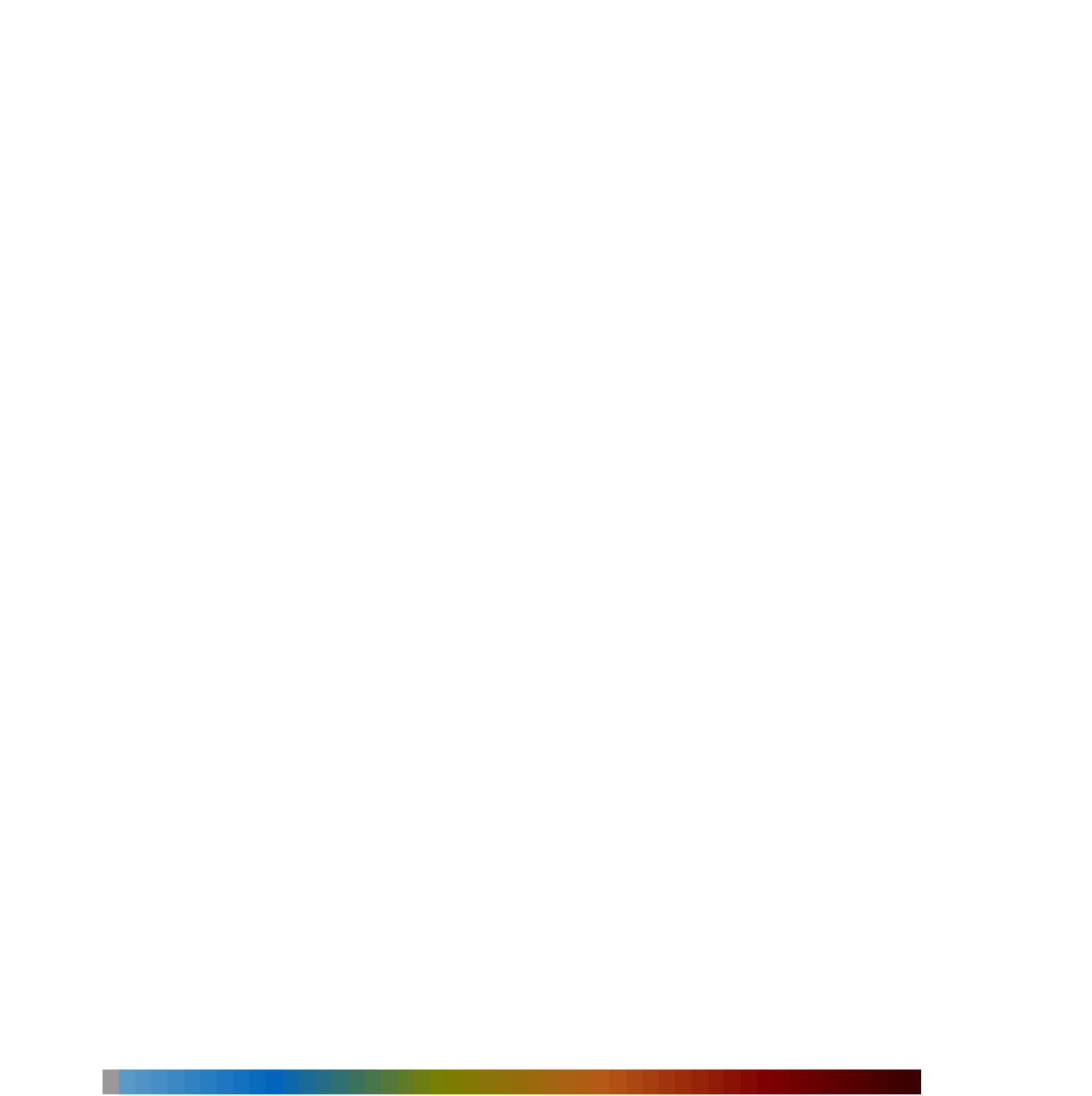 ~~~~~~~~~~~
    \def\svgwidth{5.6cm}
    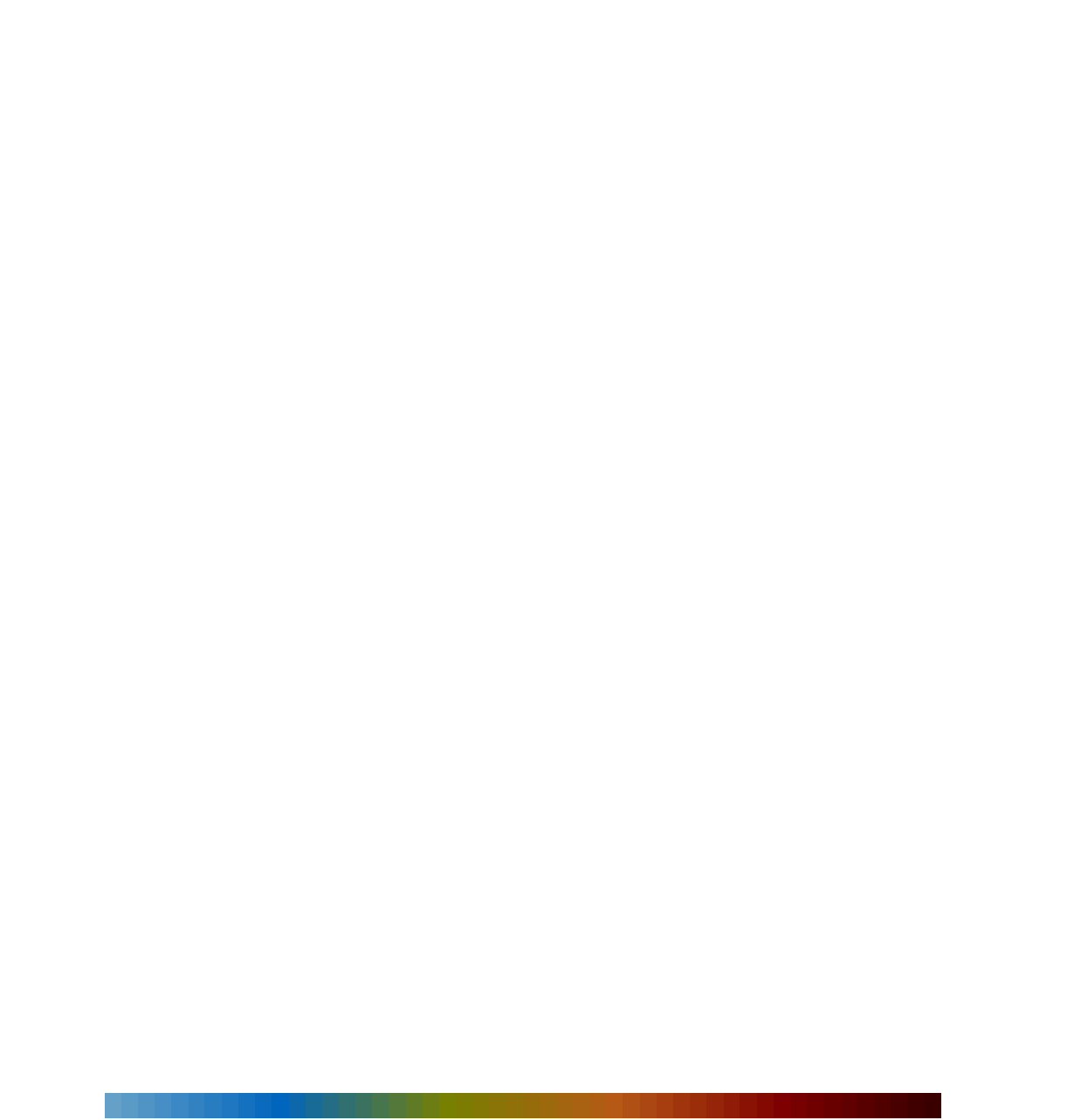 \\[2em]
    \def\svgwidth{5.8cm}
    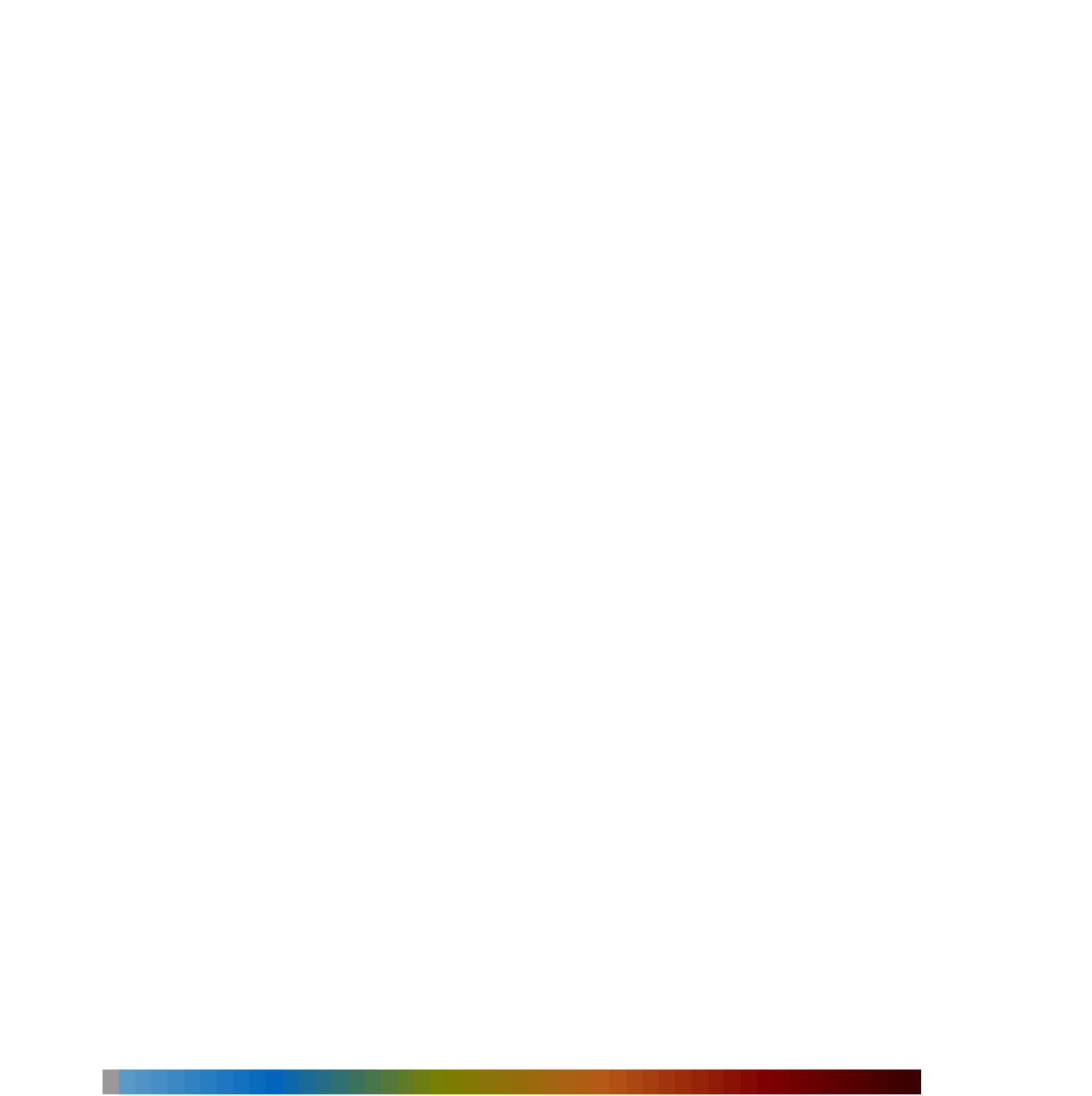 ~~~~~~~~~~~
    \def\svgwidth{5.6cm}
    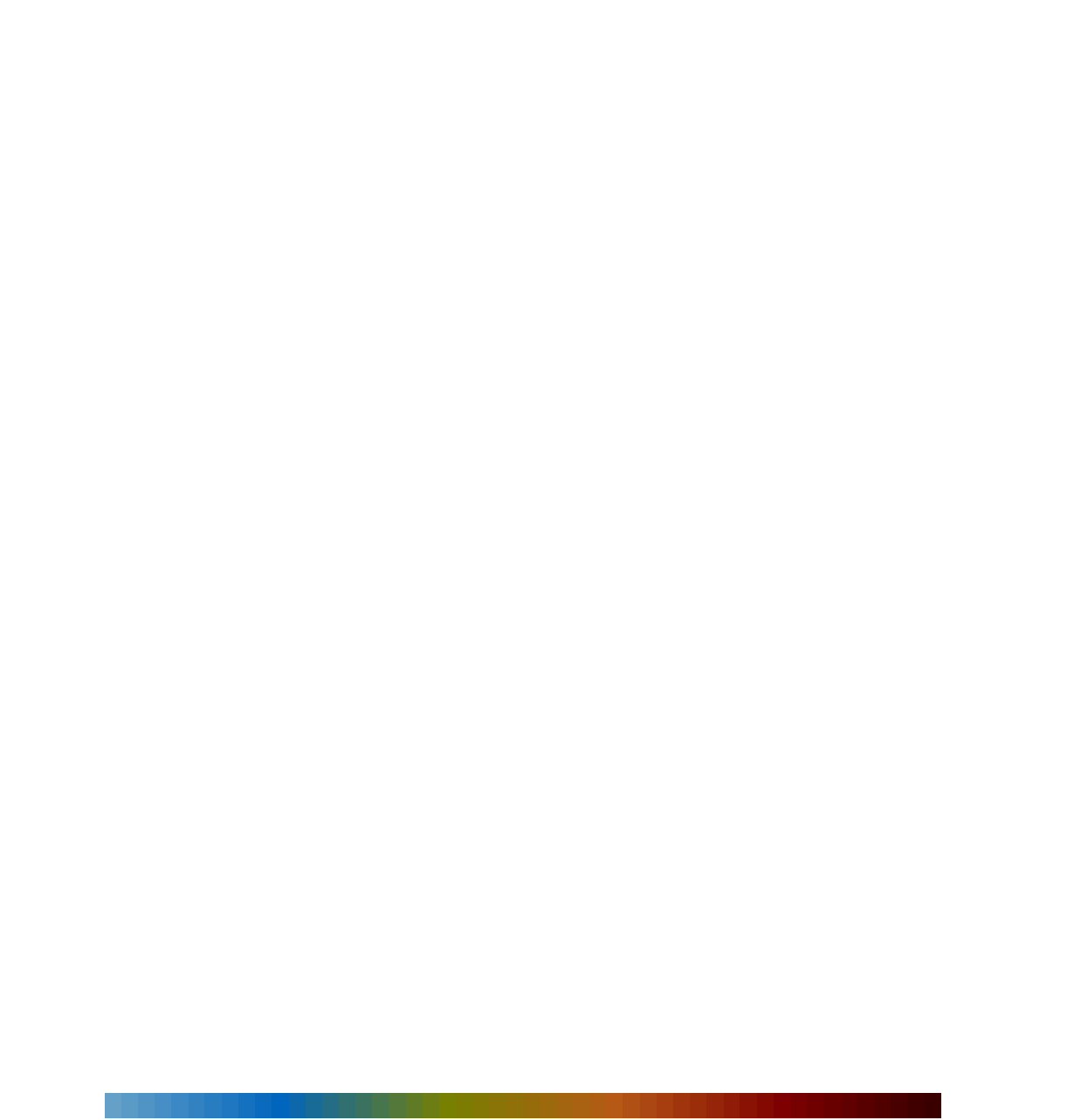
    \caption{Comparison of ERRW and RWRE}
    \parbox{\textwidth}{All trees are Galton-Watson trees with $p_1 = 0.7, p_2 = 0.15, p_3 = 0.1, p_4 = 0.05$.
    The ERRW was simulated for $100000$ steps or until the boundary of the
    visible part was reached, but not directly: the values of the variables $A_v$ of the corresponding
    RWRE were randomly selected according to the distributions which arise from the
    independent urn construction, and then the random walk was run on the resulting MC.
    On the left, the weights $w_n$ of the edges at the end of the simulations are color coded;
    on the right, $\log\left(1 + 100A_v\right)$ is color coded.}
    \label{fig:errw_rwre}
  \end{center}
\end{figure}

\begin{example}
  Consider \autoref{fig:errw_rwre} on page \pageref{fig:errw_rwre}. On the right, it shows an
  instantiation of the environment of the RWRE which corresponds to the the ERRW
  with different $\Delta$ values. The example demonstrates that for large values of $\Delta$, the
  the $A_v$ mostly take values at the extreme ends of the value range, while the $A_v$ cluster more
  around an average value for smaller values of $\Delta$.
  
  The intuition behind this observation is
  simple: for large $\Delta$ values, the urns often end up with a distribution of balls which
  is biased very strongly towards one color since after only a single ball of that color was drawn,
  many other balls of the same color are added. On the other hand, for small $\Delta$ values, the
  proportion of balls of a certain color only increases by little if a ball of that color is drawn.
  As a result, the distribution of balls is more balanced.
\end{example}

\begin{theorem}[ERRW on Trees]
  \label{thm:errw}
  Let $T$ be a an infinite tree. Consider the ERRW on
  $T$ with reinforcement parameter $\gls{reinforcement} > 0$. Let $\Delta_0$ be the solution to the equation
  \begin{align*}
    \frac{\Gamma\left(\frac{2 + \Delta}{4\Delta}\right)^2}{\Gamma\left(\frac{1}{2\Delta}\right)\Gamma\left(\frac{1 + \Delta}{2\Delta}\right)} = \frac{1}{\gls{br}}
  \end{align*}
  Then the ERRW on $T$ is
  \begin{enumerate}[(i)]
    \item \label{thm:errw_trans} a.s.~transient if $\Delta < \Delta_0$
    \item \label{thm:errw_rec} a.s.~recurrent if $\Delta > \Delta_0$.
  \end{enumerate}
  The expected return time to the root $\rho$ is always infinite.
\end{theorem}

\begin{tproof}
  The argumentation in the proof of \autoref{thm:errwgw} shows that a RWRE equivalent
  to the ERRW can be constructed using the tools (the independent urns) from above, and that
  for the RWRE, the $p$ from \autoref{thm:transrecrwre} is given by \Cref{equ:p_ex_at}.
  The assumptions of \autoref{thm:transrecrwre} are also shown to hold in the proof of \autoref{thm:errwgw}.
  Hence, \autoref{thm:transrecrwre} \ref{thm:rwre_trans} and \ref{thm:rwre_rec} are applicable and prove
  \autoref{thm:errw} \ref{thm:errw_trans} and \ref{thm:errw_rec}, respectively.

  The expected return time is proven to be infinite at the beginning of \cite[Section 2 on page 1230]{errwpemantle}.
\end{tproof}

\subsection{ERRW on Galton-Watson Trees}
\label{ssec:errwgw}

\begin{theorem}[ERRW on Galton-Watson Trees]
  \label{thm:errwgw}
  Let $T$ be a Galton-Watson tree with mean $\gls{gwmean} > 1$. Consider the ERRW on
  $T$ with reinforcement parameter $\gls{reinforcement} > 0$. Let $\Delta_0$ be the solution to the equation
  \begin{align*}
    \frac{\Gamma\left(\frac{2 + \Delta}{4\Delta}\right)^2}{\Gamma\left(\frac{1}{2\Delta}\right)\Gamma\left(\frac{1 + \Delta}{2\Delta}\right)} = \frac{1}{m}
  \end{align*}
  Then, conditional upon non-extinction of $T$, the ERRW on $T$ is
  \begin{enumerate}[(i)]
    \item \label{thm:errwgw_trans} a.s.~transient if $\Delta < \Delta_0$
    \item \label{thm:errwgw_rec} a.s.~recurrent if $\Delta \geq \Delta_0$.
  \end{enumerate}
  The expected return time to the root $\rho$ is always infinite.
\end{theorem}

\begin{remark}
  In the case $\Delta > \Delta_0$, the RWRE equivalent to the ERRW is
  a.s.~positive recurrent. This case is not mentioned in \autoref{thm:errwgw} since the
  a.s.~positive recurrence does not carry over to the ERRW, at least not in the
  usual sense: the expected return time to the root is always infinite.
\end{remark}

\begin{tproof}
  The proof is based on \autoref{thm:rwregw}. In order to use \autoref{thm:rwregw}, a RWRE
  with the same distribution as the ERRW is necessary. The construction was described in
  the introductory paragraphs of \autoref{sec:transrecerrw} and it specifies the distribution of
  the random environment $\omega$ of the RWRE which is equivalent to the ERRW.
  \cite[Lemma 1]{errwpemantle} directly gives the distribution of the environment, but \autoref{thm:rwregw}
  relies on the random variable $\glslink{av}{A}$, which can be recovered from $\omega$ as in the discussion
  following \autoref{def:rwre}: $A_{\gls{child}} = \glslink{omuv}{\omega\left(v, v_i\right)} \cdot \omega\left(v, \gls{par}\right)^{-1}$.

  By the construction with the independent urns, the random variables $A_v$ fulfil the independence
  conditions of \autoref{thm:rwregw}. In addition, they are also identically distributed for
  $v \in V \setminus \left(\left\{\gls{root}\right\} \cup \glslink{distset}{T_1}\right)$. To see
  this, note that for $v \neq \rho$, $A_{v_i}$ is represented as the ratio of the probability to go
  from $v$ to $v_i$ divided by the probability to go from $v$ to $\overleftarrow{v}$. These
  transitions are simulated by the independent urns in the RWRE.
  
  Now, consider another node $u \neq \rho$ and $A_{u_j}$. If $u$ and $v$ have the same number of
  children, it is immediately clear that $A_{v_i}$ and $A_{u_j}$ are identically distributed since
  the same type of urn is used at both $u$ and $v$. But if the number of children is different,
  the distributions are still identical: the considered probability ratio only depends on
  how many balls of the colors corresponding to the respective edges are currently contained in the
  urn. All other colors are irrelevant for this ratio, so it is in particular irrelevant how many
  other colors are used, i.e.~how many edges are incident to $u$ and $v$.

  It remains to calculate the expectation of $A^t$. As the number of children is not relevant,
  it suffices to consider an urn with balls of only two different colors, and calculate the
  expectation of the probability ratio using the distribution given in \cite[Lemma 1]{errwpemantle}.
  This is done in \autoref{sec:errw_z}. Indeed, on \cite[page 25]{rrwrenlund} the expectation is
  also calculated for a variable number of children and seen to be the same for each number of
  children. Both calculations yield the result:
  \begin{align*}
    \bEx{A^t} = \frac{\Gamma\left(\frac{1}{2\Delta} + t\right)\Gamma\left(\frac{1 + \Delta}{2\Delta} - t\right)}{\Gamma\left(\frac{1}{2\Delta}\right)\Gamma\left(\frac{1 + \Delta}{2\Delta}\right)}
    \textrm{ for }-\frac{1}{2\Delta} < t < \frac{1 + \Delta}{2\Delta}
    \textrm{, otherwise }\infty
  \end{align*}
  The $p = \min_{0 \leq t \leq 1} \bEx{A^t}$ is also calculated in \autoref{sec:errw_z}, in
  \Cref{equ:edelta}. Depending on the reinforcement paramater $\Delta$, it holds that
  \begin{align}
    \label{equ:p_ex_at}
    p = \min_{0 \leq t \leq 1} \bEx{A^t} = E\left(\Delta\right) := \frac{\Gamma\left(\frac{2 + \Delta}{4\Delta}\right)^2}{\Gamma\left(\frac{1}{2\Delta}\right)\Gamma\left(\frac{1 + \Delta}{2\Delta}\right)}
  \end{align}
  $E\left(\Delta\right)$ was already analyzed on \cite[page 27]{rrwrenlund}, and is later analyzed
  in \autoref{sec:errw_z} after \Cref{equ:edelta}. It is a strictly decreasing function with
  $E\left(\Delta\right) \overset{\Delta \to 0}{\longrightarrow} 1$ and $E\left(\Delta\right) \overset{\Delta \to \infty}{\longrightarrow} 0$.
  Also see
  \autoref{fig:e_delta} on page \pageref{fig:e_delta}
  for a plot.

  By \autoref{thm:rwregw} \ref{thm:rwregw_trans}, the RWRE (and therefore the
  equivalent ERRW) is a.s.~transient for $pm > 1$, i.e.~for
  \begin{align*}
    E\left(\Delta\right) = \frac{\Gamma\left(\frac{2 + \Delta}{4\Delta}\right)^2}{\Gamma\left(\frac{1}{2\Delta}\right)\Gamma\left(\frac{1 + \Delta}{2\Delta}\right)} > \frac{1}{m} \in \left(0,1\right)
  \end{align*}
  By \autoref{thm:rwregw} \ref{thm:rwregw_rec}, the RWRE is a.s.~recurrent
  for $E\left(\Delta\right) \leq \frac{1}{m}$. Since $E$ is decreasing, the value $\Delta_0$, defined
  to be the solution of the equation $E\left(\Delta\right) = \frac{1}{m}$, is a critical threshold
  value: for $\Delta < \Delta_0$, the random walk is transient, and for $\Delta \geq \Delta_0$, it
  is recurrent.

  The expected return time is proven to be infinite at the beginning of \cite[Section 2 on page 1230]{errwpemantle},
  and in \autoref{ssec:zexpret} for the special case $\Delta \geq 1$ for a Galton-Watson tree with
  exactly one child at every node (i.e.~the graph $\bbN$).
\end{tproof}

\subsection{The Critical $\Delta$ Value for Transience and Recurrence}
\label{ssec:critdel}

The critical $\glslink{reinforcement}{\Delta_0}$ from \autoref{thm:errwgw} be computed for
Galton-Watson trees with different means $\gls{gwmean} > 1$ (for $\Delta < \Delta_0$, the
ERRW is a.s.~transient, for $\Delta \geq \Delta_0$ a.s.~recurrent):
\begin{table}[H]
  \begin{center}
    \def\arraystretch{1.1}
    \begin{tabular}{c|cccccccccccc}
      $m =$ & $1.1$ & $1.18$ & $1.3$ & $1.5$ & $1.7$ & $2$ & $2.3$ & $2.5$ & $2.7$ & $3$ & $4$ & $5$ \\ \hline
      $\Delta_0\approx$ & $0.61$ & $1.00$ & $1.53$ & $2.35$ & $3.14$ & $4.29$ & $5.43$ & $6.19$ & $6.94$ & $8.07$ & $11.80$ & $15.52$
    \end{tabular}\\[0.5em]
    \begin{tabular}{c|cccccccccc}
      $m =$ & $6$ & $7$ & $8$ & $9$ & $10$ & $10^2$ & $10^3$ & $10^4$ & $10^5$ \\ \hline
      $\Delta_0\approx$ & $19.24$ & $22.97$ & $26.67$ & $30.38$ & $34.09$ & $367.85$ & $3705.19$ & $37078.53$ & $370811.97$
    \end{tabular}
    \caption{Critical $\Delta$ values for different Galton-Watson trees}
    \label{tab:critdel}
  \end{center}
\end{table}
\noindent where the values for $m \geq 4$ are taken from \cite{rrwrenlund}. The values are also
plotted in the following figure.

\begin{figure}[H]
  \begin{center}
    \begin{tikzpicture}
      \begin{axis}[
        xlabel={Mean $m$ of Galton-Watson tree},
        ylabel={Critical value $\Delta_0$},
        grid=major
      ]
        \addplot[color=tumOrange,dashed] coordinates {
          (1, 0.66) (4, 11.85)
        } node[below=1.3cm, left=2.1cm, rotate=40] {straight line for comparison};
        \addplot[color=tumBlue,mark=*] coordinates {
          (1.1, 0.61) (1.18, 1) (1.3, 1.53) (1.5, 2.35) (1.7, 3.14) (2, 4.29) (2.3, 5.43)
          (2.5, 6.19) (2.7, 6.94) (3, 8.07) (4, 11.8)
        };
      \end{axis}
    \end{tikzpicture}
    \begin{tikzpicture}
      \begin{axis}[
        xlabel={Mean $m$ of Galton-Watson tree (logarithmic scale)},
        ylabel={Critical value $\Delta_0$ (logarithmic scale)},
        grid=major,
        xmode=log,
        ymode=log
      ]
      \addplot[color=tumBlue,mark=*] coordinates {
        (1.1, 0.61) (1.18, 1) (1.3, 1.53) (1.5, 2.35) (1.7, 3.14) (2, 4.29) (2.3, 5.43)
        (2.5, 6.19) (2.7, 6.94) (3, 8.07) (4, 11.8) (5, 15.52) (6, 19.24) (7, 22.97) (8, 26.67)
        (9, 30.38) (10, 34.09) (100, 367.85) (1000, 3705.19) (10000, 37078.53) (100000, 370811.97)
      };
      \end{axis}
    \end{tikzpicture}
    \caption{Plot of the critical $\Delta$ values in dependence of the mean}
    \label{fig:critdel}
  \end{center}
\end{figure}
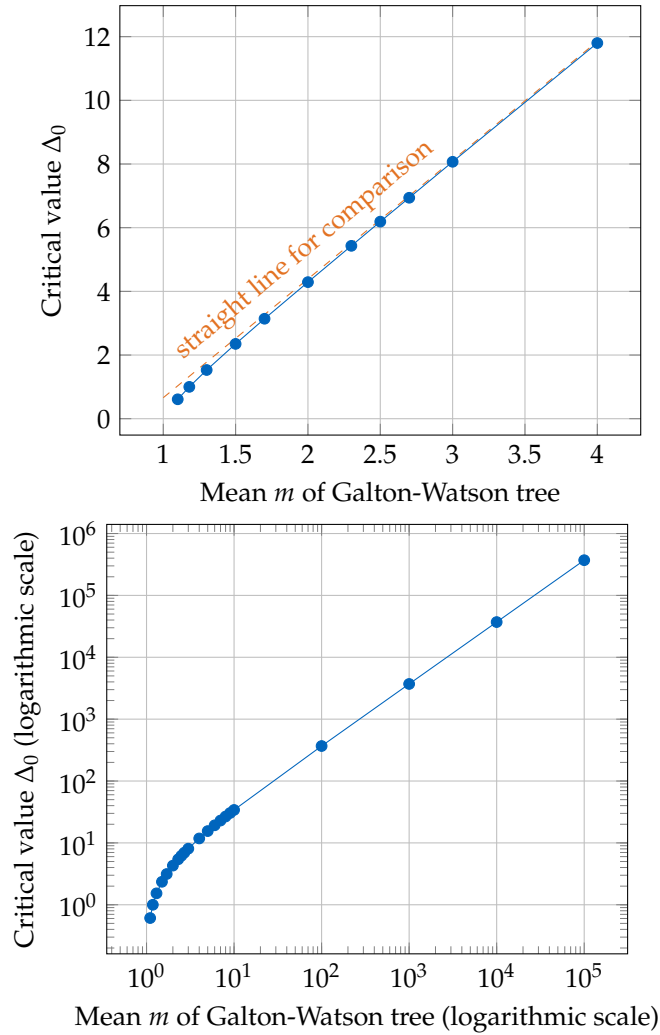

\begin{example}
  Consider \autoref{fig:errw_rwre} on page \pageref{fig:errw_rwre}. On the left, a path taken by
  the ERRW on the trees is shown. The plotted Galton-Watson
  trees have mean $1.5$, hence the critical value $\Delta_0$ is approximately $2.35$. For
  smaller values of $\Delta$, such as depicted in the top row, the ERRW is transient.
  For larger values of $\Delta$, such as depicted in the bottom row, the ERRW is recurrent.
  Indeed, the random walk on the top leaves the visible area, while the random walk on the bottom
  stays close to the root for the simulated time period.
\end{example}

\clearpage
\newpage

\section{ERRW on $\bbZ$}
\label{sec:errw_z}

In this section, the ERRW with reinforcement parameter $\gls{reinforcement}$ on $\bbZ$ is considered,
where $\bbZ$ is interpreted as a tree with root $\rho = 0$ where the root has two children, $1$ and
$-1$, and every other vertex $0 \neq z \in \bbZ$ has one child, $z + 1$ if $z > 0$ and $z - 1$
otherwise. The distribution of the possible paths taken by the ERRW is equal to the
distribution of the RWRE with the following environment $\omega$, where
$\omega_x := \glslink{omuv}{\omega\left(x, x + 1\right)}$ for $x \in \bbZ$ is used for conciseness:

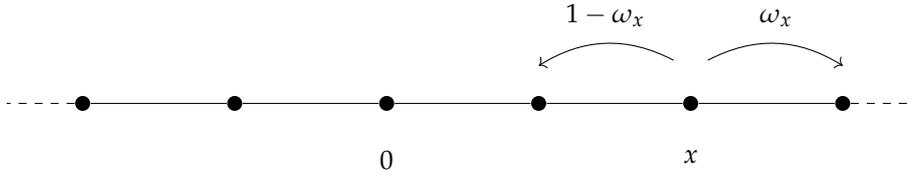
\begin{figure}[H]
  \begin{center}
    \begin{tikzpicture}
      \node[circle,fill=black,inner sep=0.7mm] (N0) at (0, 0) {};
      \node[circle,fill=black,inner sep=0.7mm] (N1) at (2, 0) {};
      \node[circle,fill=black,inner sep=0.7mm] (N2) at (4, 0) {};
      \node[circle,fill=black,inner sep=0.7mm] (N3) at (6, 0) {};
      \node[circle,fill=black,inner sep=0.7mm] (Nm1) at (-2, 0) {};
      \node[circle,fill=black,inner sep=0.7mm] (Nm2) at (-4, 0) {};
      \node[below=5mm] at (N0) {$0$};
      \node[below=5mm] at (N2) {$x$};
      \draw (Nm2) -- (Nm1) -- (N0) -- (N1) -- (N2) -- (N3);
      \draw[dashed] (N3) -- (7, 0);
      \draw[dashed] (Nm2) -- (-5, 0);
      \node (aN2l) at (3.9, 0.5) {};
      \node (aN2r) at (4.1, 0.5) {};
      \draw[->] (aN2l) edge[bend right] node[above=1mm] {$1 - \omega_x$} (2, 0.5);
      \draw[->] (aN2r) edge[bend left] node[above=1mm] {$\omega_x$} (6, 0.5);
    \end{tikzpicture}
    \caption{The ERRW on $\bbZ$}
    \label{fig:errw_z}
  \end{center}
\end{figure}

Using the construction described in \autoref{sec:transrecerrw} and \cite[Lemma 1]{errwpemantle}, one
sees that the $\omega_x$ are Beta-distributed:
\begin{align*}
  \textrm{for }x > 0\textrm{: }& \qquad \omega_x \sim \glslink{beta}{B\left(\frac{1}{2\Delta}, \frac{1+\Delta}{2\Delta}\right)} \\
  \textrm{for }x = 0\textrm{: }& \qquad \omega_x \sim B\left(\frac{1}{2\Delta}, \frac{1}{2\Delta}\right) \\
  \textrm{for }x < 0\textrm{: }& \qquad \omega_x \sim B\left(\frac{1+\Delta}{2\Delta}, \frac{1}{2\Delta}\right)
\end{align*}
i.e., $\omega_x \sim B\left(\frac{1}{2}, 1\right)$ for $\Delta = 1$ and $\omega_x \sim B\left(1, \frac{3}{2}\right)$
for $\Delta = 0.5$, for example (for $x > 0$). The random variables $\omega_x, x \in \bbZ$ are independent.

\begin{figure}
  \begin{center}
    \begin{tabular}{cc}
      ~~~~~~~~~~~~~ $\Delta = 0.01$ & ~~~~~~~~~~~~~ $\Delta = 0.1$ \\
      \input{simulations/cdf-d0.01.tex} & \input{simulations/cdf-d0.10.tex} \\
      & \\
      ~~~~~~~~~~~~~ $\Delta = 0.2$ & ~~~~~~~~~~~~~ $\Delta = 0.5$ \\
      \input{simulations/cdf-d0.20.tex} & \input{simulations/cdf-d0.50.tex} \\
      & \\
      ~~~~~~~~~~~~~ $\Delta = 1.0$ & ~~~~~~~~~~~~~ $\Delta = 1.5$ \\
      \input{simulations/cdf-d1.00.tex} & \input{simulations/cdf-d1.50.tex}
    \end{tabular}
    \caption{Cumulative distribution function for the $A_x$}
    \label{fig:cdf}
  \end{center}
\end{figure}
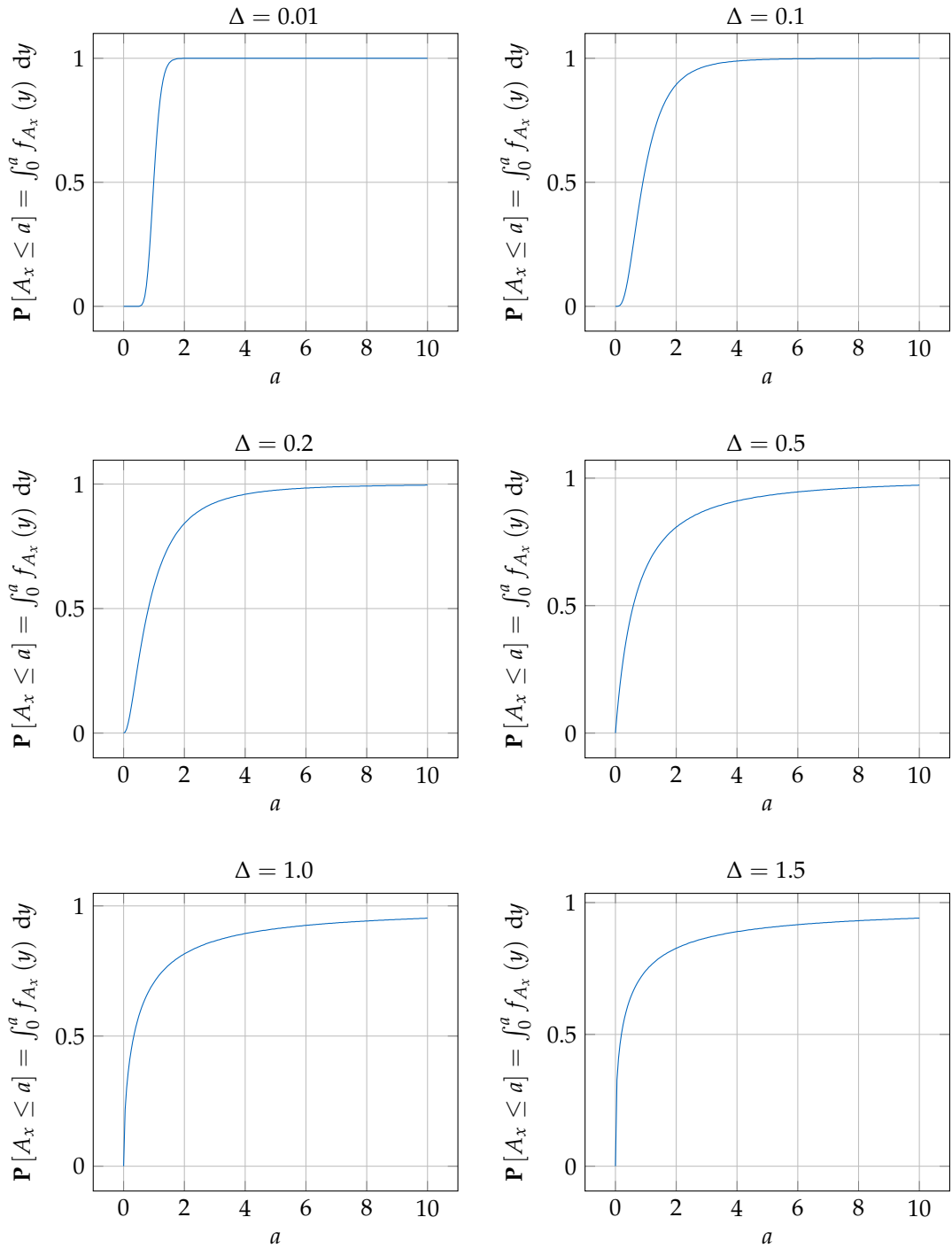

With the notation from \autoref{def:rwre}, $\glslink{av}{A_x} = \frac{\omega_{x-1}}{1 - \omega_{x-1}}$ for $x > 1$,
$A_x = \frac{1 - \omega_{x+1}}{\omega_{x+1}}$ for $x < -1$. The $A_x$ are iid for
$x \in \bbZ \setminus \left\{-1,0,1\right\}$ and their density on $\left[0, \infty\right)$ is
given by
\begin{align*}
  f_{A_x} : \left[0, \infty\right) \to \left[0, \infty\right), f_{A_x}\left(y\right) =
  \frac{\Gamma\left(\frac{2 + \Delta}{2\Delta}\right)}{\Gamma\left(\frac{1}{2\Delta}\right)\Gamma\left(\frac{1 + \Delta}{2\Delta}\right)} \cdot
  y^{\frac{1}{2\Delta} - 1} \cdot \left(1 + y\right)^{-\frac{2 + \Delta}{2\Delta}}
\end{align*}
where the density can easily be calculated using a change of variables and the density of the Beta
distribution. However, it is usually easier to express any expectations using the Beta distribution
directly, so this expression for the density will no longer be used and is only a side note. For
$x \in \left\{-1, 1\right\}$, it is possible to choose $A_1 = \omega_0$ and $A_{-1} = 1 - \omega_0$.
\autoref{fig:cdf} on page \pageref{fig:cdf} shows the cumulative distribution function of $f_{A_x}$.

Now, for $x > 0$ and $t \in \bbR$, it holds that ($\bEx{\cdot}$ is the expectation under the probability measure $\gls{ompmeas}$,
which defines the distribution of $\omega$ as in \autoref{def:rwre})
\begin{align*}
  \bEx{A_{x+1}^t} &= \bEx{\left(\frac{\omega_x}{1 - \omega_x}\right)^t}
  = \frac{\Gamma\left(\frac{2 + \Delta}{2\Delta}\right)}{\Gamma\left(\frac{1}{2\Delta}\right)\Gamma\left(\frac{1 + \Delta}{2\Delta}\right)} \int_0^1 \frac{x^t}{\left(1 - x\right)^t} x^{\frac{1}{2\Delta} - 1} \left(1 - x\right)^{\frac{1+\Delta}{2\Delta} - 1} \dx{x} \\
  &= \frac{\Gamma\left(\frac{2 + \Delta}{2\Delta}\right)}{\Gamma\left(\frac{1}{2\Delta}\right)\Gamma\left(\frac{1 + \Delta}{2\Delta}\right)} \int_0^1 x^{\frac{1}{2\Delta}+t-1} \left(1 - x\right)^{\frac{1+\Delta}{2\Delta}-t-1} \dx{x} \\
  &\overset{\textrm{for }-\frac{1}{2\Delta} < t < \frac{1 + \Delta}{2\Delta}}{=}\;\; \frac{\Gamma\left(\frac{2 + \Delta}{2\Delta}\right)}{\Gamma\left(\frac{1}{2\Delta}\right)\Gamma\left(\frac{1 + \Delta}{2\Delta}\right)} \cdot
  \frac{\Gamma\left(\frac{1}{2\Delta} + t\right)\Gamma\left(\frac{1 + \Delta}{2\Delta} - t\right)}{\Gamma\left(\frac{2 + \Delta}{2\Delta}\right)}\\
  &= \frac{\Gamma\left(\frac{1}{2\Delta} + t\right)\Gamma\left(\frac{1 + \Delta}{2\Delta} - t\right)}{\Gamma\left(\frac{1}{2\Delta}\right)\Gamma\left(\frac{1 + \Delta}{2\Delta}\right)}
  \;\;\overset{\textrm{for }t = 1}{=}\;\; \frac{\Gamma\left(\frac{1}{2\Delta} + 1\right)\Gamma\left(\frac{1 - \Delta}{2\Delta}\right)}{\Gamma\left(\frac{1}{2\Delta}\right)\Gamma\left(\frac{1 - \Delta}{2\Delta} + 1\right)}
  = \frac{1}{2\Delta} \cdot \frac{2\Delta}{1 - \Delta} = \frac{1}{1 - \Delta}
\end{align*}
for $t \geq \frac{1 + \Delta}{2\Delta}$ and $t \leq -\frac{1}{2\Delta}$, the integral is
$\infty$. Setting $t = 1$ as in the last line above gives $\bEx{A_{x+1}} = \frac{1}{1 - \Delta}$,
where the condition for the integral to be $\infty$ simplifies to $\Delta \geq 1$. The expectation is
always strictly larger than $1$ in this case.

In order to analyze the stationary distribution, it will be useful to see for which $t$
the expectation $\bEx{A_{x+1}^t}$ attains the minimal value. With
$K_\Delta := \left(\Gamma\left(\frac{1}{2\Delta}\right)\Gamma\left(\frac{1 + \Delta}{2\Delta}\right)\right)^{-1}$
and assuming $-\frac{1}{2\Delta} < t < \frac{1 + \Delta}{2\Delta}$, it holds that
$f\left(t\right) := \bEx{A_{x+1}^t} = K_\Delta \Gamma\left(\frac{1}{2\Delta} + t\right)\Gamma\left(\frac{1 + \Delta}{2\Delta} - t\right)$.
Note that $f$ is symmetric about $t_0 = \frac{1}{2} \cdot \left(\frac{1 + \Delta}{2\Delta} - \frac{1}{2\Delta}\right) = \frac{1}{4}$,
where $\frac{1}{2\Delta} + t_0 = \frac{1 + \Delta}{2\Delta} - t_0$.
$t_0$ is therefore a likely candidate for the point where $f$ attains its minimum.

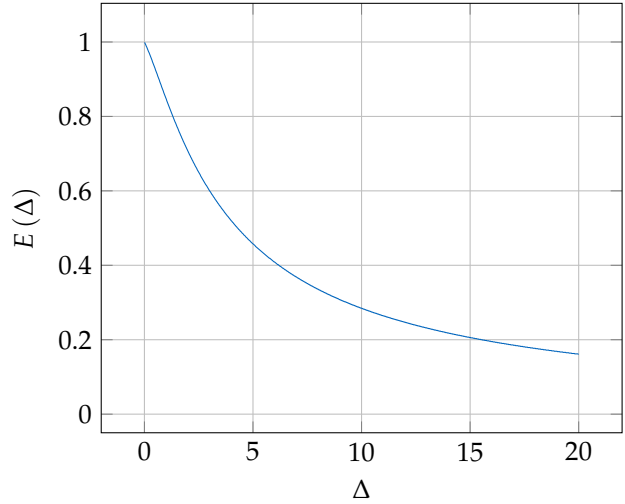
\begin{figure}
  \begin{center}
    \input{simulations/e_delta.tex}
    \caption{Plot of $E\left(\Delta\right)$}
    \label{fig:e_delta}
  \end{center}
\end{figure}

Indeed, with $\psi$ being the digamma function, the derivative of $f$ can be expressed as:
\begin{align*}
  f'\left(t\right) &= \underbrace{K_\Delta \Gamma\left(\frac{1}{2\Delta} + t\right)\Gamma\left(\frac{1 + \Delta}{2\Delta} - t\right)}_{> 0}
  \cdot \underbrace{\left( \psi\left(\frac{1}{2\Delta} + t\right) - \psi\left(\frac{1 + \Delta}{2\Delta} - t\right) \right)}_{\begin{cases}
    = 0 & \textrm{ for }t = t_0 \\
    > 0 & \textrm{ for }\frac{1 + \Delta}{2\Delta} > t > t_0 \\
    < 0 & \textrm{ for }-\frac{1}{2\Delta} < t < t_0
  \end{cases}}
\end{align*}
since $\psi$ is increasing on $\bbR_+$. Hence, $t_0 = \frac{1}{4}$ is the (strict) global
minimum of $\bEx{A_{x+1}^t}$ with
\begin{align}
  \label{equ:edelta}
  \bEx{A_{x+1}^\frac{1}{4}} &= \bEx{\left(\frac{\omega_x}{1 - \omega_x}\right)^\frac{1}{4}}
  = \frac{\Gamma\left(\frac{2 + \Delta}{4\Delta}\right)^2}{\Gamma\left(\frac{1}{2\Delta}\right)\Gamma\left(\frac{1 + \Delta}{2\Delta}\right)} =: E\left(\Delta\right)
\end{align}
It is easy to check that $E\left(\Delta\right) \overset{\Delta \to 0}{\longrightarrow} 1$,
$E\left(\Delta\right) \overset{\Delta \to \infty}{\longrightarrow} 0$ and
that $E$ is strictly decreasing for $\Delta > 0$. For a plot of $E\left(\Delta\right)$, see
\autoref{fig:e_delta} on page \pageref{fig:e_delta}.
In particular, $E\left(\Delta\right) < 1$ for
$\Delta > 0$, which directly follows from the following fact:
$E\left(\Delta\right) = f\left(\frac{1}{4}\right) = K_\Delta \Gamma\left(\frac{1}{2\Delta} + \frac{1}{4}\right)\Gamma\left(\frac{1 + \Delta}{2\Delta} - \frac{1}{4}\right)$
and $f\left(0\right) = K_\Delta \cdot K_\Delta^{-1} = 1$; as $t_0 = \frac{1}{4}$ is the
(strict) global minimum of $f$, it follows that
$E\left(\Delta\right) = f\left(t_0\right) < f\left(0\right) = 1$.

By \autoref{thm:transrecrwre} \ref{thm:rwre_posrec}, the considered RWRE is
a.s.~positive recurrent ($p = E\left(\Delta\right) < 1$ for $\Delta > 0$ and
$\varlimsup_{n \to \infty} \left|\gls{distset}\right|^{\frac{1}{n}} = \varlimsup_{n \to \infty} \sqrt[n]{2} = 1$).
It is quite intuitive that the RWRE equivalent to the ERRW on $\bbZ$ is
recurrent since the simple random walk on $\bbZ$ is already recurrent, and adding a reinforcement
parameter which favors traversing already visited edges makes the random walk even more recurrent.

\subsection{Stationary Distribution}
\label{ssec:zstatdistr}

Since the correspondent RWRE is a.s.~positive recurrent (see last paragraph),
a stationary distribution exists for almost all MCs.
Denote by $\gls{statmeas} : \bbZ \to \left[0, \infty\right)$ a stationary measure (not
necessarily a distribution) for the MC given by the environment $\omega$. Setting
$\mu_\omega\left(0\right) = 1$, the unique stationary measure with that property is given by
\begin{align*}
  \textrm{for }x > 0\textrm{: }& \qquad \mu_{\omega}\left(x\right) = \glslink{cv}{C_x} + C_{x+1} = \left(1 + A_{x+1}\right) \cdot \prod_{i=1}^x A_x = \frac{\omega_0}{1 - \omega_x} \cdot \prod_{i=1}^{x-1} \frac{\omega_i}{1 - \omega_i}\\
  \textrm{for }x < 0\textrm{: }& \qquad \mu_{\omega}\left(x\right) = C_x + C_{x-1} = \left(1 + A_{x-1}\right) \cdot \prod_{i=x}^{-1} A_x = \frac{1 - \omega_0}{\omega_x} \cdot \prod_{i=x+1}^{-1} \frac{1 - \omega_i}{\omega_i}
\end{align*}

Now, consider, for $x > 0$
\begin{align*}
  \bEx{\mu_{\omega}\left(x\right)} &= \bEx{\frac{\omega_0}{1 - \omega_x} \cdot \prod_{i=1}^{x-1} \frac{\omega_i}{1 - \omega_i}}
  \overset{\textrm{independence}}{=} \bEx{\omega_0} \cdot \bEx{\frac{1}{1 - \omega_x}} \cdot \prod_{i=1}^{x-1} \bEx{\frac{\omega_i}{1 - \omega_i}} \\
  &= \bEx{\omega_0} \cdot \bEx{\frac{1}{1 - \omega_x}} \cdot \begin{cases}
    \left(\frac{1}{1 - \Delta}\right)^{x-1} & \textrm{ if } \Delta < 1 \\
    \infty & \textrm{ if } \Delta \geq 1
  \end{cases}
  \qquad \overset{x \to \infty}{\longrightarrow} \infty
\end{align*}
Note that $\bEx{\frac{1}{1 - \omega_x}}$ is constant for $x > 0$ since the $\omega_x$ are
iid for $x > 0$.

Since $\mu_\omega$ is not a probability measure in general, it is necessary to divide by the total
measure to obtain a stationary distribution. Call the total measure $Z_\omega$, i.e.
\begin{align*}
  Z_\omega := \sum_{x \in \bbZ} \mu_{\omega}\left(x\right)
\end{align*}
As the expectation of $\mu_{\omega}\left(x\right)$ tends to $\infty$ for $x \to \infty$, it holds
that $\bEx{Z_\omega} = \infty$ even though almost all MCs do have a stationary measure
in the environment $\omega$ which satisfies $\mu_\omega\left(0\right) = 1$ and which is, of course,
finite (because of the a.s.~positive recurrence), i.e.~$\bPrb{Z_\omega < \infty} = 1$.

Dividing by $Z_\omega$ and mixing the respective stationary distributions according to the law of
$\omega$ yields a mixed stationary distribution:
\begin{align}
  \label{equ:zstat}
  \mu\left(x\right) := \int \mu_\omega\left(x\right) \cdot Z_\omega^{-1} \mudx{\bP}{\omega} \qquad\qquad \left(x \in \bbZ\right)
\end{align}
As can easily be seen with the dominated convergence theorem, $\mu$ is a probability measure on
$\bbZ$, which can be interpreted as a kind of stationary distribution for the ERRW on $\bbZ$:
\begin{equation}
  \label{equ:zlongterm}
  \begin{split}
    &\lim\limits_{m \to \infty} \frac{1}{m + 1} \sum_{n = 0}^m \Prb{X_n = x}\\
    &\quad= \lim\limits_{m \to \infty} \frac{1}{m + 1} \sum_{n = 0}^m \int \Prbc{\omega}{X_n = x} \mudx{\bP}{\omega}
    = \lim\limits_{m \to \infty} \int \frac{1}{m + 1} \sum_{n = 0}^m \Prbc{\omega}{X_n = x} \mudx{\bP}{\omega}\\
    &\quad\overset{\textrm{dom. conv.}}{=} \int \lim\limits_{m \to \infty} \frac{1}{m + 1} \sum_{n = 0}^m \Prbc{\omega}{X_n = x} \mudx{\bP}{\omega}\\
    &\quad\overset{\textrm{MC given by }\omega\textrm{ a.s.~positive rec.}}{=}\;\; \int \mu_\omega\left(x\right) \cdot Z_\omega^{-1} \mudx{\bP}{\omega} = \mu\left(x\right)
  \end{split}
\end{equation}
Due to the fact that the random walk can only be at even numbers at even time steps (the
MCs given by the environments $\omega$ are all periodic), it is necessary to take the
average of the probabilities to be in $x$ from time $0$ to time $m$, instead of directly taking
the limit of the probability to be in $x$ at time $n$ (which does not exist).

To analyze $\mu$ further, and to circumvent the problem of the infinite expectation of the total
measure of the measures $\mu_\omega$, \Cref{equ:edelta} is useful. Consider, for $x > 0$,
\begin{align*}
  \bEx{\mu_{\omega}\left(x\right)^{\frac{1}{4}}}
  &= \underbrace{\bEx{\omega_0^{\frac{1}{4}}} \cdot \bEx{\left(1 - \omega_x\right)^{-\frac{1}{4}}}}_{=: C_{1,\Delta}} \cdot \prod_{i=2}^{x} \bEx{A_x^{\frac{1}{4}}}
  = C_{1,\Delta} \cdot E\left(\Delta\right)^{x-1}
\end{align*}
Note that $C_{1,\Delta}$ does not depend on $x$ since the $\omega_x$ are iid for $x > 0$.
It is also easily verifiable that $C_{1,\Delta} < \infty$ for all $\Delta > 0$ ($\omega_0^{\frac{1}{4}} \leq 1$
and for $\left(1 - \omega_x\right)^{-\frac{1}{4}}$, a calculation as for $\bEx{A_x^t}$ as well
as noting that $\frac{1}{4} < \frac{1 + \Delta}{2\Delta}$ shows finiteness).
Hence, while the expectation of $\mu_{\omega}\left(x\right)$ tends to infinity, the expectation of
$\mu_{\omega}\left(x\right)^{\frac{1}{4}}$ decreases exponentially in $x$ (recall that
$E\left(\Delta\right) < 1$ for $\Delta > 0$).

As a next step, the probability that $\mu_\omega\left(x\right) \cdot Z_\omega^{-1}$ is greater than
$e^{-\alpha x}$ is analyzed for some $\alpha > 0$.
\begin{align*}
  \bPrb{\mu_\omega\left(x\right) \cdot Z_\omega^{-1} \geq e^{-\alpha x}}
  &= \bPrb{\left(\mu_\omega\left(x\right) \cdot Z_\omega^{-1}\right)^{\frac{1}{8}} \geq e^{-\frac{\alpha}{8} x}}
  \;\; \overset{\textrm{Markov inequ.}}{\leq} \;\;
  \bEx{ \mu_\omega\left(x\right)^\frac{1}{8} \cdot Z_\omega^{-\frac{1}{8}} } \cdot e^{\frac{\alpha}{8} x} \\
  &\hskip-1.3cm \overset{\textrm{Cauchy-Schwarz inequ.}}{\leq} \;\;
  \bEx{ \mu_\omega\left(x\right)^\frac{1}{4} }^{\frac{1}{2}} \cdot \bEx{ Z_\omega^{-\frac{1}{4}} }^{\frac{1}{2}} \cdot e^{\frac{\alpha}{8} x}
  \leq C_{1,\Delta}^{\frac{1}{2}} \cdot E\left(\Delta\right)^{\frac{x - 1}{2}} \cdot e^{\frac{\alpha}{8} x}
\end{align*}
Here, the inequality derived for the expectation of $\mu_{\omega}\left(x\right)^{\frac{1}{4}}$ above
is used. Also note that $Z_\omega \geq 1$ since $\mu_\omega\left(0\right) = 1$ and hence
$\bEx{ Z_\omega^{-\frac{1}{4}} } \leq 1$. The newly derived inequality can be reformulated as
\begin{align*}
  \bPrb{\mu_\omega\left(x\right) \cdot Z_\omega^{-1} \geq e^{-\alpha x}}
  &\leq \overbrace{\left(C_{1,\Delta} \cdot E\left(\Delta\right)^{-1}\right)^\frac{1}{2}}^{=: C_{2,\Delta}} \cdot E\left(\Delta\right)^{\frac{x}{2}} \cdot e^{\frac{\alpha}{8} x}
  &&= C_{2,\Delta} \cdot e^{\frac{x}{2} \log E\left(\Delta\right) + \frac{\alpha}{8} x}\\
  &&&= C_{2,\Delta} \cdot e^{-\left(\frac{\log\left(E\left(\Delta\right)^{-1}\right)}{2} - \frac{\alpha}{8}\right) x}
\end{align*}
Now, since $E\left(\Delta\right)^{-1} > 1$, it holds that $\frac{\log\left(E\left(\Delta\right)^{-1}\right)}{2} > 0$,
so choosing $0 < \alpha < \log\left(E\left(\Delta\right)^{-4}\right)$ implies
$\frac{\log\left(E\left(\Delta\right)^{-1}\right)}{2} - \frac{\alpha}{8} =: c_{\alpha,\Delta} > 0$.
The inequality then reads
\begin{align*}
  \bPrb{\mu_\omega\left(x\right) \cdot Z_\omega^{-1} \geq e^{-\alpha x}}
  &\leq C_{2,\Delta} \cdot e^{-c_{\alpha,\Delta} x} \qquad \textrm{where } c_{\alpha,\Delta} > 0
\end{align*}

With these tools, $\mu\left(x\right) = \bEx{\mu_\omega\left(x\right) \cdot Z_\omega^{-1}}$ can
be examined more closely.
\begin{align*}
  \mu\left(x\right) &= \bEx{\mu_\omega\left(x\right) \cdot Z_\omega^{-1}}
  = \bEx{\mu_\omega\left(x\right) Z_\omega^{-1} \cdot \mathbbm{1}_{\mu_\omega\left(x\right) Z_\omega^{-1} \geq e^{-\alpha x}}}
  + \bEx{\mu_\omega\left(x\right) Z_\omega^{-1} \cdot \mathbbm{1}_{\mu_\omega\left(x\right) Z_\omega^{-1} < e^{-\alpha x}}}\\
  &\leq \bPrb{\mu_\omega\left(x\right) \cdot Z_\omega^{-1} \geq e^{-\alpha x}} + e^{-\alpha x}
  \leq C_{2,\Delta} \cdot e^{-c_{\alpha,\Delta} x} + e^{-\alpha x} \qquad \textrm{where } \alpha, c_{\alpha,\Delta} > 0
\end{align*}
This final inequality shows that $\mu\left(x\right)$ decreases exponentially in $\left|x\right|$ (the inequality
was derived for $x > 0$, but it also holds for $x < 0$ by symmetry with $x$ replaced by
$\left|x\right|$ in the exponents on the right hand side).

Now, \Cref{equ:zlongterm} showed $\mu\left(x\right) = \lim\limits_{m \to \infty} \frac{1}{m + 1} \sum_{n = 0}^m \Prb{X_n = x}$
for the ERRW on $\bbZ$. $\mu$ therefore captures the long-term average distribution of
the ERRW on $\bbZ$. Consider an integer-valued random variable $X$ with law $\mu$,
which represents this long-term average. Using the inequality from above, reformulated one last time as
\begin{align*}
  \mu\left(x\right) \leq \overbrace{\left(1 + C_{2,\Delta}\right)}^{=:C_{3,\Delta}} e^{-\overbrace{\scriptstyle \min\left\{c_{\alpha,\Delta}, \alpha\right\}}^{=: d_{\alpha,\Delta} > 0}\left|x\right|}
  = C_{3,\Delta} e^{-d_{\alpha,\Delta}\left|x\right|}
\end{align*}
it is easy to see that
\begin{align*}
  \Ex{X} = \sum_{x \in \bbZ \setminus \left\{0\right\}} x \mu\left(x\right)
  = \sum_{x \in \bbN_{> 0}} x \mu\left(x\right) - \sum_{x \in \bbN_{> 0}} x \mu\left(x\right) &= 0 \\
  \textrm{since } \sum_{x \in \bbN_{> 0}} x \mu\left(x\right) \leq C_{3,\Delta} \sum_{x \in \bbN_{> 0}} x e^{-d_{\alpha,\Delta}x} &< \infty
\end{align*}

Indeed, all moments exist by the same argument (but they are $\neq 0$ for even exponents, of course).
Next to the expectation, which must be $0$ by symmetry (if existent), the variance is also
interesting. The inequality from above can be used to bound the variance, even though the
bound seems to be far from the actual value and its usefulness is therefore questionable.
\begin{align*}
  \Var{X} &= \Ex{X^2} = 2 \sum_{x \in \bbN_{> 0}} x^2 \mu\left(x\right)
  \leq 2C_{3,\Delta} \sum_{x \in \bbN_{> 0}} x^2 e^{-d_{\alpha,\Delta}x}
  = 2C_{3,\Delta} \cdot \frac{e^{d_{\alpha,\Delta}}\left(1 + e^{d_{\alpha,\Delta}}\right)}{\left(e^{d_{\alpha,\Delta}} - 1\right)^3}
\end{align*}
For the above calculations, $\alpha = \log\left(E\left(\Delta\right)^{-2}\right)$ can be chosen,
and then $d_{\alpha,\Delta} = \log\left(E\left(\Delta\right)^{-\frac{1}{4}}\right)$ holds.
In this case, the variables take the following values:
\begin{table}[H]
  \begin{center}
    \begin{tabular}{cccccccc}
      $\Delta =$ & $0.1$ & $0.2$ & $0.5$ & $1$ & $2$ & $5$ & check on WolframAlpha \\ \hline
      $E\left(\Delta\right) \approx$ & $0.9870$ & $0.9730$ & $0.9270$ & $0.8472$ & $0.7071$ & $0.4576$ & \href{https://www.wolframalpha.com/input/?i=%28+gamma%28%282+%2B+p%29+%2F+%284+*+p%29%29%5E2+%2F+%28gamma%281+%2F+%282+*+p%29%29+*+gamma%28%281+%2B+p%29+%2F+%282+*+p%29%29%29+%29+for+p+%3D+1}{https://bit.ly/3j0273B} \\
      $C_{1,\Delta} \approx$ & $0.9933$ & $0.9858$ & $0.9600$ & $0.9139$ & $0.8346$ & $0.7040$ & \href{https://www.wolframalpha.com/input/?i=%28+%28gamma%28%282+%2B+p%29+%2F+%284+*+p%29%29%5E2+*+gamma%281+%2F+p%29+*+gamma%28%282+%2B+p%29+%2F+%282+*+p%29%29%29+%2F+%28gamma%28%284+%2B+p%29+%2F+%284+*+p%29%29%5E2+*+gamma%281+%2F+%282+*+p%29%29+*+gamma%28%281+%2B+p%29+%2F+%282+*+p%29%29%29+%29+for+p+%3D+1}{https://bit.ly/3hbKMEF} \\
      $C_{3,\Delta} \approx$ & $2.003$ & $2.007$ & $2.018$ & $2.039$ & $2.086$ & $2.240$ &  \\
      $e^{d_{\alpha,\Delta}} \approx$ & $1.003$ & $1.007$ & $1.019$ & $1.042$ & $1.091$ & $1.216$ &  \\
      $\Var{X} \lesssim$ & $3 \cdot 10^8$ & $2 \cdot 10^7$ & $10^6$ & $10^5$ & $10^4$ & $10^3$ & 
    \end{tabular}
    \caption{Variable values for different $\Delta$}
    \label{tab:deltavars}
  \end{center}
\end{table}
Simulations seem to show that the bounds on $\Var{X}$ are much higher than the actual value which
seems to be close to around $\frac{4}{\Delta}$ for $0.01 \leq \Delta \leq 4$ (the range for which
simulations were run).

Summarizing the results of this section, the following theorem can be given. All assertions were
already proved in the preceding paragraphs.
\begin{theorem}[Stationary Distribution for the ERRW on $\bbZ$]
  \label{thm:zstat}
  Let $\mu$ be the probability measure on $\bbZ$ as given in \Cref{equ:zstat}, representing the
  long-term behavior of the ERRW on $\bbZ$ in the sense of \Cref{equ:zlongterm}.
  Then $\mu$ is symmetric about zero, i.e.~$\mu\left(-x\right) = \mu\left(x\right)$ for $x \in \bbZ$,
  and $\mu\left(x\right) \leq C e^{-d\left|x\right|}$ where $C,d > 0$ are constants depending only on $\Delta$
  (with $d \to 0$ for $\Delta \to 0$ and $d \to \infty$ for $\Delta \to \infty$).

  Let $X$ be an integer-valued random variable with law $\mu$. Then $\Ex{X} = 0$ and all moments
  of $X$ exist and are finite ($\Ex{X^n} = 0$ for $n \geq 0$ odd, and $0 < \Ex{X^n} < \infty$ for
  $n \geq 0$ even).
\end{theorem}

\subsection{Expected Return Time}
\label{ssec:zexpret}

The expected return time to the root is always $\infty$, as mentioned in \cite[page 1230]{errwpemantle}. To see this,
consider the number of steps $S$ until the root $x = 0$ is visited again after starting
at the root. Then $\Ex{S} = \sum_{n \geq 0} \Prb{S > n}$. Now, $\Prb{S > n}$ can be bounded from below
by the probability to first go to $x = 1$ and then go from $x = 1$ to $x = 2$ at least $n + 1$ times
before returning to the root (this is a very crude bound). The probability to first go to $x = 1$ is $\frac{1}{2}$.
At this point, the graph looks like this:

\begin{figure}[H]
  \begin{center}
    \begin{tikzpicture}
      \node[circle,fill=black,inner sep=0.7mm] (N0) at (0, 0) {};
      \node[circle,fill=black,inner sep=0.7mm] (N1) at (2, 0) {};
      \node[circle,fill=black,inner sep=0.7mm] (N2) at (4, 0) {};
      \node[circle,fill=black,inner sep=0.7mm] (N3) at (6, 0) {};
      \node[circle,fill=black,inner sep=0.7mm] (Nm1) at (-2, 0) {};
      \node[circle,fill=black,inner sep=0.7mm] (Nm2) at (-4, 0) {};
      \node[below=5mm] at (N0) {$0$};
      \node[below=5mm] at (N1) {$1$};
      \draw (Nm2) -- node[above] {$1$} (Nm1) -- node[above] {$1$} (N0) -- node[above] {$1 + \Delta$} (N1) -- node[above] {$1$} (N2) -- node[above] {$1$} (N3);
      \draw[dashed] (N3) -- (7, 0);
      \draw[dashed] (Nm2) -- (-5, 0);
    \end{tikzpicture}
    \caption{Edge weights after one step of the ERRW on $\bbZ$}
    \label{fig:errw_z_onestep}
  \end{center}
\end{figure}
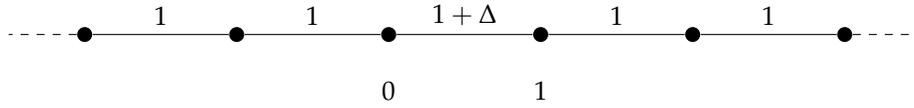

The probability to go to $x = 2$ in the next step is $\frac{1}{2 + \Delta}$. After having returned to $x = 1$,
the graph looks like this:

\begin{figure}[H]
  \begin{center}
    \begin{tikzpicture}
      \node[circle,fill=black,inner sep=0.7mm] (N0) at (0, 0) {};
      \node[circle,fill=black,inner sep=0.7mm] (N1) at (2, 0) {};
      \node[circle,fill=black,inner sep=0.7mm] (N2) at (4, 0) {};
      \node[circle,fill=black,inner sep=0.7mm] (N3) at (6, 0) {};
      \node[circle,fill=black,inner sep=0.7mm] (Nm1) at (-2, 0) {};
      \node[circle,fill=black,inner sep=0.7mm] (Nm2) at (-4, 0) {};
      \node[below=5mm] at (N0) {$0$};
      \node[below=5mm] at (N1) {$1$};
      \draw (Nm2) -- node[above] {$1$} (Nm1) -- node[above] {$1$} (N0) -- node[above] {$1 + \Delta$} (N1) -- node[above] {$1 + 2\Delta$} (N2) -- node[above] {?} (N3);
      \draw[dashed] (N3) -- (7, 0);
      \draw[dashed] (Nm2) -- (-5, 0);
    \end{tikzpicture}
    \caption{Edge weights after multiple steps of the ERRW on $\bbZ$}
    \label{fig:errw_z_steps}
  \end{center}
\end{figure}
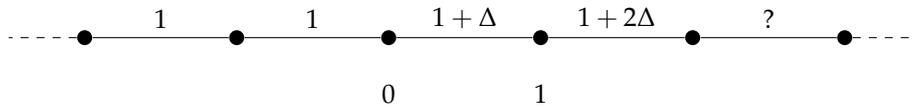

The probability to go to $x = 2$ again is thus $\frac{1 + 2\Delta}{2 + 3\Delta}$ now. Hence, the probability
to go to $x=1$ and then from $x = 1$ to $x = 2$ at least $n + 1$ times is at least
\begin{align*}
  \frac{1}{2} \cdot \frac{1}{2 + \Delta} \cdot \frac{1 + 2\Delta}{2 + 3\Delta} \cdot \ldots \cdot \frac{1 + 2n\Delta}{2 + \left(2n + 1\right)\Delta}
  \;\; \overset{\textrm{for }\Delta \geq 1}{\geq} \;\; \frac{1}{2} \cdot \frac{1}{2 + \left(2n + 1\right)\Delta} = \frac{1}{4 + 4n\Delta + 2\Delta}
\end{align*}
Using this as a lower bound for $\Prb{S > n}$, one gets for $\Delta \geq 1$
\begin{align*}
  \Ex{S} = \sum_{n \geq 0} \Prb{S > n} \geq \sum_{n \geq 0} \frac{1}{4 + 4n\Delta + 2\Delta} = \infty
\end{align*}
\cite{errwpemantle} also deals with the case $\Delta < 1$, where the expected return time can also
be shown to be infinite.

Together with \autoref{ssec:zstatdistr}, this shows that the ERRW on $\bbZ$ does have a
kind of stationary distribution, while the expected return time to the root is infinite. Hence, the
ERRW behaves not at all like a MC. This section only examined the
ERRW on $\bbZ$ in order to give a very concrete example which already demonstrates how
difficult it can be to analyze the ERRW on simple trees: the mixed stationary distribution
from \autoref{ssec:zstatdistr} is very hard to analyze because of the need to normalize the stationary
measures for every $\omega$, and the fact that the expected return time is infinite might seem
counterintuitive when comparing the ERRW with MCs. The results actually
generalize to the ERRW on general trees in the recurrent case. More details can be found
in \cite{errwpemantle}.

\newpage

\clearpage
\newpage
\thispagestyle{empty}
{
	\makeatletter
	\centering
	~\\[19em]
	\def\svgwidth{100pt}
	\import{svg-inkscape/}{TUM_svg-tex.pdf_tex}\\[10em]
	\def\svgwidth{60pt}
	\import{svg-inkscape/}{tum_mathematik_svg-tex.pdf_tex}\\

	\AddToShipoutPicture*{\put(0,0){%
		\parbox[b][\paperheight]{\paperwidth}{%
		\vfill
		\centering
		{\transparent{0.15}\includegraphics[width=\paperwidth]{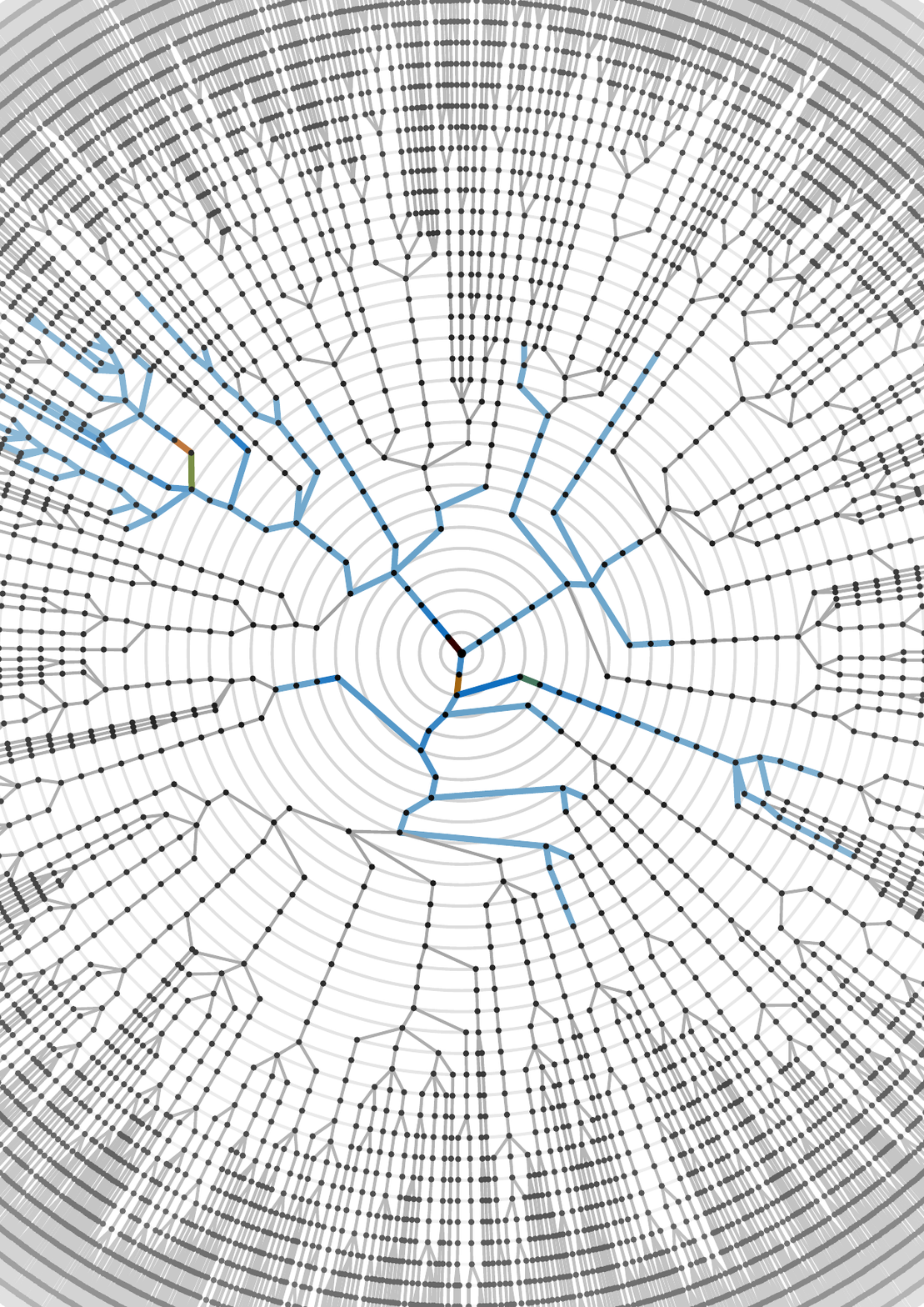}}%
		\vfill
	}}}
	\makeatother
}

\end{document}

%% file: svg-inkscape/tum_mathematik_svg-tex.pdf_tex
\begingroup%
  \makeatletter%
  \providecommand\color[2][]{%
    \errmessage{(Inkscape) Color is used for the text in Inkscape, but the package 'color.sty' is not loaded}%
    \renewcommand\color[2][]{}%
  }%
  \providecommand\transparent[1]{%
    \errmessage{(Inkscape) Transparency is used (non-zero) for the text in Inkscape, but the package 'transparent.sty' is not loaded}%
    \renewcommand\transparent[1]{}%
  }%
  \providecommand\rotatebox[2]{#2}%
  \newcommand*\fsize{\dimexpr\f@size pt\relax}%
  \newcommand*\lineheight[1]{\fontsize{\fsize}{#1\fsize}\selectfont}%
  \ifx\svgwidth\undefined%
    \setlength{\unitlength}{29.76198926bp}%
    \ifx\svgscale\undefined%
      \relax%
    \else%
      \setlength{\unitlength}{\unitlength * \real{\svgscale}}%
    \fi%
  \else%
    \setlength{\unitlength}{\svgwidth}%
  \fi%
  \global\let\svgwidth\undefined%
  \global\let\svgscale\undefined%
  \makeatother%
  \begin{picture}(1,1.00003394)%
    \lineheight{1}%
    \setlength\tabcolsep{0pt}%
    \put(0,0){\includegraphics[width=\unitlength,page=1]{tum_mathematik_svg-tex.pdf}}%
  \end{picture}%
\endgroup%

%% file: svg-inkscape/TUM_svg-tex.pdf_tex
\begingroup%
  \makeatletter%
  \providecommand\color[2][]{%
    \errmessage{(Inkscape) Color is used for the text in Inkscape, but the package 'color.sty' is not loaded}%
    \renewcommand\color[2][]{}%
  }%
  \providecommand\transparent[1]{%
    \errmessage{(Inkscape) Transparency is used (non-zero) for the text in Inkscape, but the package 'transparent.sty' is not loaded}%
    \renewcommand\transparent[1]{}%
  }%
  \providecommand\rotatebox[2]{#2}%
  \newcommand*\fsize{\dimexpr\f@size pt\relax}%
  \newcommand*\lineheight[1]{\fontsize{\fsize}{#1\fsize}\selectfont}%
  \ifx\svgwidth\undefined%
    \setlength{\unitlength}{54.75bp}%
    \ifx\svgscale\undefined%
      \relax%
    \else%
      \setlength{\unitlength}{\unitlength * \real{\svgscale}}%
    \fi%
  \else%
    \setlength{\unitlength}{\svgwidth}%
  \fi%
  \global\let\svgwidth\undefined%
  \global\let\svgscale\undefined%
  \makeatother%
  \begin{picture}(1,0.52054795)%
    \lineheight{1}%
    \setlength\tabcolsep{0pt}%
    \put(0,0){\includegraphics[width=\unitlength,page=1]{TUM_svg-tex.pdf}}%
  \end{picture}%
\endgroup%

%% file: figures/tree_notation.tex
\begin{tikzpicture}
  \node[circle,fill=black,inner sep=0.7mm] (rho) at (0, 0) {};
  \node[circle,fill=black,inner sep=0.7mm] (rho1) at (2, 1) {};
  \node[circle,fill=black,inner sep=0.7mm] (rho2) at (2, 0) {};
  \node[circle,fill=black,inner sep=0.7mm] (rho3) at (2, -1) {};
  \node[circle,fill=black,inner sep=0.7mm] (vpar) at (4, 1) {};
  \node[circle,fill=black,inner sep=0.7mm] (u1) at (6, 2.5) {};
  \node[circle,fill=black,inner sep=0.7mm] (u2) at (6, -0.5) {};
  \node[circle,fill=tumBlue,inner sep=0.7mm] (v) at (6, 1) {};
  \node[circle,fill=black,inner sep=0.7mm] (v1) at (8, 2) {};
  \node[circle,fill=black,inner sep=0.7mm] (v2) at (8, 0) {};

  \draw (rho) -- node[above=2mm] {$e_{\rho_1}$} (rho1) -- (vpar) -- node[above] {$e_v$} (v) -- (v1);
  \draw (v) -- (v2);
  \draw (vpar) -- (u1);
  \draw (vpar) -- (u2);
  \draw (rho) -- (rho2);
  \draw (rho) -- (rho3);
  \draw[rounded corners,dashed,tumBlue] (1.7, 1.8) -- (1.7, -1.3) -- (2.3, -1.3) -- (2.3, 1.8) -- node[above] {$T_1$} cycle;
  \draw[dashed,tumGray] (rho1) edge[bend left=20] node[above=1mm] {$\rho_1 < u$} (u1);

  \node[left=1mm] at (rho) {$\rho$};
  \node[above=1mm] at (rho1) {$\rho_1$};
  \node[above=1mm] at (rho2) {$\rho_2$};
  \node[above=1mm] at (rho3) {$\rho_3$};
  \node[right=3mm] at (rho2) {$\cdots$};
  \node[right=3mm] at (rho3) {$\cdots$};
  \node[above=2mm] at (u1) {$u$};
  \node[right=3mm] at (u1) {$\cdots$};
  \node[right=3mm] at (u2) {$\cdots$};
  \node[above=1mm] at (vpar) {$\overleftarrow{v}$};
  \node[above=1mm,tumBlue] at (v) {$v$};
  \node[right=4mm,tumBlue] at (v) {$b_v = 2$};
  \node[below=3mm,tumBlue] at (v) {$\left|v\right| = 3$};
  \node[above=1mm] at (v1) {$v_1$};
  \node[below=1mm] at (v2) {$v_2$};
  \node[right=3mm] at (v1) {$\cdots$};
  \node[right=3mm] at (v2) {$\cdots$};
\end{tikzpicture}

%% file: figures/rwre_instantiation.tex
\begin{tikzpicture}
  \fill[fill=linkred!15!white,rounded corners] (2.3, 2.7) -- (2.7, 2.7) -- (5.2, 1.2)
    -- (7.7, 0.2) -- (7.7, -0.6) -- (5.25, -1.2) -- (5.25, -2.2)
    -- node[below,linkred] {a cutset $\Pi$} (4.75, -2.2) -- (4.75, -0.8)
    -- (7.3, -0.2) -- (4.8, 0.8) -- (2.3, 2.3) -- cycle;

  \node[circle,fill=black,inner sep=0.7mm] (rho) at (0, 0) {};

  \node[circle,fill=linkred,inner sep=0.7mm] (rho1) at (2.5, 2.5) {};
  \node[circle,fill=black,inner sep=0.7mm] (rho2) at (2.5, 0) {};
  \node[circle,fill=black,inner sep=0.7mm] (rho3) at (2.5, -2) {};

  \node[circle,fill=black,inner sep=0.7mm] (gen2_1) at (5, 3.5) {};
  \node[circle,fill=black,inner sep=0.7mm] (gen2_2) at (5, 2) {};

  \node[circle,fill=linkred,inner sep=0.7mm] (gen2_3) at (5, 1) {};
  \node[circle,fill=black,inner sep=0.7mm] (gen2_4) at (5, 0) {};
  \node[circle,fill=linkred,inner sep=0.7mm] (gen2_5) at (5, -1) {};

  \node[circle,fill=linkred,inner sep=0.7mm] (gen2_6) at (5, -2) {};

  \node[circle,fill=black,inner sep=0.7mm] (gen3_1) at (7.5, 3.8) {};
  \node[circle,fill=black,inner sep=0.7mm] (gen3_2) at (7.5, 3.2) {};

  \node[circle,fill=black,inner sep=0.7mm] (gen3_3) at (7.5, 2) {};

  \node[circle,fill=black,inner sep=0.7mm] (gen3_4) at (7.5, 1.6) {};
  \node[circle,fill=black,inner sep=0.7mm] (gen3_5) at (7.5, 1.2) {};
  \node[circle,fill=black,inner sep=0.7mm] (gen3_6) at (7.5, 0.8) {};
  \node[circle,fill=black,inner sep=0.7mm] (gen3_7) at (7.5, 0.4) {};
  
  \node[circle,fill=linkred,inner sep=0.7mm] (gen3_8) at (7.5, 0) {};
  \node[circle,fill=linkred,inner sep=0.7mm] (gen3_9) at (7.5, -0.4) {};
  
  \node[circle,fill=black,inner sep=0.7mm] (gen3_10) at (7.5, -0.8) {};
  \node[circle,fill=black,inner sep=0.7mm] (gen3_11) at (7.5, -1.2) {};
  
  \node[circle,fill=black,inner sep=0.7mm] (gen3_12) at (7.5, -1.6) {};
  \node[circle,fill=black,inner sep=0.7mm] (gen3_13) at (7.5, -2) {};
  \node[circle,fill=black,inner sep=0.7mm] (gen3_14) at (7.5, -2.4) {};
  
  \node (rhoabove) at (0, 0.5) {};
  \node (vabovel) at (2.3, 2.9) {};
  \node (vabove) at (2.5, 3) {};
  \node (v1above) at (4.8, 4.2) {};

  \draw (rho) -- node[below=2mm,right=-1mm,tumOrange] {$A_v = 0.5$} (rho1) -- node[above=1mm,tumOrange] {$A_{v_1} = 2$} (gen2_1) -- node[above=1mm,tumOrange] {$A_{u} = 3$} (gen3_1);
    \draw (gen2_1) -- (gen3_2);
  \draw (rho1) -- node[above,tumOrange] {$A_{v_2} = 3$} (gen2_2) -- (gen3_3);
  \draw (rho) -- (rho2) -- (gen2_3) -- (gen3_4); \draw (gen2_3) -- (gen3_5); \draw (gen2_3) -- (gen3_6); \draw (gen2_3) -- (gen3_7);
  \draw (rho2) -- (gen2_4) -- (gen3_8); \draw (gen2_4) -- (gen3_9);
  \draw (rho2) -- (gen2_5) -- (gen3_10); \draw (gen2_5) -- (gen3_11);
  \draw (rho) -- (rho3) -- (gen2_6) -- (gen3_12); \draw (gen2_6) -- (gen3_13); \draw (gen2_6) -- (gen3_14);
  
  \node[left=1mm] at (rho) {$\overleftarrow{v} = \rho$};
  \node[above=1mm] at (rho1) {$v$};
  \node[below=2mm,tumBlue] at (rho1) {$0.5$};
  \node[above=12mm,left=5mm] at (rho1) {$\mu_\omega\left(v\right) = {\color{tumBlue} 0.5 + 1 + 1.5} = 3$};
  \node[above=1mm] at (gen2_1) {$v_1$};
  \node[below=1mm,tumBlue] at (gen2_1) {$1$};
  \node[above=1mm] at (gen2_2) {$v_2$};
  \node[below=1mm,tumBlue] at (gen2_2) {$1.5$};
  \node[above=1mm] at (gen3_1) {$u$};
  \node[right=1mm,tumBlue] at (gen3_1) {$C_u = {\color{tumOrange} 0.5 \cdot 2 \cdot 3} = 3$};

  \node[right=10mm] at (gen3_8) {$\cdots \qquad \infty$};

  \draw[->] (vabove) edge[bend left] node[above=3mm] {$\omega\left(v, v_1\right) = \frac{1}{3}$} (v1above);
  \draw[->] (vabovel) edge[bend right] node[above,left=3mm] {$\omega\left(v, \overleftarrow{v}\right) = \frac{1}{6}$} (rhoabove);
\end{tikzpicture}

%% file: figures/errw_example.tex
\begin{tikzpicture}[scale=0.6]
  \node[circle,fill=black,inner sep=0.7mm] (rho) at (0, 0) {};

  \node[circle,fill=black,inner sep=0.7mm] (rho1) at (2.5, 2.5) {};
  \node[circle,fill=black,inner sep=0.7mm] (rho2) at (2.5, 0) {};
  \node[circle,fill=black,inner sep=0.7mm] (rho3) at (2.5, -2) {};

  \node[circle,fill=black,inner sep=0.7mm] (gen2_1) at (5, 3) {};
  \node[circle,fill=black,inner sep=0.7mm] (gen2_2) at (5, 2) {};

  \node[circle,fill=black,inner sep=0.7mm] (gen2_3) at (5, 1) {};
  \node[circle,fill=black,inner sep=0.7mm] (gen2_4) at (5, 0) {};
  \node[circle,fill=black,inner sep=0.7mm] (gen2_5) at (5, -1) {};

  \node[circle,fill=black,inner sep=0.7mm] (gen2_6) at (5, -2) {};

  \node[circle,fill=black,inner sep=0.7mm] (gen3_1) at (7.5, 3.3) {};
  \node[circle,fill=black,inner sep=0.7mm] (gen3_2) at (7.5, 2.7) {};

  \node[circle,fill=black,inner sep=0.7mm] (gen3_3) at (7.5, 2) {};

  \node[circle,fill=black,inner sep=0.7mm] (gen3_4) at (7.5, 1.6) {};
  \node[circle,fill=black,inner sep=0.7mm] (gen3_5) at (7.5, 1.2) {};
  \node[circle,fill=black,inner sep=0.7mm] (gen3_6) at (7.5, 0.8) {};
  \node[circle,fill=black,inner sep=0.7mm] (gen3_7) at (7.5, 0.4) {};
  
  \node[circle,fill=black,inner sep=0.7mm] (gen3_8) at (7.5, 0) {};
  \node[circle,fill=black,inner sep=0.7mm] (gen3_9) at (7.5, -0.4) {};
  
  \node[circle,fill=black,inner sep=0.7mm] (gen3_10) at (7.5, -0.8) {};
  \node[circle,fill=black,inner sep=0.7mm] (gen3_11) at (7.5, -1.2) {};
  
  \node[circle,fill=black,inner sep=0.7mm] (gen3_12) at (7.5, -1.6) {};
  \node[circle,fill=black,inner sep=0.7mm] (gen3_13) at (7.5, -2) {};
  \node[circle,fill=black,inner sep=0.7mm] (gen3_14) at (7.5, -2.4) {};

  \draw[->,tumOrange] (0, 0.5) -- node[above=3mm,tumBlue] {$1 + \Delta$} (2.5, 3) -- node[above,tumBlue] {$1 + \Delta$} (5, 3.5);
  
  \draw (rho) -- (rho1) -- (gen2_1) -- (gen3_1); \draw (gen2_1) -- (gen3_2);
  \draw (rho1) -- node[below,tumBlue] {$1$} (gen2_2) -- (gen3_3);
  \draw (rho) -- node[above,tumBlue] {$1$} (rho2) -- (gen2_3) -- (gen3_4); \draw (gen2_3) -- (gen3_5); \draw (gen2_3) -- (gen3_6); \draw (gen2_3) -- (gen3_7);
  \draw (rho2) -- (gen2_4) -- (gen3_8); \draw (gen2_4) -- (gen3_9);
  \draw (rho2) -- (gen2_5) -- (gen3_10); \draw (gen2_5) -- (gen3_11);
  \draw (rho) -- node[above,tumBlue] {$1$} (rho3) -- (gen2_6) -- (gen3_12); \draw (gen2_6) -- (gen3_13); \draw (gen2_6) -- (gen3_14);
  
  \node[right=2mm] at (gen3_8) {$\cdots \; \infty$};
  \node[right,tumBlue] at (0, -2) {$w_2$};
\end{tikzpicture} ~~
\begin{tikzpicture}[scale=0.6]
  \node[white] at (0, 3.2) {$\Delta$};
  \node at (0, 0) {$\to$};
  \node[white] at (0, -2.1) {$\Delta$};
\end{tikzpicture} ~~
\begin{tikzpicture}[scale=0.6]
  \node[circle,fill=black,inner sep=0.7mm] (rho) at (0, 0) {};

  \node[circle,fill=black,inner sep=0.7mm] (rho1) at (2.5, 2.5) {};
  \node[circle,fill=black,inner sep=0.7mm] (rho2) at (2.5, 0) {};
  \node[circle,fill=black,inner sep=0.7mm] (rho3) at (2.5, -2) {};

  \node[circle,fill=black,inner sep=0.7mm] (gen2_1) at (5, 3) {};
  \node[circle,fill=black,inner sep=0.7mm] (gen2_2) at (5, 2) {};

  \node[circle,fill=black,inner sep=0.7mm] (gen2_3) at (5, 1) {};
  \node[circle,fill=black,inner sep=0.7mm] (gen2_4) at (5, 0) {};
  \node[circle,fill=black,inner sep=0.7mm] (gen2_5) at (5, -1) {};

  \node[circle,fill=black,inner sep=0.7mm] (gen2_6) at (5, -2) {};

  \node[circle,fill=black,inner sep=0.7mm] (gen3_1) at (7.5, 3.3) {};
  \node[circle,fill=black,inner sep=0.7mm] (gen3_2) at (7.5, 2.7) {};

  \node[circle,fill=black,inner sep=0.7mm] (gen3_3) at (7.5, 2) {};

  \node[circle,fill=black,inner sep=0.7mm] (gen3_4) at (7.5, 1.6) {};
  \node[circle,fill=black,inner sep=0.7mm] (gen3_5) at (7.5, 1.2) {};
  \node[circle,fill=black,inner sep=0.7mm] (gen3_6) at (7.5, 0.8) {};
  \node[circle,fill=black,inner sep=0.7mm] (gen3_7) at (7.5, 0.4) {};
  
  \node[circle,fill=black,inner sep=0.7mm] (gen3_8) at (7.5, 0) {};
  \node[circle,fill=black,inner sep=0.7mm] (gen3_9) at (7.5, -0.4) {};
  
  \node[circle,fill=black,inner sep=0.7mm] (gen3_10) at (7.5, -0.8) {};
  \node[circle,fill=black,inner sep=0.7mm] (gen3_11) at (7.5, -1.2) {};
  
  \node[circle,fill=black,inner sep=0.7mm] (gen3_12) at (7.5, -1.6) {};
  \node[circle,fill=black,inner sep=0.7mm] (gen3_13) at (7.5, -2) {};
  \node[circle,fill=black,inner sep=0.7mm] (gen3_14) at (7.5, -2.4) {};
  
  \draw[->,tumOrange] (0, 0.5) -- node[above=3mm,tumBlue] {$1 + \Delta$} (2.5, 3) -- node[above,tumBlue] {$1 + 2\Delta$} (5, 3.5);
  \draw[->,tumOrange] (5, 3.3) -- (2.5, 2.8);
  \draw[->,tumOrange] (2.5, 2.2) -- (5, 1.7);
  \draw[->,tumOrange] (5, 1.5) -- node[below,tumBlue] {$1 + 2\Delta$} (2.5, 2);
  
  \draw (rho) -- (rho1) -- (gen2_1) -- (gen3_1); \draw (gen2_1) -- (gen3_2);
  \draw (rho1) -- (gen2_2) -- (gen3_3);
  \draw (rho) -- node[above,tumBlue] {$1$} (rho2) -- (gen2_3) -- (gen3_4); \draw (gen2_3) -- (gen3_5); \draw (gen2_3) -- (gen3_6); \draw (gen2_3) -- (gen3_7);
  \draw (rho2) -- (gen2_4) -- (gen3_8); \draw (gen2_4) -- (gen3_9);
  \draw (rho2) -- (gen2_5) -- (gen3_10); \draw (gen2_5) -- (gen3_11);
  \draw (rho) -- node[above,tumBlue] {$1$} (rho3) -- (gen2_6) -- (gen3_12); \draw (gen2_6) -- (gen3_13); \draw (gen2_6) -- (gen3_14);
  
  \node[right=2mm] at (gen3_8) {$\cdots \; \infty$};
  \node[right,tumBlue] at (0, -2) {$w_5$};
\end{tikzpicture}

%% file: svg-inkscape/tree-m1.20-g23-d0.20-max100000_svg-tex.pdf_tex
\begingroup%
  \makeatletter%
  \providecommand\color[2][]{%
    \errmessage{(Inkscape) Color is used for the text in Inkscape, but the package 'color.sty' is not loaded}%
    \renewcommand\color[2][]{}%
  }%
  \providecommand\transparent[1]{%
    \errmessage{(Inkscape) Transparency is used (non-zero) for the text in Inkscape, but the package 'transparent.sty' is not loaded}%
    \renewcommand\transparent[1]{}%
  }%
  \providecommand\rotatebox[2]{#2}%
  \newcommand*\fsize{\dimexpr\f@size pt\relax}%
  \newcommand*\lineheight[1]{\fontsize{\fsize}{#1\fsize}\selectfont}%
  \ifx\svgwidth\undefined%
    \setlength{\unitlength}{650.725314bp}%
    \ifx\svgscale\undefined%
      \relax%
    \else%
      \setlength{\unitlength}{\unitlength * \real{\svgscale}}%
    \fi%
  \else%
    \setlength{\unitlength}{\svgwidth}%
  \fi%
  \global\let\svgwidth\undefined%
  \global\let\svgscale\undefined%
  \makeatother%
  \begin{picture}(1,0.94213573)%
    \lineheight{1}%
    \setlength\tabcolsep{0pt}%
    \put(0.48289854,0.92558754){\makebox(0,0)[t]{\lineheight{1.25}\smash{\begin{tabular}[t]{c}$\Delta = 0.2$, $m = 1.2$, $100000$ steps or until boundary reached\end{tabular}}}}%
    \put(0,0){\includegraphics[width=\unitlength,page=1]{tree-m1.20-g23-d0.20-max100000_svg-tex.pdf}}%
    \put(0.11436975,0.00513675){\makebox(0,0)[rt]{\lineheight{1.25}\smash{\begin{tabular}[t]{r}$w_n = 1$\end{tabular}}}}%
    \put(0.85099177,0.00513675){\makebox(0,0)[lt]{\lineheight{1.25}\smash{\begin{tabular}[t]{l}$w_n = 15.2$\end{tabular}}}}%
    \put(0,0){\includegraphics[width=\unitlength,page=2]{tree-m1.20-g23-d0.20-max100000_svg-tex.pdf}}%
  \end{picture}%
\endgroup%

%% file: svg-inkscape/tree-m1.20-g23-d1.10-max100000_svg-tex.pdf_tex
\begingroup%
  \makeatletter%
  \providecommand\color[2][]{%
    \errmessage{(Inkscape) Color is used for the text in Inkscape, but the package 'color.sty' is not loaded}%
    \renewcommand\color[2][]{}%
  }%
  \providecommand\transparent[1]{%
    \errmessage{(Inkscape) Transparency is used (non-zero) for the text in Inkscape, but the package 'transparent.sty' is not loaded}%
    \renewcommand\transparent[1]{}%
  }%
  \providecommand\rotatebox[2]{#2}%
  \newcommand*\fsize{\dimexpr\f@size pt\relax}%
  \newcommand*\lineheight[1]{\fontsize{\fsize}{#1\fsize}\selectfont}%
  \ifx\svgwidth\undefined%
    \setlength{\unitlength}{677.77763314bp}%
    \ifx\svgscale\undefined%
      \relax%
    \else%
      \setlength{\unitlength}{\unitlength * \real{\svgscale}}%
    \fi%
  \else%
    \setlength{\unitlength}{\svgwidth}%
  \fi%
  \global\let\svgwidth\undefined%
  \global\let\svgscale\undefined%
  \makeatother%
  \begin{picture}(1,0.90453202)%
    \lineheight{1}%
    \setlength\tabcolsep{0pt}%
    \put(0.46362448,0.88864431){\makebox(0,0)[t]{\lineheight{1.25}\smash{\begin{tabular}[t]{c}$\Delta = 1.1$, $m = 1.2$, $100000$ steps or until boundary reached\end{tabular}}}}%
    \put(0,0){\includegraphics[width=\unitlength,page=1]{tree-m1.20-g23-d1.10-max100000_svg-tex.pdf}}%
    \put(0.10980488,0.00493172){\makebox(0,0)[rt]{\lineheight{1.25}\smash{\begin{tabular}[t]{r}$w_n = 1$\end{tabular}}}}%
    \put(0.81702591,0.00493172){\makebox(0,0)[lt]{\lineheight{1.25}\smash{\begin{tabular}[t]{l}$w_n = 14446.2$\end{tabular}}}}%
    \put(0,0){\includegraphics[width=\unitlength,page=2]{tree-m1.20-g23-d1.10-max100000_svg-tex.pdf}}%
  \end{picture}%
\endgroup%

%% file: simulations/plot_exat-d1.00.tex
\begin{tikzpicture}
  \begin{axis}[
    xlabel={$t$},
    ylabel={$\bEx{A^t}$},
    width=7cm,
    height=6.5cm,
    grid=major
  ]
    \addplot[color=tumBlue] coordinates {
      (-0.350, 3.127) (-0.338, 2.886) (-0.326, 2.679) (-0.314, 2.499) (-0.302, 2.343) (-0.289, 2.205) (-0.277, 2.083) (-0.265, 1.974) (-0.253, 1.877) (-0.241, 1.789) (-0.229, 1.710) (-0.217, 1.638) (-0.205, 1.573) (-0.192, 1.513) (-0.180, 1.458) (-0.168, 1.408) (-0.156, 1.362) (-0.144, 1.319) (-0.132, 1.280) (-0.120, 1.244) (-0.108, 1.210) (-0.095, 1.179) (-0.083, 1.150) (-0.071, 1.123) (-0.059, 1.099) (-0.047, 1.075) (-0.035, 1.054) (-0.023, 1.034) (-0.011, 1.015) (0.002, 0.998) (0.014, 0.982) (0.026, 0.967) (0.038, 0.953) (0.050, 0.940) (0.062, 0.929) (0.074, 0.918) (0.086, 0.908) (0.098, 0.899) (0.111, 0.891) (0.123, 0.883) (0.135, 0.877) (0.147, 0.871) (0.159, 0.865) (0.171, 0.861) (0.183, 0.857) (0.195, 0.854) (0.208, 0.851) (0.220, 0.849) (0.232, 0.848) (0.244, 0.847) (0.256, 0.847) (0.268, 0.848) (0.280, 0.849) (0.292, 0.851) (0.305, 0.854) (0.317, 0.857) (0.329, 0.861) (0.341, 0.865) (0.353, 0.871) (0.365, 0.877) (0.377, 0.883) (0.389, 0.891) (0.402, 0.899) (0.414, 0.908) (0.426, 0.918) (0.438, 0.929) (0.450, 0.940) (0.462, 0.953) (0.474, 0.967) (0.486, 0.982) (0.498, 0.998) (0.511, 1.015) (0.523, 1.034) (0.535, 1.054) (0.547, 1.075) (0.559, 1.099) (0.571, 1.123) (0.583, 1.150) (0.595, 1.179) (0.608, 1.210) (0.620, 1.244) (0.632, 1.280) (0.644, 1.319) (0.656, 1.362) (0.668, 1.408) (0.680, 1.458) (0.692, 1.513) (0.705, 1.573) (0.717, 1.638) (0.729, 1.710) (0.741, 1.789) (0.753, 1.877) (0.765, 1.974) (0.777, 2.083) (0.789, 2.205) (0.802, 2.343) (0.814, 2.499) (0.826, 2.679) (0.838, 2.886) (0.850, 3.127)
    };
    \addplot[color=tumBlue,mark=*] coordinates {
      (0.250, 0.847)
    };
  \end{axis}
\end{tikzpicture}

%% file: simulations/plot_exat-d0.20.tex
\begin{tikzpicture}
  \begin{axis}[
    xlabel={$t$},
    ylabel={$\bEx{A^t}$},
    width=7cm,
    height=6.5cm,
    ymax=10,
    grid=major
  ]
    \addplot[color=tumBlue] coordinates {
      (-1.950, 13.534) (-1.906, 11.773) (-1.861, 10.321) (-1.817, 9.111) (-1.772, 8.093) (-1.728, 7.231) (-1.683, 6.494) (-1.639, 5.861) (-1.594, 5.314) (-1.550, 4.839) (-1.506, 4.423) (-1.461, 4.059) (-1.417, 3.738) (-1.372, 3.454) (-1.328, 3.201) (-1.283, 2.977) (-1.239, 2.777) (-1.194, 2.597) (-1.150, 2.436) (-1.106, 2.291) (-1.061, 2.160) (-1.017, 2.041) (-0.972, 1.934) (-0.928, 1.837) (-0.883, 1.749) (-0.839, 1.669) (-0.794, 1.596) (-0.750, 1.529) (-0.706, 1.468) (-0.661, 1.413) (-0.617, 1.362) (-0.572, 1.316) (-0.528, 1.274) (-0.483, 1.236) (-0.439, 1.201) (-0.394, 1.170) (-0.350, 1.141) (-0.306, 1.115) (-0.261, 1.092) (-0.217, 1.071) (-0.172, 1.052) (-0.128, 1.036) (-0.083, 1.022) (-0.039, 1.009) (0.006, 0.999) (0.050, 0.990) (0.094, 0.983) (0.139, 0.978) (0.183, 0.975) (0.228, 0.973) (0.272, 0.973) (0.317, 0.975) (0.361, 0.978) (0.406, 0.983) (0.450, 0.990) (0.494, 0.999) (0.539, 1.009) (0.583, 1.022) (0.628, 1.036) (0.672, 1.052) (0.717, 1.071) (0.761, 1.092) (0.806, 1.115) (0.850, 1.141) (0.894, 1.170) (0.939, 1.201) (0.983, 1.236) (1.028, 1.274) (1.072, 1.316) (1.117, 1.362) (1.161, 1.413) (1.206, 1.468) (1.250, 1.529) (1.294, 1.596) (1.339, 1.669) (1.383, 1.749) (1.428, 1.837) (1.472, 1.934) (1.517, 2.041) (1.561, 2.160) (1.606, 2.291) (1.650, 2.436) (1.694, 2.597) (1.739, 2.777) (1.783, 2.977) (1.828, 3.201) (1.872, 3.454) (1.917, 3.738) (1.961, 4.059) (2.006, 4.423) (2.050, 4.839) (2.094, 5.314) (2.139, 5.861) (2.183, 6.494) (2.228, 7.231) (2.272, 8.093) (2.317, 9.111) (2.361, 10.321) (2.406, 11.773) (2.450, 13.534)
    };
    \addplot[color=tumBlue,mark=*] coordinates {
      (0.250, 0.973)
    };
  \end{axis}
\end{tikzpicture}

%% file: simulations/plot_exatmod-d1.00.tex
\begin{tikzpicture}
  \begin{axis}[
    xlabel={$t$},
    ylabel={$s^{1-t} \bEx{A^t}$},
    width=7cm,
    height=6.5cm,
    grid=major
  ]
    \addplot[color=tumBlue] coordinates {
      (0.000, 1.000) (0.009, 0.988) (0.017, 0.977) (0.026, 0.967) (0.034, 0.957) (0.043, 0.948) (0.052, 0.939) (0.060, 0.931) (0.069, 0.923) (0.077, 0.915) (0.086, 0.908) (0.094, 0.902) (0.103, 0.896) (0.112, 0.890) (0.120, 0.885) (0.129, 0.880) (0.137, 0.875) (0.146, 0.871) (0.155, 0.867) (0.163, 0.864) (0.172, 0.861) (0.180, 0.858) (0.189, 0.855) (0.197, 0.853) (0.206, 0.851) (0.215, 0.850) (0.223, 0.849) (0.232, 0.848) (0.240, 0.847) (0.249, 0.847) (0.258, 0.847) (0.266, 0.848) (0.275, 0.849) (0.283, 0.850) (0.292, 0.851) (0.301, 0.853) (0.309, 0.855) (0.318, 0.857) (0.326, 0.860) (0.335, 0.863) (0.343, 0.866) (0.352, 0.870) (0.361, 0.874) (0.369, 0.879) (0.378, 0.884) (0.386, 0.889) (0.395, 0.894) (0.404, 0.900) (0.412, 0.907) (0.421, 0.914) (0.429, 0.921) (0.438, 0.929) (0.446, 0.937) (0.455, 0.946) (0.464, 0.955) (0.472, 0.965) (0.481, 0.975) (0.489, 0.986) (0.498, 0.997) (0.507, 1.009) (0.515, 1.022) (0.524, 1.035) (0.532, 1.050) (0.541, 1.064) (0.549, 1.080) (0.558, 1.097) (0.567, 1.114) (0.575, 1.132) (0.584, 1.151) (0.592, 1.172) (0.601, 1.193) (0.610, 1.216) (0.618, 1.240) (0.627, 1.265) (0.635, 1.291) (0.644, 1.319) (0.653, 1.349) (0.661, 1.381) (0.670, 1.414) (0.678, 1.450) (0.687, 1.487) (0.695, 1.527) (0.704, 1.570) (0.713, 1.616) (0.721, 1.664) (0.730, 1.716) (0.738, 1.772) (0.747, 1.832) (0.756, 1.896) (0.764, 1.966) (0.773, 2.041) (0.781, 2.122) (0.790, 2.210) (0.798, 2.307) (0.807, 2.412) (0.816, 2.527) (0.824, 2.655) (0.833, 2.796) (0.841, 2.952) (0.850, 3.127)
    };
    \addplot[color=tumBlue,mark=*] coordinates {
      (0.249, 0.847)
    };
    \addplot[color=tumBlue!85] coordinates {
      (0.000, 0.750) (0.009, 0.743) (0.017, 0.737) (0.026, 0.731) (0.034, 0.725) (0.043, 0.720) (0.052, 0.715) (0.060, 0.710) (0.069, 0.706) (0.077, 0.702) (0.086, 0.698) (0.094, 0.695) (0.103, 0.692) (0.112, 0.689) (0.120, 0.687) (0.129, 0.685) (0.137, 0.683) (0.146, 0.681) (0.155, 0.680) (0.163, 0.679) (0.172, 0.678) (0.180, 0.678) (0.189, 0.677) (0.197, 0.677) (0.206, 0.678) (0.215, 0.678) (0.223, 0.679) (0.232, 0.680) (0.240, 0.681) (0.249, 0.683) (0.258, 0.684) (0.266, 0.686) (0.275, 0.689) (0.283, 0.691) (0.292, 0.694) (0.301, 0.697) (0.309, 0.701) (0.318, 0.704) (0.326, 0.708) (0.335, 0.713) (0.343, 0.717) (0.352, 0.722) (0.361, 0.727) (0.369, 0.733) (0.378, 0.739) (0.386, 0.745) (0.395, 0.751) (0.404, 0.758) (0.412, 0.766) (0.421, 0.773) (0.429, 0.782) (0.438, 0.790) (0.446, 0.799) (0.455, 0.808) (0.464, 0.818) (0.472, 0.829) (0.481, 0.840) (0.489, 0.851) (0.498, 0.863) (0.507, 0.876) (0.515, 0.889) (0.524, 0.903) (0.532, 0.917) (0.541, 0.933) (0.549, 0.949) (0.558, 0.966) (0.567, 0.983) (0.575, 1.002) (0.584, 1.022) (0.592, 1.042) (0.601, 1.064) (0.610, 1.087) (0.618, 1.111) (0.627, 1.136) (0.635, 1.163) (0.644, 1.191) (0.653, 1.221) (0.661, 1.253) (0.670, 1.286) (0.678, 1.322) (0.687, 1.359) (0.695, 1.399) (0.704, 1.442) (0.713, 1.487) (0.721, 1.536) (0.730, 1.588) (0.738, 1.644) (0.747, 1.703) (0.756, 1.768) (0.764, 1.837) (0.773, 1.912) (0.781, 1.993) (0.790, 2.081) (0.798, 2.177) (0.807, 2.282) (0.816, 2.397) (0.824, 2.524) (0.833, 2.664) (0.841, 2.821) (0.850, 2.995)
    };
    \addplot[color=tumBlue!85,mark=*] coordinates {
      (0.197, 0.677)
    };
    \addplot[color=tumBlue!70] coordinates {
      (0.000, 0.500) (0.009, 0.497) (0.017, 0.495) (0.026, 0.492) (0.034, 0.490) (0.043, 0.488) (0.052, 0.487) (0.060, 0.485) (0.069, 0.484) (0.077, 0.483) (0.086, 0.482) (0.094, 0.481) (0.103, 0.481) (0.112, 0.481) (0.120, 0.481) (0.129, 0.481) (0.137, 0.481) (0.146, 0.482) (0.155, 0.483) (0.163, 0.484) (0.172, 0.485) (0.180, 0.486) (0.189, 0.487) (0.197, 0.489) (0.206, 0.491) (0.215, 0.493) (0.223, 0.495) (0.232, 0.498) (0.240, 0.501) (0.249, 0.503) (0.258, 0.506) (0.266, 0.510) (0.275, 0.513) (0.283, 0.517) (0.292, 0.521) (0.301, 0.525) (0.309, 0.530) (0.318, 0.534) (0.326, 0.539) (0.335, 0.544) (0.343, 0.550) (0.352, 0.555) (0.361, 0.561) (0.369, 0.567) (0.378, 0.574) (0.386, 0.581) (0.395, 0.588) (0.404, 0.595) (0.412, 0.603) (0.421, 0.612) (0.429, 0.620) (0.438, 0.629) (0.446, 0.638) (0.455, 0.648) (0.464, 0.658) (0.472, 0.669) (0.481, 0.680) (0.489, 0.692) (0.498, 0.704) (0.507, 0.717) (0.515, 0.730) (0.524, 0.744) (0.532, 0.759) (0.541, 0.774) (0.549, 0.790) (0.558, 0.807) (0.567, 0.825) (0.575, 0.843) (0.584, 0.863) (0.592, 0.883) (0.601, 0.905) (0.610, 0.928) (0.618, 0.951) (0.627, 0.976) (0.635, 1.003) (0.644, 1.031) (0.653, 1.060) (0.661, 1.092) (0.670, 1.125) (0.678, 1.160) (0.687, 1.197) (0.695, 1.237) (0.704, 1.279) (0.713, 1.324) (0.721, 1.372) (0.730, 1.423) (0.738, 1.478) (0.747, 1.537) (0.756, 1.601) (0.764, 1.669) (0.773, 1.743) (0.781, 1.824) (0.790, 1.911) (0.798, 2.006) (0.807, 2.110) (0.816, 2.224) (0.824, 2.350) (0.833, 2.490) (0.841, 2.645) (0.850, 2.819)
    };
    \addplot[color=tumBlue!70,mark=*] coordinates {
      (0.120, 0.481)
    };
    \addplot[color=tumBlue!55] coordinates {
      (0.000, 0.300) (0.009, 0.300) (0.017, 0.299) (0.026, 0.299) (0.034, 0.299) (0.043, 0.299) (0.052, 0.300) (0.060, 0.300) (0.069, 0.301) (0.077, 0.301) (0.086, 0.302) (0.094, 0.303) (0.103, 0.304) (0.112, 0.305) (0.120, 0.307) (0.129, 0.308) (0.137, 0.310) (0.146, 0.312) (0.155, 0.313) (0.163, 0.315) (0.172, 0.317) (0.180, 0.320) (0.189, 0.322) (0.197, 0.325) (0.206, 0.327) (0.215, 0.330) (0.223, 0.333) (0.232, 0.336) (0.240, 0.340) (0.249, 0.343) (0.258, 0.347) (0.266, 0.350) (0.275, 0.354) (0.283, 0.359) (0.292, 0.363) (0.301, 0.367) (0.309, 0.372) (0.318, 0.377) (0.326, 0.382) (0.335, 0.387) (0.343, 0.393) (0.352, 0.399) (0.361, 0.405) (0.369, 0.411) (0.378, 0.418) (0.386, 0.425) (0.395, 0.432) (0.404, 0.439) (0.412, 0.447) (0.421, 0.455) (0.429, 0.463) (0.438, 0.472) (0.446, 0.481) (0.455, 0.491) (0.464, 0.501) (0.472, 0.511) (0.481, 0.522) (0.489, 0.533) (0.498, 0.545) (0.507, 0.557) (0.515, 0.570) (0.524, 0.584) (0.532, 0.598) (0.541, 0.612) (0.549, 0.628) (0.558, 0.644) (0.567, 0.661) (0.575, 0.679) (0.584, 0.698) (0.592, 0.717) (0.601, 0.738) (0.610, 0.760) (0.618, 0.783) (0.627, 0.807) (0.635, 0.832) (0.644, 0.859) (0.653, 0.888) (0.661, 0.918) (0.670, 0.950) (0.678, 0.984) (0.687, 1.020) (0.695, 1.059) (0.704, 1.099) (0.713, 1.143) (0.721, 1.190) (0.730, 1.240) (0.738, 1.293) (0.747, 1.351) (0.756, 1.413) (0.764, 1.480) (0.773, 1.552) (0.781, 1.631) (0.790, 1.716) (0.798, 1.810) (0.807, 1.912) (0.816, 2.024) (0.824, 2.148) (0.833, 2.286) (0.841, 2.439) (0.850, 2.611)
    };
    \addplot[color=tumBlue!55,mark=*] coordinates {
      (0.026, 0.299)
    };
    \addplot[color=tumBlue!40] coordinates {
      (0.000, 0.100) (0.009, 0.101) (0.017, 0.102) (0.026, 0.103) (0.034, 0.104) (0.043, 0.105) (0.052, 0.106) (0.060, 0.107) (0.069, 0.108) (0.077, 0.109) (0.086, 0.111) (0.094, 0.112) (0.103, 0.114) (0.112, 0.115) (0.120, 0.117) (0.129, 0.118) (0.137, 0.120) (0.146, 0.122) (0.155, 0.124) (0.163, 0.126) (0.172, 0.128) (0.180, 0.130) (0.189, 0.132) (0.197, 0.134) (0.206, 0.137) (0.215, 0.139) (0.223, 0.142) (0.232, 0.145) (0.240, 0.147) (0.249, 0.150) (0.258, 0.153) (0.266, 0.156) (0.275, 0.160) (0.283, 0.163) (0.292, 0.167) (0.301, 0.170) (0.309, 0.174) (0.318, 0.178) (0.326, 0.182) (0.335, 0.187) (0.343, 0.191) (0.352, 0.196) (0.361, 0.201) (0.369, 0.206) (0.378, 0.211) (0.386, 0.216) (0.395, 0.222) (0.404, 0.228) (0.412, 0.234) (0.421, 0.241) (0.429, 0.247) (0.438, 0.255) (0.446, 0.262) (0.455, 0.270) (0.464, 0.278) (0.472, 0.286) (0.481, 0.295) (0.489, 0.304) (0.498, 0.314) (0.507, 0.324) (0.515, 0.335) (0.524, 0.346) (0.532, 0.358) (0.541, 0.370) (0.549, 0.383) (0.558, 0.396) (0.567, 0.411) (0.575, 0.426) (0.584, 0.442) (0.592, 0.458) (0.601, 0.476) (0.610, 0.495) (0.618, 0.515) (0.627, 0.536) (0.635, 0.558) (0.644, 0.581) (0.653, 0.606) (0.661, 0.633) (0.670, 0.661) (0.678, 0.691) (0.687, 0.723) (0.695, 0.758) (0.704, 0.794) (0.713, 0.834) (0.721, 0.876) (0.730, 0.921) (0.738, 0.970) (0.747, 1.023) (0.756, 1.080) (0.764, 1.142) (0.773, 1.209) (0.781, 1.283) (0.790, 1.363) (0.798, 1.450) (0.807, 1.547) (0.816, 1.653) (0.824, 1.771) (0.833, 1.902) (0.841, 2.049) (0.850, 2.214)
    };
    \addplot[color=tumBlue!40,mark=*] coordinates {
      (0.000, 0.100)
    };
  \end{axis}
\end{tikzpicture}

%% file: simulations/plot_exatmod-d0.20.tex
\begin{tikzpicture}
  \begin{axis}[
    xlabel={$t$},
    ylabel={$s^{1-t} \bEx{A^t}$},
    width=7cm,
    height=6.5cm,
    ymax=10,
    grid=major
  ]
    \addplot[color=tumBlue] coordinates {
      (0.000, 1.000) (0.025, 0.995) (0.049, 0.990) (0.074, 0.986) (0.099, 0.983) (0.124, 0.980) (0.148, 0.977) (0.173, 0.975) (0.198, 0.974) (0.223, 0.973) (0.247, 0.973) (0.272, 0.973) (0.297, 0.974) (0.322, 0.975) (0.346, 0.977) (0.371, 0.979) (0.396, 0.982) (0.421, 0.985) (0.445, 0.989) (0.470, 0.994) (0.495, 0.999) (0.520, 1.005) (0.544, 1.011) (0.569, 1.017) (0.594, 1.025) (0.619, 1.033) (0.643, 1.041) (0.668, 1.051) (0.693, 1.061) (0.718, 1.071) (0.742, 1.083) (0.767, 1.095) (0.792, 1.108) (0.817, 1.121) (0.841, 1.136) (0.866, 1.151) (0.891, 1.167) (0.916, 1.184) (0.940, 1.202) (0.965, 1.221) (0.990, 1.242) (1.015, 1.263) (1.039, 1.285) (1.064, 1.308) (1.089, 1.333) (1.114, 1.359) (1.138, 1.386) (1.163, 1.415) (1.188, 1.446) (1.213, 1.478) (1.237, 1.511) (1.262, 1.547) (1.287, 1.584) (1.312, 1.623) (1.336, 1.664) (1.361, 1.708) (1.386, 1.754) (1.411, 1.802) (1.435, 1.853) (1.460, 1.907) (1.485, 1.964) (1.510, 2.024) (1.534, 2.087) (1.559, 2.154) (1.584, 2.225) (1.609, 2.300) (1.633, 2.380) (1.658, 2.464) (1.683, 2.553) (1.708, 2.648) (1.732, 2.749) (1.757, 2.856) (1.782, 2.970) (1.807, 3.091) (1.831, 3.220) (1.856, 3.358) (1.881, 3.506) (1.906, 3.663) (1.930, 3.832) (1.955, 4.012) (1.980, 4.206) (2.005, 4.414) (2.029, 4.638) (2.054, 4.879) (2.079, 5.139) (2.104, 5.420) (2.128, 5.724) (2.153, 6.053) (2.178, 6.410) (2.203, 6.799) (2.227, 7.222) (2.252, 7.684) (2.277, 8.190) (2.302, 8.745) (2.326, 9.354) (2.351, 10.027) (2.376, 10.770) (2.401, 11.594) (2.425, 12.511) (2.450, 13.534)
    };
    \addplot[color=tumBlue,mark=*] coordinates {
      (0.247, 0.973)
    };
    \addplot[color=tumBlue!85] coordinates {
      (0.000, 0.750) (0.025, 0.751) (0.049, 0.753) (0.074, 0.756) (0.099, 0.758) (0.124, 0.761) (0.148, 0.765) (0.173, 0.769) (0.198, 0.773) (0.223, 0.778) (0.247, 0.784) (0.272, 0.789) (0.297, 0.796) (0.322, 0.802) (0.346, 0.810) (0.371, 0.817) (0.396, 0.825) (0.421, 0.834) (0.445, 0.843) (0.470, 0.853) (0.495, 0.864) (0.520, 0.875) (0.544, 0.887) (0.569, 0.899) (0.594, 0.912) (0.619, 0.926) (0.643, 0.940) (0.668, 0.955) (0.693, 0.971) (0.718, 0.988) (0.742, 1.005) (0.767, 1.024) (0.792, 1.043) (0.817, 1.064) (0.841, 1.085) (0.866, 1.108) (0.891, 1.131) (0.916, 1.156) (0.940, 1.182) (0.965, 1.209) (0.990, 1.238) (1.015, 1.268) (1.039, 1.300) (1.064, 1.333) (1.089, 1.368) (1.114, 1.404) (1.138, 1.443) (1.163, 1.483) (1.188, 1.526) (1.213, 1.571) (1.237, 1.618) (1.262, 1.668) (1.287, 1.720) (1.312, 1.775) (1.336, 1.833) (1.361, 1.895) (1.386, 1.960) (1.411, 2.028) (1.435, 2.100) (1.460, 2.177) (1.485, 2.258) (1.510, 2.343) (1.534, 2.434) (1.559, 2.530) (1.584, 2.632) (1.609, 2.740) (1.633, 2.855) (1.658, 2.977) (1.683, 3.107) (1.708, 3.246) (1.732, 3.393) (1.757, 3.551) (1.782, 3.719) (1.807, 3.898) (1.831, 4.091) (1.856, 4.296) (1.881, 4.517) (1.906, 4.753) (1.930, 5.008) (1.955, 5.281) (1.980, 5.576) (2.005, 5.894) (2.029, 6.237) (2.054, 6.608) (2.079, 7.009) (2.104, 7.445) (2.128, 7.918) (2.153, 8.434) (2.178, 8.995) (2.203, 9.609) (2.227, 10.280) (2.252, 11.016) (2.277, 11.825) (2.302, 12.716) (2.326, 13.700) (2.351, 14.789) (2.376, 15.999) (2.401, 17.346) (2.425, 18.852) (2.450, 20.540)
    };
    \addplot[color=tumBlue!85,mark=*] coordinates {
      (0.000, 0.750)
    };
    \addplot[color=tumBlue!70] coordinates {
      (0.000, 0.500) (0.025, 0.506) (0.049, 0.512) (0.074, 0.519) (0.099, 0.526) (0.124, 0.534) (0.148, 0.542) (0.173, 0.550) (0.198, 0.559) (0.223, 0.568) (0.247, 0.578) (0.272, 0.588) (0.297, 0.598) (0.322, 0.609) (0.346, 0.621) (0.371, 0.633) (0.396, 0.646) (0.421, 0.660) (0.445, 0.674) (0.470, 0.688) (0.495, 0.704) (0.520, 0.720) (0.544, 0.737) (0.569, 0.755) (0.594, 0.773) (0.619, 0.793) (0.643, 0.813) (0.668, 0.835) (0.693, 0.857) (0.718, 0.881) (0.742, 0.906) (0.767, 0.932) (0.792, 0.959) (0.817, 0.988) (0.841, 1.018) (0.866, 1.049) (0.891, 1.082) (0.916, 1.117) (0.940, 1.154) (0.965, 1.192) (0.990, 1.233) (1.015, 1.276) (1.039, 1.321) (1.064, 1.368) (1.089, 1.418) (1.114, 1.470) (1.138, 1.526) (1.163, 1.585) (1.188, 1.647) (1.213, 1.712) (1.237, 1.781) (1.262, 1.855) (1.287, 1.932) (1.312, 2.014) (1.336, 2.101) (1.361, 2.193) (1.386, 2.291) (1.411, 2.395) (1.435, 2.506) (1.460, 2.623) (1.485, 2.748) (1.510, 2.881) (1.534, 3.023) (1.559, 3.174) (1.584, 3.335) (1.609, 3.507) (1.633, 3.691) (1.658, 3.888) (1.683, 4.099) (1.708, 4.324) (1.732, 4.567) (1.757, 4.827) (1.782, 5.106) (1.807, 5.406) (1.831, 5.730) (1.856, 6.079) (1.881, 6.456) (1.906, 6.862) (1.930, 7.302) (1.955, 7.779) (1.980, 8.296) (2.005, 8.857) (2.029, 9.467) (2.054, 10.131) (2.079, 10.855) (2.104, 11.646) (2.128, 12.512) (2.153, 13.460) (2.178, 14.501) (2.203, 15.646) (2.227, 16.908) (2.252, 18.302) (2.277, 19.844) (2.302, 21.554) (2.326, 23.456) (2.351, 25.577) (2.376, 27.948) (2.401, 30.607) (2.425, 33.599) (2.450, 36.977)
    };
    \addplot[color=tumBlue!70,mark=*] coordinates {
      (0.000, 0.500)
    };
    \addplot[color=tumBlue!55] coordinates {
      (0.000, 0.300) (0.025, 0.307) (0.049, 0.315) (0.074, 0.324) (0.099, 0.332) (0.124, 0.341) (0.148, 0.351) (0.173, 0.361) (0.198, 0.371) (0.223, 0.382) (0.247, 0.393) (0.272, 0.405) (0.297, 0.418) (0.322, 0.431) (0.346, 0.445) (0.371, 0.459) (0.396, 0.475) (0.421, 0.491) (0.445, 0.507) (0.470, 0.525) (0.495, 0.544) (0.520, 0.563) (0.544, 0.584) (0.569, 0.606) (0.594, 0.629) (0.619, 0.653) (0.643, 0.678) (0.668, 0.705) (0.693, 0.733) (0.718, 0.763) (0.742, 0.794) (0.767, 0.827) (0.792, 0.862) (0.817, 0.899) (0.841, 0.938) (0.866, 0.980) (0.891, 1.024) (0.916, 1.070) (0.940, 1.119) (0.965, 1.171) (0.990, 1.226) (1.015, 1.285) (1.039, 1.347) (1.064, 1.413) (1.089, 1.484) (1.114, 1.558) (1.138, 1.638) (1.163, 1.722) (1.188, 1.813) (1.213, 1.909) (1.237, 2.011) (1.262, 2.120) (1.287, 2.237) (1.312, 2.362) (1.336, 2.495) (1.361, 2.638) (1.386, 2.791) (1.411, 2.954) (1.435, 3.130) (1.460, 3.318) (1.485, 3.520) (1.510, 3.738) (1.534, 3.971) (1.559, 4.223) (1.584, 4.494) (1.609, 4.786) (1.633, 5.101) (1.658, 5.441) (1.683, 5.809) (1.708, 6.207) (1.732, 6.638) (1.757, 7.105) (1.782, 7.612) (1.807, 8.163) (1.831, 8.762) (1.856, 9.413) (1.881, 10.124) (1.906, 10.898) (1.930, 11.745) (1.955, 12.670) (1.980, 13.684) (2.005, 14.795) (2.029, 16.016) (2.054, 17.357) (2.079, 18.835) (2.104, 20.465) (2.128, 22.265) (2.153, 24.258) (2.178, 26.466) (2.203, 28.920) (2.227, 31.649) (2.252, 34.693) (2.277, 38.095) (2.302, 41.905) (2.326, 46.184) (2.351, 51.000) (2.376, 56.437) (2.401, 62.593) (2.425, 69.586) (2.450, 77.555)
    };
    \addplot[color=tumBlue!55,mark=*] coordinates {
      (0.000, 0.300)
    };
    \addplot[color=tumBlue!40] coordinates {
      (0.000, 0.100) (0.025, 0.105) (0.049, 0.111) (0.074, 0.117) (0.099, 0.123) (0.124, 0.130) (0.148, 0.138) (0.173, 0.145) (0.198, 0.154) (0.223, 0.163) (0.247, 0.172) (0.272, 0.182) (0.297, 0.193) (0.322, 0.205) (0.346, 0.217) (0.371, 0.230) (0.396, 0.244) (0.421, 0.260) (0.445, 0.276) (0.470, 0.293) (0.495, 0.312) (0.520, 0.332) (0.544, 0.354) (0.569, 0.377) (0.594, 0.402) (0.619, 0.429) (0.643, 0.458) (0.668, 0.489) (0.693, 0.523) (0.718, 0.559) (0.742, 0.598) (0.767, 0.641) (0.792, 0.686) (0.817, 0.735) (0.841, 0.788) (0.866, 0.846) (0.891, 0.908) (0.916, 0.975) (0.940, 1.048) (0.965, 1.127) (0.990, 1.213) (1.015, 1.306) (1.039, 1.407) (1.064, 1.517) (1.089, 1.636) (1.114, 1.766) (1.138, 1.907) (1.163, 2.061) (1.188, 2.228) (1.213, 2.411) (1.237, 2.610) (1.262, 2.828) (1.287, 3.066) (1.312, 3.326) (1.336, 3.611) (1.361, 3.922) (1.386, 4.264) (1.411, 4.638) (1.435, 5.049) (1.460, 5.501) (1.485, 5.997) (1.510, 6.542) (1.534, 7.143) (1.559, 7.805) (1.584, 8.534) (1.609, 9.340) (1.633, 10.229) (1.658, 11.212) (1.683, 12.300) (1.708, 13.505) (1.732, 14.841) (1.757, 16.323) (1.782, 17.970) (1.807, 19.801) (1.831, 21.839) (1.856, 24.110) (1.881, 26.643) (1.906, 29.473) (1.930, 32.637) (1.955, 36.180) (1.980, 40.151) (2.005, 44.608) (2.029, 49.618) (2.054, 55.257) (2.079, 61.614) (2.104, 68.791) (2.128, 76.905) (2.153, 86.097) (2.178, 96.525) (2.203, 108.378) (2.227, 121.877) (2.252, 137.281) (2.277, 154.897) (2.302, 175.086) (2.326, 198.279) (2.351, 224.992) (2.376, 255.840) (2.401, 291.567) (2.425, 333.072) (2.450, 381.448)
    };
    \addplot[color=tumBlue!40,mark=*] coordinates {
      (0.000, 0.100)
    };
  \end{axis}
\end{tikzpicture}

%% file: simulations/plot_maxs-d1.00.tex
\begin{tikzpicture}
  \begin{axis}[
    xlabel={$s$},
    ylabel={$\inf_{t \geq 0} s^{1-t} \bEx{A^t}$},
    width=7cm,
    height=6.5cm,
    grid=major
  ]
    \addplot[color=tumBlue] coordinates {
      (0.020, 0.020) (0.030, 0.030) (0.040, 0.040) (0.050, 0.050) (0.060, 0.060) (0.069, 0.069) (0.079, 0.079) (0.089, 0.089) (0.099, 0.099) (0.109, 0.109) (0.119, 0.119) (0.129, 0.129) (0.139, 0.139) (0.149, 0.149) (0.159, 0.159) (0.168, 0.168) (0.178, 0.178) (0.188, 0.188) (0.198, 0.198) (0.208, 0.208) (0.218, 0.218) (0.228, 0.228) (0.238, 0.238) (0.248, 0.248) (0.258, 0.258) (0.267, 0.267) (0.277, 0.277) (0.287, 0.287) (0.297, 0.296) (0.307, 0.306) (0.317, 0.316) (0.327, 0.325) (0.337, 0.334) (0.347, 0.344) (0.357, 0.353) (0.366, 0.362) (0.376, 0.371) (0.386, 0.380) (0.396, 0.390) (0.406, 0.398) (0.416, 0.407) (0.426, 0.416) (0.436, 0.425) (0.446, 0.434) (0.456, 0.442) (0.465, 0.451) (0.475, 0.460) (0.485, 0.468) (0.495, 0.477) (0.505, 0.485) (0.515, 0.493) (0.525, 0.502) (0.535, 0.510) (0.545, 0.518) (0.555, 0.526) (0.564, 0.534) (0.574, 0.542) (0.584, 0.550) (0.594, 0.558) (0.604, 0.566) (0.614, 0.574) (0.624, 0.582) (0.634, 0.590) (0.644, 0.597) (0.654, 0.605) (0.663, 0.613) (0.673, 0.620) (0.683, 0.628) (0.693, 0.635) (0.703, 0.643) (0.713, 0.650) (0.723, 0.657) (0.733, 0.665) (0.743, 0.672) (0.753, 0.679) (0.762, 0.686) (0.772, 0.693) (0.782, 0.701) (0.792, 0.708) (0.802, 0.715) (0.812, 0.722) (0.822, 0.729) (0.832, 0.735) (0.842, 0.742) (0.852, 0.749) (0.861, 0.756) (0.871, 0.763) (0.881, 0.769) (0.891, 0.776) (0.901, 0.783) (0.911, 0.789) (0.921, 0.796) (0.931, 0.802) (0.941, 0.809) (0.951, 0.815) (0.960, 0.822) (0.970, 0.828) (0.980, 0.835) (0.990, 0.841) (1.000, 0.847)
    };
    \addplot[color=tumBlue,mark=*] coordinates {
      (1.000, 0.847)
    };
  \end{axis}
\end{tikzpicture}

%% file: simulations/plot_maxs-d0.20.tex
\begin{tikzpicture}
  \begin{axis}[
    xlabel={$s$},
    ylabel={$\inf_{t \geq 0} s^{1-t} \bEx{A^t}$},
    width=7cm,
    height=6.5cm,
    grid=major
  ]
    \addplot[color=tumBlue] coordinates {
      (0.020, 0.020) (0.030, 0.030) (0.040, 0.040) (0.050, 0.050) (0.060, 0.060) (0.069, 0.069) (0.079, 0.079) (0.089, 0.089) (0.099, 0.099) (0.109, 0.109) (0.119, 0.119) (0.129, 0.129) (0.139, 0.139) (0.149, 0.149) (0.159, 0.159) (0.168, 0.168) (0.178, 0.178) (0.188, 0.188) (0.198, 0.198) (0.208, 0.208) (0.218, 0.218) (0.228, 0.228) (0.238, 0.238) (0.248, 0.248) (0.258, 0.258) (0.267, 0.267) (0.277, 0.277) (0.287, 0.287) (0.297, 0.297) (0.307, 0.307) (0.317, 0.317) (0.327, 0.327) (0.337, 0.337) (0.347, 0.347) (0.357, 0.357) (0.366, 0.366) (0.376, 0.376) (0.386, 0.386) (0.396, 0.396) (0.406, 0.406) (0.416, 0.416) (0.426, 0.426) (0.436, 0.436) (0.446, 0.446) (0.456, 0.456) (0.465, 0.465) (0.475, 0.475) (0.485, 0.485) (0.495, 0.495) (0.505, 0.505) (0.515, 0.515) (0.525, 0.525) (0.535, 0.535) (0.545, 0.545) (0.555, 0.555) (0.564, 0.564) (0.574, 0.574) (0.584, 0.584) (0.594, 0.594) (0.604, 0.604) (0.614, 0.614) (0.624, 0.624) (0.634, 0.634) (0.644, 0.644) (0.654, 0.654) (0.663, 0.663) (0.673, 0.673) (0.683, 0.683) (0.693, 0.693) (0.703, 0.703) (0.713, 0.713) (0.723, 0.723) (0.733, 0.733) (0.743, 0.743) (0.753, 0.753) (0.762, 0.762) (0.772, 0.772) (0.782, 0.782) (0.792, 0.792) (0.802, 0.802) (0.812, 0.812) (0.822, 0.822) (0.832, 0.831) (0.842, 0.841) (0.852, 0.850) (0.861, 0.859) (0.871, 0.868) (0.881, 0.877) (0.891, 0.886) (0.901, 0.894) (0.911, 0.903) (0.921, 0.911) (0.931, 0.919) (0.941, 0.927) (0.951, 0.935) (0.960, 0.943) (0.970, 0.951) (0.980, 0.958) (0.990, 0.966) (1.000, 0.973)
    };
    \addplot[color=tumBlue,mark=*] coordinates {
      (1.000, 0.973)
    };
  \end{axis}
\end{tikzpicture}

%% file: svg-inkscape/tree-m1.50-g12-d0.70-max100000-avstree_svg-tex.pdf_tex
\begingroup%
  \makeatletter%
  \providecommand\color[2][]{%
    \errmessage{(Inkscape) Color is used for the text in Inkscape, but the package 'color.sty' is not loaded}%
    \renewcommand\color[2][]{}%
  }%
  \providecommand\transparent[1]{%
    \errmessage{(Inkscape) Transparency is used (non-zero) for the text in Inkscape, but the package 'transparent.sty' is not loaded}%
    \renewcommand\transparent[1]{}%
  }%
  \providecommand\rotatebox[2]{#2}%
  \newcommand*\fsize{\dimexpr\f@size pt\relax}%
  \newcommand*\lineheight[1]{\fontsize{\fsize}{#1\fsize}\selectfont}%
  \ifx\svgwidth\undefined%
    \setlength{\unitlength}{586.66173206bp}%
    \ifx\svgscale\undefined%
      \relax%
    \else%
      \setlength{\unitlength}{\unitlength * \real{\svgscale}}%
    \fi%
  \else%
    \setlength{\unitlength}{\svgwidth}%
  \fi%
  \global\let\svgwidth\undefined%
  \global\let\svgscale\undefined%
  \makeatother%
  \begin{picture}(1,1.04866843)%
    \lineheight{1}%
    \setlength\tabcolsep{0pt}%
    \put(0.49091321,1.03031316){\makebox(0,0)[t]{\lineheight{1.25}\smash{\begin{tabular}[t]{c}$\Delta = 0.7$, $m = 1.5$, $n = 4$, $s = 1.0$, $\varepsilon = 0.01$\end{tabular}}}}%
    \put(0,0){\includegraphics[width=\unitlength,page=1]{tree-m1.50-g12-d0.70-max100000-avstree_svg-tex.pdf}}%
    \put(0.08020826,0.00355073){\makebox(0,0)[rt]{\lineheight{1.25}\smash{\begin{tabular}[t]{r}$0$\end{tabular}}}}%
    \put(0.90210132,0.00355073){\makebox(0,0)[lt]{\lineheight{1.25}\smash{\begin{tabular}[t]{l}$10.94$\end{tabular}}}}%
    \put(0,0){\includegraphics[width=\unitlength,page=2]{tree-m1.50-g12-d0.70-max100000-avstree_svg-tex.pdf}}%
  \end{picture}%
\endgroup%

%% file: svg-inkscape/tree-m1.50-g11-d0.70-max100000_svg-tex.pdf_tex
\begingroup%
  \makeatletter%
  \providecommand\color[2][]{%
    \errmessage{(Inkscape) Color is used for the text in Inkscape, but the package 'color.sty' is not loaded}%
    \renewcommand\color[2][]{}%
  }%
  \providecommand\transparent[1]{%
    \errmessage{(Inkscape) Transparency is used (non-zero) for the text in Inkscape, but the package 'transparent.sty' is not loaded}%
    \renewcommand\transparent[1]{}%
  }%
  \providecommand\rotatebox[2]{#2}%
  \newcommand*\fsize{\dimexpr\f@size pt\relax}%
  \newcommand*\lineheight[1]{\fontsize{\fsize}{#1\fsize}\selectfont}%
  \ifx\svgwidth\undefined%
    \setlength{\unitlength}{577.077376bp}%
    \ifx\svgscale\undefined%
      \relax%
    \else%
      \setlength{\unitlength}{\unitlength * \real{\svgscale}}%
    \fi%
  \else%
    \setlength{\unitlength}{\svgwidth}%
  \fi%
  \global\let\svgwidth\undefined%
  \global\let\svgscale\undefined%
  \makeatother%
  \begin{picture}(1,1.06706752)%
    \lineheight{1}%
    \setlength\tabcolsep{0pt}%
    \put(0.49808413,1.0484074){\makebox(0,0)[t]{\lineheight{1.25}\smash{\begin{tabular}[t]{c}$\Delta = 0.7$, $m = 1.5$, $100000$ steps\end{tabular}}}}%
    \put(0,0){\includegraphics[width=\unitlength,page=1]{tree-m1.50-g11-d0.70-max100000_svg-tex.pdf}}%
    \put(0.07957557,0.0036097){\makebox(0,0)[rt]{\lineheight{1.25}\smash{\begin{tabular}[t]{r}$1$\end{tabular}}}}%
    \put(0.91610143,0.0036097){\makebox(0,0)[lt]{\lineheight{1.25}\smash{\begin{tabular}[t]{l}$89.9$\end{tabular}}}}%
    \put(0,0){\includegraphics[width=\unitlength,page=2]{tree-m1.50-g11-d0.70-max100000_svg-tex.pdf}}%
  \end{picture}%
\endgroup%

%% file: svg-inkscape/tree-m1.50-g11-d0.70-max100000-av_svg-tex.pdf_tex
\begingroup%
  \makeatletter%
  \providecommand\color[2][]{%
    \errmessage{(Inkscape) Color is used for the text in Inkscape, but the package 'color.sty' is not loaded}%
    \renewcommand\color[2][]{}%
  }%
  \providecommand\transparent[1]{%
    \errmessage{(Inkscape) Transparency is used (non-zero) for the text in Inkscape, but the package 'transparent.sty' is not loaded}%
    \renewcommand\transparent[1]{}%
  }%
  \providecommand\rotatebox[2]{#2}%
  \newcommand*\fsize{\dimexpr\f@size pt\relax}%
  \newcommand*\lineheight[1]{\fontsize{\fsize}{#1\fsize}\selectfont}%
  \ifx\svgwidth\undefined%
    \setlength{\unitlength}{577.077376bp}%
    \ifx\svgscale\undefined%
      \relax%
    \else%
      \setlength{\unitlength}{\unitlength * \real{\svgscale}}%
    \fi%
  \else%
    \setlength{\unitlength}{\svgwidth}%
  \fi%
  \global\let\svgwidth\undefined%
  \global\let\svgscale\undefined%
  \makeatother%
  \begin{picture}(1,1.06706752)%
    \lineheight{1}%
    \setlength\tabcolsep{0pt}%
    \put(0.49808413,1.0484074){\makebox(0,0)[t]{\lineheight{1.25}\smash{\begin{tabular}[t]{c}$\Delta = 0.7$, $m = 1.5$\end{tabular}}}}%
    \put(0,0){\includegraphics[width=\unitlength,page=1]{tree-m1.50-g11-d0.70-max100000-av_svg-tex.pdf}}%
    \put(0.07957557,0.0036097){\makebox(0,0)[rt]{\lineheight{1.25}\smash{\begin{tabular}[t]{r}$0$\end{tabular}}}}%
    \put(0.91610143,0.0036097){\makebox(0,0)[lt]{\lineheight{1.25}\smash{\begin{tabular}[t]{l}$9.96$\end{tabular}}}}%
    \put(0,0){\includegraphics[width=\unitlength,page=2]{tree-m1.50-g11-d0.70-max100000-av_svg-tex.pdf}}%
  \end{picture}%
\endgroup%

%% file: svg-inkscape/tree-m1.50-g11-d2.40-max100000_svg-tex.pdf_tex
\begingroup%
  \makeatletter%
  \providecommand\color[2][]{%
    \errmessage{(Inkscape) Color is used for the text in Inkscape, but the package 'color.sty' is not loaded}%
    \renewcommand\color[2][]{}%
  }%
  \providecommand\transparent[1]{%
    \errmessage{(Inkscape) Transparency is used (non-zero) for the text in Inkscape, but the package 'transparent.sty' is not loaded}%
    \renewcommand\transparent[1]{}%
  }%
  \providecommand\rotatebox[2]{#2}%
  \newcommand*\fsize{\dimexpr\f@size pt\relax}%
  \newcommand*\lineheight[1]{\fontsize{\fsize}{#1\fsize}\selectfont}%
  \ifx\svgwidth\undefined%
    \setlength{\unitlength}{613.14713486bp}%
    \ifx\svgscale\undefined%
      \relax%
    \else%
      \setlength{\unitlength}{\unitlength * \real{\svgscale}}%
    \fi%
  \else%
    \setlength{\unitlength}{\svgwidth}%
  \fi%
  \global\let\svgwidth\undefined%
  \global\let\svgscale\undefined%
  \makeatother%
  \begin{picture}(1,1.00429487)%
    \lineheight{1}%
    \setlength\tabcolsep{0pt}%
    \put(0.46878321,0.98673248){\makebox(0,0)[t]{\lineheight{1.25}\smash{\begin{tabular}[t]{c}$\Delta = 2.4$, $m = 1.5$, $100000$ steps\end{tabular}}}}%
    \put(0,0){\includegraphics[width=\unitlength,page=1]{tree-m1.50-g11-d2.40-max100000_svg-tex.pdf}}%
    \put(0.07489436,0.00339735){\makebox(0,0)[rt]{\lineheight{1.25}\smash{\begin{tabular}[t]{r}$1$\end{tabular}}}}%
    \put(0.8622097,0.00339735){\makebox(0,0)[lt]{\lineheight{1.25}\smash{\begin{tabular}[t]{l}$102584.2$\end{tabular}}}}%
    \put(0,0){\includegraphics[width=\unitlength,page=2]{tree-m1.50-g11-d2.40-max100000_svg-tex.pdf}}%
  \end{picture}%
\endgroup%

%% file: svg-inkscape/tree-m1.50-g11-d2.40-max100000-av_svg-tex.pdf_tex
\begingroup%
  \makeatletter%
  \providecommand\color[2][]{%
    \errmessage{(Inkscape) Color is used for the text in Inkscape, but the package 'color.sty' is not loaded}%
    \renewcommand\color[2][]{}%
  }%
  \providecommand\transparent[1]{%
    \errmessage{(Inkscape) Transparency is used (non-zero) for the text in Inkscape, but the package 'transparent.sty' is not loaded}%
    \renewcommand\transparent[1]{}%
  }%
  \providecommand\rotatebox[2]{#2}%
  \newcommand*\fsize{\dimexpr\f@size pt\relax}%
  \newcommand*\lineheight[1]{\fontsize{\fsize}{#1\fsize}\selectfont}%
  \ifx\svgwidth\undefined%
    \setlength{\unitlength}{586.09481571bp}%
    \ifx\svgscale\undefined%
      \relax%
    \else%
      \setlength{\unitlength}{\unitlength * \real{\svgscale}}%
    \fi%
  \else%
    \setlength{\unitlength}{\svgwidth}%
  \fi%
  \global\let\svgwidth\undefined%
  \global\let\svgscale\undefined%
  \makeatother%
  \begin{picture}(1,1.05065001)%
    \lineheight{1}%
    \setlength\tabcolsep{0pt}%
    \put(0.49042079,1.03227699){\makebox(0,0)[t]{\lineheight{1.25}\smash{\begin{tabular}[t]{c}$\Delta = 2.4$, $m = 1.5$\end{tabular}}}}%
    \put(0,0){\includegraphics[width=\unitlength,page=1]{tree-m1.50-g11-d2.40-max100000-av_svg-tex.pdf}}%
    \put(0.07835125,0.00355416){\makebox(0,0)[rt]{\lineheight{1.25}\smash{\begin{tabular}[t]{r}$0$\end{tabular}}}}%
    \put(0.90200662,0.00355416){\makebox(0,0)[lt]{\lineheight{1.25}\smash{\begin{tabular}[t]{l}$11.64$\end{tabular}}}}%
    \put(0,0){\includegraphics[width=\unitlength,page=2]{tree-m1.50-g11-d2.40-max100000-av_svg-tex.pdf}}%
  \end{picture}%
\endgroup%

%% file: svg-inkscape/tree-m1.50-g11-d5.00-max100000_svg-tex.pdf_tex
\begingroup%
  \makeatletter%
  \providecommand\color[2][]{%
    \errmessage{(Inkscape) Color is used for the text in Inkscape, but the package 'color.sty' is not loaded}%
    \renewcommand\color[2][]{}%
  }%
  \providecommand\transparent[1]{%
    \errmessage{(Inkscape) Transparency is used (non-zero) for the text in Inkscape, but the package 'transparent.sty' is not loaded}%
    \renewcommand\transparent[1]{}%
  }%
  \providecommand\rotatebox[2]{#2}%
  \newcommand*\fsize{\dimexpr\f@size pt\relax}%
  \newcommand*\lineheight[1]{\fontsize{\fsize}{#1\fsize}\selectfont}%
  \ifx\svgwidth\undefined%
    \setlength{\unitlength}{613.14713486bp}%
    \ifx\svgscale\undefined%
      \relax%
    \else%
      \setlength{\unitlength}{\unitlength * \real{\svgscale}}%
    \fi%
  \else%
    \setlength{\unitlength}{\svgwidth}%
  \fi%
  \global\let\svgwidth\undefined%
  \global\let\svgscale\undefined%
  \makeatother%
  \begin{picture}(1,1.00429487)%
    \lineheight{1}%
    \setlength\tabcolsep{0pt}%
    \put(0.46878321,0.98673248){\makebox(0,0)[t]{\lineheight{1.25}\smash{\begin{tabular}[t]{c}$\Delta = 5.0$, $m = 1.5$, $100000$ steps\end{tabular}}}}%
    \put(0,0){\includegraphics[width=\unitlength,page=1]{tree-m1.50-g11-d5.00-max100000_svg-tex.pdf}}%
    \put(0.07489436,0.00339735){\makebox(0,0)[rt]{\lineheight{1.25}\smash{\begin{tabular}[t]{r}$1$\end{tabular}}}}%
    \put(0.8622097,0.00339735){\makebox(0,0)[lt]{\lineheight{1.25}\smash{\begin{tabular}[t]{l}$463141.0$\end{tabular}}}}%
    \put(0,0){\includegraphics[width=\unitlength,page=2]{tree-m1.50-g11-d5.00-max100000_svg-tex.pdf}}%
  \end{picture}%
\endgroup%

%% file: svg-inkscape/tree-m1.50-g11-d5.00-max100000-av_svg-tex.pdf_tex
\begingroup%
  \makeatletter%
  \providecommand\color[2][]{%
    \errmessage{(Inkscape) Color is used for the text in Inkscape, but the package 'color.sty' is not loaded}%
    \renewcommand\color[2][]{}%
  }%
  \providecommand\transparent[1]{%
    \errmessage{(Inkscape) Transparency is used (non-zero) for the text in Inkscape, but the package 'transparent.sty' is not loaded}%
    \renewcommand\transparent[1]{}%
  }%
  \providecommand\rotatebox[2]{#2}%
  \newcommand*\fsize{\dimexpr\f@size pt\relax}%
  \newcommand*\lineheight[1]{\fontsize{\fsize}{#1\fsize}\selectfont}%
  \ifx\svgwidth\undefined%
    \setlength{\unitlength}{586.09481571bp}%
    \ifx\svgscale\undefined%
      \relax%
    \else%
      \setlength{\unitlength}{\unitlength * \real{\svgscale}}%
    \fi%
  \else%
    \setlength{\unitlength}{\svgwidth}%
  \fi%
  \global\let\svgwidth\undefined%
  \global\let\svgscale\undefined%
  \makeatother%
  \begin{picture}(1,1.05065001)%
    \lineheight{1}%
    \setlength\tabcolsep{0pt}%
    \put(0.49042079,1.03227699){\makebox(0,0)[t]{\lineheight{1.25}\smash{\begin{tabular}[t]{c}$\Delta = 5.0$, $m = 1.5$\end{tabular}}}}%
    \put(0,0){\includegraphics[width=\unitlength,page=1]{tree-m1.50-g11-d5.00-max100000-av_svg-tex.pdf}}%
    \put(0.07835125,0.00355416){\makebox(0,0)[rt]{\lineheight{1.25}\smash{\begin{tabular}[t]{r}$0$\end{tabular}}}}%
    \put(0.90200662,0.00355416){\makebox(0,0)[lt]{\lineheight{1.25}\smash{\begin{tabular}[t]{l}$13.03$\end{tabular}}}}%
    \put(0,0){\includegraphics[width=\unitlength,page=2]{tree-m1.50-g11-d5.00-max100000-av_svg-tex.pdf}}%
  \end{picture}%
\endgroup%

%% file: simulations/cdf-d0.01.tex
\begin{tikzpicture}
  \begin{axis}[
    xlabel={$a$},
    ylabel={$\bPrb{A_x \leq a} = \int_0^a f_{A_x}\left(y\right) \dx{y}$},
    width=7cm,
    grid=major
  ]
    \addplot[color=tumBlue] coordinates {
      (0.000, 0.000) (0.050, 0.000) (0.100, 0.000) (0.150, 0.000) (0.200, 0.000) (0.250, 0.000) (0.300, 0.000) (0.350, 0.000) (0.400, 0.000) (0.450, 0.000) (0.500, 0.000) (0.550, 0.002) (0.600, 0.006) (0.650, 0.018) (0.700, 0.042) (0.750, 0.083) (0.800, 0.144) (0.850, 0.223) (0.900, 0.317) (0.950, 0.418) (1.000, 0.520) (1.050, 0.616) (1.100, 0.701) (1.150, 0.773) (1.200, 0.832) (1.250, 0.878) (1.300, 0.913) (1.350, 0.939) (1.400, 0.958) (1.450, 0.971) (1.500, 0.981) (1.550, 0.987) (1.600, 0.992) (1.650, 0.994) (1.700, 0.996) (1.750, 0.998) (1.800, 0.999) (1.850, 0.999) (1.900, 0.999) (1.950, 1.000) (2.000, 1.000) (2.050, 1.000) (2.100, 1.000) (2.150, 1.000) (2.200, 1.000) (2.250, 1.000) (2.300, 1.000) (2.350, 1.000) (2.400, 1.000) (2.450, 1.000) (2.500, 1.000) (2.550, 1.000) (2.600, 1.000) (2.650, 1.000) (2.700, 1.000) (2.750, 1.000) (2.800, 1.000) (2.850, 1.000) (2.900, 1.000) (2.950, 1.000) (3.000, 1.000) (3.050, 1.000) (3.100, 1.000) (3.150, 1.000) (3.200, 1.000) (3.250, 1.000) (3.300, 1.000) (3.350, 1.000) (3.400, 1.000) (3.450, 1.000) (3.500, 1.000) (3.550, 1.000) (3.600, 1.000) (3.650, 1.000) (3.700, 1.000) (3.750, 1.000) (3.800, 1.000) (3.850, 1.000) (3.900, 1.000) (3.950, 1.000) (4.000, 1.000) (4.050, 1.000) (4.100, 1.000) (4.150, 1.000) (4.200, 1.000) (4.250, 1.000) (4.300, 1.000) (4.350, 1.000) (4.400, 1.000) (4.450, 1.000) (4.500, 1.000) (4.550, 1.000) (4.600, 1.000) (4.650, 1.000) (4.700, 1.000) (4.750, 1.000) (4.800, 1.000) (4.850, 1.000) (4.900, 1.000) (4.950, 1.000) (5.000, 1.000) (5.050, 1.000) (5.100, 1.000) (5.150, 1.000) (5.200, 1.000) (5.250, 1.000) (5.300, 1.000) (5.350, 1.000) (5.400, 1.000) (5.450, 1.000) (5.500, 1.000) (5.550, 1.000) (5.600, 1.000) (5.650, 1.000) (5.700, 1.000) (5.750, 1.000) (5.800, 1.000) (5.850, 1.000) (5.900, 1.000) (5.950, 1.000) (6.000, 1.000) (6.050, 1.000) (6.100, 1.000) (6.150, 1.000) (6.200, 1.000) (6.250, 1.000) (6.300, 1.000) (6.350, 1.000) (6.400, 1.000) (6.450, 1.000) (6.500, 1.000) (6.550, 1.000) (6.600, 1.000) (6.650, 1.000) (6.700, 1.000) (6.750, 1.000) (6.800, 1.000) (6.850, 1.000) (6.900, 1.000) (6.950, 1.000) (7.000, 1.000) (7.050, 1.000) (7.100, 1.000) (7.150, 1.000) (7.200, 1.000) (7.250, 1.000) (7.300, 1.000) (7.350, 1.000) (7.400, 1.000) (7.450, 1.000) (7.500, 1.000) (7.550, 1.000) (7.600, 1.000) (7.650, 1.000) (7.700, 1.000) (7.750, 1.000) (7.800, 1.000) (7.850, 1.000) (7.900, 1.000) (7.950, 1.000) (8.000, 1.000) (8.050, 1.000) (8.100, 1.000) (8.150, 1.000) (8.200, 1.000) (8.250, 1.000) (8.300, 1.000) (8.350, 1.000) (8.400, 1.000) (8.450, 1.000) (8.500, 1.000) (8.550, 1.000) (8.600, 1.000) (8.650, 1.000) (8.700, 1.000) (8.750, 1.000) (8.800, 1.000) (8.850, 1.000) (8.900, 1.000) (8.950, 1.000) (9.000, 1.000) (9.050, 1.000) (9.100, 1.000) (9.150, 1.000) (9.200, 1.000) (9.250, 1.000) (9.300, 1.000) (9.350, 1.000) (9.400, 1.000) (9.450, 1.000) (9.500, 1.000) (9.550, 1.000) (9.600, 1.000) (9.650, 1.000) (9.700, 1.000) (9.750, 1.000) (9.800, 1.000) (9.850, 1.000) (9.900, 1.000) (9.950, 1.000) (10.000, 1.000)
    };
  \end{axis}
\end{tikzpicture}

%% file: simulations/cdf-d0.10.tex
\begin{tikzpicture}
  \begin{axis}[
    xlabel={$a$},
    ylabel={$\bPrb{A_x \leq a} = \int_0^a f_{A_x}\left(y\right) \dx{y}$},
    width=7cm,
    grid=major
  ]
    \addplot[color=tumBlue] coordinates {
      (0.000, 0.000) (0.050, 0.000) (0.100, 0.001) (0.150, 0.004) (0.200, 0.012) (0.250, 0.026) (0.300, 0.046) (0.350, 0.072) (0.400, 0.103) (0.450, 0.139) (0.500, 0.178) (0.550, 0.219) (0.600, 0.260) (0.650, 0.302) (0.700, 0.344) (0.750, 0.385) (0.800, 0.424) (0.850, 0.462) (0.900, 0.498) (0.950, 0.532) (1.000, 0.564) (1.050, 0.594) (1.100, 0.622) (1.150, 0.649) (1.200, 0.673) (1.250, 0.696) (1.300, 0.717) (1.350, 0.737) (1.400, 0.756) (1.450, 0.773) (1.500, 0.788) (1.550, 0.803) (1.600, 0.816) (1.650, 0.829) (1.700, 0.840) (1.750, 0.851) (1.800, 0.861) (1.850, 0.870) (1.900, 0.879) (1.950, 0.887) (2.000, 0.894) (2.050, 0.901) (2.100, 0.907) (2.150, 0.913) (2.200, 0.918) (2.250, 0.923) (2.300, 0.928) (2.350, 0.933) (2.400, 0.937) (2.450, 0.940) (2.500, 0.944) (2.550, 0.947) (2.600, 0.950) (2.650, 0.953) (2.700, 0.956) (2.750, 0.958) (2.800, 0.961) (2.850, 0.963) (2.900, 0.965) (2.950, 0.967) (3.000, 0.969) (3.050, 0.970) (3.100, 0.972) (3.150, 0.973) (3.200, 0.975) (3.250, 0.976) (3.300, 0.977) (3.350, 0.979) (3.400, 0.980) (3.450, 0.981) (3.500, 0.982) (3.550, 0.983) (3.600, 0.983) (3.650, 0.984) (3.700, 0.985) (3.750, 0.986) (3.800, 0.986) (3.850, 0.987) (3.900, 0.988) (3.950, 0.988) (4.000, 0.989) (4.050, 0.989) (4.100, 0.990) (4.150, 0.990) (4.200, 0.991) (4.250, 0.991) (4.300, 0.991) (4.350, 0.992) (4.400, 0.992) (4.450, 0.993) (4.500, 0.993) (4.550, 0.993) (4.600, 0.993) (4.650, 0.994) (4.700, 0.994) (4.750, 0.994) (4.800, 0.994) (4.850, 0.995) (4.900, 0.995) (4.950, 0.995) (5.000, 0.995) (5.050, 0.995) (5.100, 0.996) (5.150, 0.996) (5.200, 0.996) (5.250, 0.996) (5.300, 0.996) (5.350, 0.996) (5.400, 0.997) (5.450, 0.997) (5.500, 0.997) (5.550, 0.997) (5.600, 0.997) (5.650, 0.997) (5.700, 0.997) (5.750, 0.997) (5.800, 0.997) (5.850, 0.998) (5.900, 0.998) (5.950, 0.998) (6.000, 0.998) (6.050, 0.998) (6.100, 0.998) (6.150, 0.998) (6.200, 0.998) (6.250, 0.998) (6.300, 0.998) (6.350, 0.998) (6.400, 0.998) (6.450, 0.998) (6.500, 0.998) (6.550, 0.998) (6.600, 0.999) (6.650, 0.999) (6.700, 0.999) (6.750, 0.999) (6.800, 0.999) (6.850, 0.999) (6.900, 0.999) (6.950, 0.999) (7.000, 0.999) (7.050, 0.999) (7.100, 0.999) (7.150, 0.999) (7.200, 0.999) (7.250, 0.999) (7.300, 0.999) (7.350, 0.999) (7.400, 0.999) (7.450, 0.999) (7.500, 0.999) (7.550, 0.999) (7.600, 0.999) (7.650, 0.999) (7.700, 0.999) (7.750, 0.999) (7.800, 0.999) (7.850, 0.999) (7.900, 0.999) (7.950, 0.999) (8.000, 0.999) (8.050, 0.999) (8.100, 0.999) (8.150, 0.999) (8.200, 0.999) (8.250, 0.999) (8.300, 0.999) (8.350, 0.999) (8.400, 0.999) (8.450, 1.000) (8.500, 1.000) (8.550, 1.000) (8.600, 1.000) (8.650, 1.000) (8.700, 1.000) (8.750, 1.000) (8.800, 1.000) (8.850, 1.000) (8.900, 1.000) (8.950, 1.000) (9.000, 1.000) (9.050, 1.000) (9.100, 1.000) (9.150, 1.000) (9.200, 1.000) (9.250, 1.000) (9.300, 1.000) (9.350, 1.000) (9.400, 1.000) (9.450, 1.000) (9.500, 1.000) (9.550, 1.000) (9.600, 1.000) (9.650, 1.000) (9.700, 1.000) (9.750, 1.000) (9.800, 1.000) (9.850, 1.000) (9.900, 1.000) (9.950, 1.000) (10.000, 1.000)
    };
  \end{axis}
\end{tikzpicture}

%% file: simulations/cdf-d0.20.tex
\begin{tikzpicture}
  \begin{axis}[
    xlabel={$a$},
    ylabel={$\bPrb{A_x \leq a} = \int_0^a f_{A_x}\left(y\right) \dx{y}$},
    width=7cm,
    grid=major
  ]
    \addplot[color=tumBlue] coordinates {
      (0.000, 0.000) (0.050, 0.004) (0.100, 0.017) (0.150, 0.040) (0.200, 0.069) (0.250, 0.104) (0.300, 0.141) (0.350, 0.180) (0.400, 0.219) (0.450, 0.258) (0.500, 0.296) (0.550, 0.333) (0.600, 0.368) (0.650, 0.402) (0.700, 0.434) (0.750, 0.464) (0.800, 0.492) (0.850, 0.519) (0.900, 0.545) (0.950, 0.569) (1.000, 0.591) (1.050, 0.612) (1.100, 0.632) (1.150, 0.651) (1.200, 0.668) (1.250, 0.684) (1.300, 0.700) (1.350, 0.714) (1.400, 0.728) (1.450, 0.741) (1.500, 0.753) (1.550, 0.764) (1.600, 0.775) (1.650, 0.785) (1.700, 0.795) (1.750, 0.804) (1.800, 0.812) (1.850, 0.820) (1.900, 0.828) (1.950, 0.835) (2.000, 0.842) (2.050, 0.848) (2.100, 0.854) (2.150, 0.860) (2.200, 0.866) (2.250, 0.871) (2.300, 0.876) (2.350, 0.880) (2.400, 0.885) (2.450, 0.889) (2.500, 0.893) (2.550, 0.897) (2.600, 0.901) (2.650, 0.904) (2.700, 0.908) (2.750, 0.911) (2.800, 0.914) (2.850, 0.917) (2.900, 0.920) (2.950, 0.922) (3.000, 0.925) (3.050, 0.927) (3.100, 0.930) (3.150, 0.932) (3.200, 0.934) (3.250, 0.936) (3.300, 0.938) (3.350, 0.940) (3.400, 0.942) (3.450, 0.943) (3.500, 0.945) (3.550, 0.947) (3.600, 0.948) (3.650, 0.950) (3.700, 0.951) (3.750, 0.953) (3.800, 0.954) (3.850, 0.955) (3.900, 0.956) (3.950, 0.958) (4.000, 0.959) (4.050, 0.960) (4.100, 0.961) (4.150, 0.962) (4.200, 0.963) (4.250, 0.964) (4.300, 0.965) (4.350, 0.966) (4.400, 0.967) (4.450, 0.968) (4.500, 0.968) (4.550, 0.969) (4.600, 0.970) (4.650, 0.971) (4.700, 0.971) (4.750, 0.972) (4.800, 0.973) (4.850, 0.973) (4.900, 0.974) (4.950, 0.975) (5.000, 0.975) (5.050, 0.976) (5.100, 0.976) (5.150, 0.977) (5.200, 0.977) (5.250, 0.978) (5.300, 0.978) (5.350, 0.979) (5.400, 0.979) (5.450, 0.980) (5.500, 0.980) (5.550, 0.981) (5.600, 0.981) (5.650, 0.981) (5.700, 0.982) (5.750, 0.982) (5.800, 0.982) (5.850, 0.983) (5.900, 0.983) (5.950, 0.984) (6.000, 0.984) (6.050, 0.984) (6.100, 0.984) (6.150, 0.985) (6.200, 0.985) (6.250, 0.985) (6.300, 0.986) (6.350, 0.986) (6.400, 0.986) (6.450, 0.986) (6.500, 0.987) (6.550, 0.987) (6.600, 0.987) (6.650, 0.987) (6.700, 0.988) (6.750, 0.988) (6.800, 0.988) (6.850, 0.988) (6.900, 0.989) (6.950, 0.989) (7.000, 0.989) (7.050, 0.989) (7.100, 0.989) (7.150, 0.990) (7.200, 0.990) (7.250, 0.990) (7.300, 0.990) (7.350, 0.990) (7.400, 0.990) (7.450, 0.991) (7.500, 0.991) (7.550, 0.991) (7.600, 0.991) (7.650, 0.991) (7.700, 0.991) (7.750, 0.991) (7.800, 0.992) (7.850, 0.992) (7.900, 0.992) (7.950, 0.992) (8.000, 0.992) (8.050, 0.992) (8.100, 0.992) (8.150, 0.992) (8.200, 0.993) (8.250, 0.993) (8.300, 0.993) (8.350, 0.993) (8.400, 0.993) (8.450, 0.993) (8.500, 0.993) (8.550, 0.993) (8.600, 0.993) (8.650, 0.994) (8.700, 0.994) (8.750, 0.994) (8.800, 0.994) (8.850, 0.994) (8.900, 0.994) (8.950, 0.994) (9.000, 0.994) (9.050, 0.994) (9.100, 0.994) (9.150, 0.994) (9.200, 0.994) (9.250, 0.995) (9.300, 0.995) (9.350, 0.995) (9.400, 0.995) (9.450, 0.995) (9.500, 0.995) (9.550, 0.995) (9.600, 0.995) (9.650, 0.995) (9.700, 0.995) (9.750, 0.995) (9.800, 0.995) (9.850, 0.995) (9.900, 0.995) (9.950, 0.996) (10.000, 0.996)
    };
  \end{axis}
\end{tikzpicture}

%% file: simulations/cdf-d0.50.tex
\begin{tikzpicture}
  \begin{axis}[
    xlabel={$a$},
    ylabel={$\bPrb{A_x \leq a} = \int_0^a f_{A_x}\left(y\right) \dx{y}$},
    width=7cm,
    grid=major
  ]
    \addplot[color=tumBlue] coordinates {
      (0.000, 0.000) (0.050, 0.071) (0.100, 0.133) (0.150, 0.189) (0.200, 0.239) (0.250, 0.284) (0.300, 0.325) (0.350, 0.362) (0.400, 0.396) (0.450, 0.427) (0.500, 0.456) (0.550, 0.482) (0.600, 0.506) (0.650, 0.528) (0.700, 0.549) (0.750, 0.568) (0.800, 0.586) (0.850, 0.603) (0.900, 0.618) (0.950, 0.633) (1.000, 0.646) (1.050, 0.659) (1.100, 0.671) (1.150, 0.683) (1.200, 0.694) (1.250, 0.704) (1.300, 0.713) (1.350, 0.722) (1.400, 0.731) (1.450, 0.739) (1.500, 0.747) (1.550, 0.754) (1.600, 0.761) (1.650, 0.768) (1.700, 0.775) (1.750, 0.781) (1.800, 0.787) (1.850, 0.792) (1.900, 0.798) (1.950, 0.803) (2.000, 0.808) (2.050, 0.812) (2.100, 0.817) (2.150, 0.821) (2.200, 0.825) (2.250, 0.829) (2.300, 0.833) (2.350, 0.837) (2.400, 0.840) (2.450, 0.844) (2.500, 0.847) (2.550, 0.850) (2.600, 0.854) (2.650, 0.857) (2.700, 0.859) (2.750, 0.862) (2.800, 0.865) (2.850, 0.868) (2.900, 0.870) (2.950, 0.873) (3.000, 0.875) (3.050, 0.877) (3.100, 0.880) (3.150, 0.882) (3.200, 0.884) (3.250, 0.886) (3.300, 0.888) (3.350, 0.890) (3.400, 0.892) (3.450, 0.893) (3.500, 0.895) (3.550, 0.897) (3.600, 0.899) (3.650, 0.900) (3.700, 0.902) (3.750, 0.903) (3.800, 0.905) (3.850, 0.906) (3.900, 0.908) (3.950, 0.909) (4.000, 0.911) (4.050, 0.912) (4.100, 0.913) (4.150, 0.914) (4.200, 0.916) (4.250, 0.917) (4.300, 0.918) (4.350, 0.919) (4.400, 0.920) (4.450, 0.921) (4.500, 0.922) (4.550, 0.924) (4.600, 0.925) (4.650, 0.926) (4.700, 0.927) (4.750, 0.927) (4.800, 0.928) (4.850, 0.929) (4.900, 0.930) (4.950, 0.931) (5.000, 0.932) (5.050, 0.933) (5.100, 0.934) (5.150, 0.934) (5.200, 0.935) (5.250, 0.936) (5.300, 0.937) (5.350, 0.938) (5.400, 0.938) (5.450, 0.939) (5.500, 0.940) (5.550, 0.940) (5.600, 0.941) (5.650, 0.942) (5.700, 0.942) (5.750, 0.943) (5.800, 0.944) (5.850, 0.944) (5.900, 0.945) (5.950, 0.945) (6.000, 0.946) (6.050, 0.947) (6.100, 0.947) (6.150, 0.948) (6.200, 0.948) (6.250, 0.949) (6.300, 0.949) (6.350, 0.950) (6.400, 0.950) (6.450, 0.951) (6.500, 0.951) (6.550, 0.952) (6.600, 0.952) (6.650, 0.953) (6.700, 0.953) (6.750, 0.954) (6.800, 0.954) (6.850, 0.955) (6.900, 0.955) (6.950, 0.955) (7.000, 0.956) (7.050, 0.956) (7.100, 0.957) (7.150, 0.957) (7.200, 0.957) (7.250, 0.958) (7.300, 0.958) (7.350, 0.959) (7.400, 0.959) (7.450, 0.959) (7.500, 0.960) (7.550, 0.960) (7.600, 0.960) (7.650, 0.961) (7.700, 0.961) (7.750, 0.961) (7.800, 0.962) (7.850, 0.962) (7.900, 0.962) (7.950, 0.963) (8.000, 0.963) (8.050, 0.963) (8.100, 0.964) (8.150, 0.964) (8.200, 0.964) (8.250, 0.964) (8.300, 0.965) (8.350, 0.965) (8.400, 0.965) (8.450, 0.966) (8.500, 0.966) (8.550, 0.966) (8.600, 0.966) (8.650, 0.967) (8.700, 0.967) (8.750, 0.967) (8.800, 0.967) (8.850, 0.968) (8.900, 0.968) (8.950, 0.968) (9.000, 0.968) (9.050, 0.969) (9.100, 0.969) (9.150, 0.969) (9.200, 0.969) (9.250, 0.970) (9.300, 0.970) (9.350, 0.970) (9.400, 0.970) (9.450, 0.970) (9.500, 0.971) (9.550, 0.971) (9.600, 0.971) (9.650, 0.971) (9.700, 0.971) (9.750, 0.972) (9.800, 0.972) (9.850, 0.972) (9.900, 0.972) (9.950, 0.972) (10.000, 0.973)
    };
  \end{axis}
\end{tikzpicture}

%% file: simulations/cdf-d1.00.tex
\begin{tikzpicture}
  \begin{axis}[
    xlabel={$a$},
    ylabel={$\bPrb{A_x \leq a} = \int_0^a f_{A_x}\left(y\right) \dx{y}$},
    width=7cm,
    grid=major
  ]
    \addplot[color=tumBlue] coordinates {
      (0.000, 0.000) (0.050, 0.217) (0.100, 0.301) (0.150, 0.360) (0.200, 0.408) (0.250, 0.446) (0.300, 0.480) (0.350, 0.508) (0.400, 0.534) (0.450, 0.556) (0.500, 0.577) (0.550, 0.595) (0.600, 0.612) (0.650, 0.627) (0.700, 0.641) (0.750, 0.654) (0.800, 0.666) (0.850, 0.677) (0.900, 0.688) (0.950, 0.697) (1.000, 0.706) (1.050, 0.715) (1.100, 0.723) (1.150, 0.731) (1.200, 0.738) (1.250, 0.745) (1.300, 0.751) (1.350, 0.757) (1.400, 0.763) (1.450, 0.769) (1.500, 0.774) (1.550, 0.779) (1.600, 0.784) (1.650, 0.788) (1.700, 0.793) (1.750, 0.797) (1.800, 0.801) (1.850, 0.805) (1.900, 0.809) (1.950, 0.812) (2.000, 0.816) (2.050, 0.819) (2.100, 0.822) (2.150, 0.825) (2.200, 0.828) (2.250, 0.831) (2.300, 0.834) (2.350, 0.837) (2.400, 0.839) (2.450, 0.842) (2.500, 0.844) (2.550, 0.847) (2.600, 0.849) (2.650, 0.851) (2.700, 0.854) (2.750, 0.856) (2.800, 0.858) (2.850, 0.860) (2.900, 0.862) (2.950, 0.863) (3.000, 0.865) (3.050, 0.867) (3.100, 0.869) (3.150, 0.870) (3.200, 0.872) (3.250, 0.874) (3.300, 0.875) (3.350, 0.877) (3.400, 0.878) (3.450, 0.880) (3.500, 0.881) (3.550, 0.883) (3.600, 0.884) (3.650, 0.885) (3.700, 0.887) (3.750, 0.888) (3.800, 0.889) (3.850, 0.890) (3.900, 0.891) (3.950, 0.893) (4.000, 0.894) (4.050, 0.895) (4.100, 0.896) (4.150, 0.897) (4.200, 0.898) (4.250, 0.899) (4.300, 0.900) (4.350, 0.901) (4.400, 0.902) (4.450, 0.903) (4.500, 0.904) (4.550, 0.905) (4.600, 0.906) (4.650, 0.906) (4.700, 0.907) (4.750, 0.908) (4.800, 0.909) (4.850, 0.910) (4.900, 0.911) (4.950, 0.911) (5.000, 0.912) (5.050, 0.913) (5.100, 0.914) (5.150, 0.914) (5.200, 0.915) (5.250, 0.916) (5.300, 0.916) (5.350, 0.917) (5.400, 0.918) (5.450, 0.918) (5.500, 0.919) (5.550, 0.920) (5.600, 0.920) (5.650, 0.921) (5.700, 0.922) (5.750, 0.922) (5.800, 0.923) (5.850, 0.923) (5.900, 0.924) (5.950, 0.925) (6.000, 0.925) (6.050, 0.926) (6.100, 0.926) (6.150, 0.927) (6.200, 0.927) (6.250, 0.928) (6.300, 0.928) (6.350, 0.929) (6.400, 0.929) (6.450, 0.930) (6.500, 0.930) (6.550, 0.931) (6.600, 0.931) (6.650, 0.932) (6.700, 0.932) (6.750, 0.933) (6.800, 0.933) (6.850, 0.933) (6.900, 0.934) (6.950, 0.934) (7.000, 0.935) (7.050, 0.935) (7.100, 0.936) (7.150, 0.936) (7.200, 0.936) (7.250, 0.937) (7.300, 0.937) (7.350, 0.937) (7.400, 0.938) (7.450, 0.938) (7.500, 0.939) (7.550, 0.939) (7.600, 0.939) (7.650, 0.940) (7.700, 0.940) (7.750, 0.940) (7.800, 0.941) (7.850, 0.941) (7.900, 0.941) (7.950, 0.942) (8.000, 0.942) (8.050, 0.942) (8.100, 0.943) (8.150, 0.943) (8.200, 0.943) (8.250, 0.944) (8.300, 0.944) (8.350, 0.944) (8.400, 0.945) (8.450, 0.945) (8.500, 0.945) (8.550, 0.945) (8.600, 0.946) (8.650, 0.946) (8.700, 0.946) (8.750, 0.947) (8.800, 0.947) (8.850, 0.947) (8.900, 0.947) (8.950, 0.948) (9.000, 0.948) (9.050, 0.948) (9.100, 0.948) (9.150, 0.949) (9.200, 0.949) (9.250, 0.949) (9.300, 0.949) (9.350, 0.950) (9.400, 0.950) (9.450, 0.950) (9.500, 0.950) (9.550, 0.951) (9.600, 0.951) (9.650, 0.951) (9.700, 0.951) (9.750, 0.952) (9.800, 0.952) (9.850, 0.952) (9.900, 0.952) (9.950, 0.953) (10.000, 0.953)
    };
  \end{axis}
\end{tikzpicture}

%% file: simulations/cdf-d1.50.tex
\begin{tikzpicture}
  \begin{axis}[
    xlabel={$a$},
    ylabel={$\bPrb{A_x \leq a} = \int_0^a f_{A_x}\left(y\right) \dx{y}$},
    width=7cm,
    grid=major
  ]
    \addplot[color=tumBlue] coordinates {
      (0.000, 0.000) (0.050, 0.327) (0.100, 0.408) (0.150, 0.462) (0.200, 0.503) (0.250, 0.536) (0.300, 0.563) (0.350, 0.586) (0.400, 0.607) (0.450, 0.625) (0.500, 0.641) (0.550, 0.655) (0.600, 0.668) (0.650, 0.680) (0.700, 0.691) (0.750, 0.701) (0.800, 0.710) (0.850, 0.719) (0.900, 0.727) (0.950, 0.735) (1.000, 0.742) (1.050, 0.748) (1.100, 0.755) (1.150, 0.760) (1.200, 0.766) (1.250, 0.771) (1.300, 0.776) (1.350, 0.781) (1.400, 0.786) (1.450, 0.790) (1.500, 0.794) (1.550, 0.798) (1.600, 0.802) (1.650, 0.805) (1.700, 0.809) (1.750, 0.812) (1.800, 0.815) (1.850, 0.818) (1.900, 0.821) (1.950, 0.824) (2.000, 0.827) (2.050, 0.829) (2.100, 0.832) (2.150, 0.834) (2.200, 0.837) (2.250, 0.839) (2.300, 0.841) (2.350, 0.843) (2.400, 0.846) (2.450, 0.848) (2.500, 0.850) (2.550, 0.851) (2.600, 0.853) (2.650, 0.855) (2.700, 0.857) (2.750, 0.859) (2.800, 0.860) (2.850, 0.862) (2.900, 0.863) (2.950, 0.865) (3.000, 0.866) (3.050, 0.868) (3.100, 0.869) (3.150, 0.871) (3.200, 0.872) (3.250, 0.873) (3.300, 0.874) (3.350, 0.876) (3.400, 0.877) (3.450, 0.878) (3.500, 0.879) (3.550, 0.880) (3.600, 0.882) (3.650, 0.883) (3.700, 0.884) (3.750, 0.885) (3.800, 0.886) (3.850, 0.887) (3.900, 0.888) (3.950, 0.889) (4.000, 0.890) (4.050, 0.891) (4.100, 0.891) (4.150, 0.892) (4.200, 0.893) (4.250, 0.894) (4.300, 0.895) (4.350, 0.896) (4.400, 0.897) (4.450, 0.897) (4.500, 0.898) (4.550, 0.899) (4.600, 0.900) (4.650, 0.900) (4.700, 0.901) (4.750, 0.902) (4.800, 0.902) (4.850, 0.903) (4.900, 0.904) (4.950, 0.905) (5.000, 0.905) (5.050, 0.906) (5.100, 0.906) (5.150, 0.907) (5.200, 0.908) (5.250, 0.908) (5.300, 0.909) (5.350, 0.909) (5.400, 0.910) (5.450, 0.911) (5.500, 0.911) (5.550, 0.912) (5.600, 0.912) (5.650, 0.913) (5.700, 0.913) (5.750, 0.914) (5.800, 0.914) (5.850, 0.915) (5.900, 0.915) (5.950, 0.916) (6.000, 0.916) (6.050, 0.917) (6.100, 0.917) (6.150, 0.918) (6.200, 0.918) (6.250, 0.919) (6.300, 0.919) (6.350, 0.920) (6.400, 0.920) (6.450, 0.920) (6.500, 0.921) (6.550, 0.921) (6.600, 0.922) (6.650, 0.922) (6.700, 0.922) (6.750, 0.923) (6.800, 0.923) (6.850, 0.924) (6.900, 0.924) (6.950, 0.924) (7.000, 0.925) (7.050, 0.925) (7.100, 0.925) (7.150, 0.926) (7.200, 0.926) (7.250, 0.927) (7.300, 0.927) (7.350, 0.927) (7.400, 0.928) (7.450, 0.928) (7.500, 0.928) (7.550, 0.929) (7.600, 0.929) (7.650, 0.929) (7.700, 0.930) (7.750, 0.930) (7.800, 0.930) (7.850, 0.930) (7.900, 0.931) (7.950, 0.931) (8.000, 0.931) (8.050, 0.932) (8.100, 0.932) (8.150, 0.932) (8.200, 0.933) (8.250, 0.933) (8.300, 0.933) (8.350, 0.933) (8.400, 0.934) (8.450, 0.934) (8.500, 0.934) (8.550, 0.934) (8.600, 0.935) (8.650, 0.935) (8.700, 0.935) (8.750, 0.935) (8.800, 0.936) (8.850, 0.936) (8.900, 0.936) (8.950, 0.936) (9.000, 0.937) (9.050, 0.937) (9.100, 0.937) (9.150, 0.937) (9.200, 0.938) (9.250, 0.938) (9.300, 0.938) (9.350, 0.938) (9.400, 0.939) (9.450, 0.939) (9.500, 0.939) (9.550, 0.939) (9.600, 0.939) (9.650, 0.940) (9.700, 0.940) (9.750, 0.940) (9.800, 0.940) (9.850, 0.941) (9.900, 0.941) (9.950, 0.941) (10.000, 0.941)
    };
  \end{axis}
\end{tikzpicture}

%% file: simulations/e_delta.tex
\begin{tikzpicture}
  \begin{axis}[
    xlabel={$\Delta$},
    ylabel={$E\left(\Delta\right)$},
    ymin=-0.05,
    grid=major
  ]
    \addplot[color=tumBlue] coordinates {
      (0.010, 0.999) (0.030, 0.996) (0.050, 0.994) (0.070, 0.991) (0.090, 0.988) (0.110, 0.986) (0.130, 0.983) (0.150, 0.980) (0.170, 0.977) (0.190, 0.974) (0.210, 0.972) (0.230, 0.969) (0.250, 0.966) (0.270, 0.963) (0.290, 0.960) (0.310, 0.957) (0.330, 0.954) (0.350, 0.951) (0.370, 0.947) (0.390, 0.944) (0.410, 0.941) (0.430, 0.938) (0.450, 0.935) (0.470, 0.932) (0.490, 0.929) (0.510, 0.925) (0.530, 0.922) (0.550, 0.919) (0.570, 0.916) (0.590, 0.913) (0.610, 0.909) (0.630, 0.906) (0.650, 0.903) (0.670, 0.900) (0.690, 0.897) (0.710, 0.893) (0.730, 0.890) (0.750, 0.887) (0.770, 0.884) (0.790, 0.880) (0.810, 0.877) (0.830, 0.874) (0.850, 0.871) (0.870, 0.868) (0.890, 0.864) (0.910, 0.861) (0.930, 0.858) (0.950, 0.855) (0.970, 0.852) (0.990, 0.849) (1.011, 0.846) (1.031, 0.842) (1.051, 0.839) (1.071, 0.836) (1.091, 0.833) (1.111, 0.830) (1.131, 0.827) (1.151, 0.824) (1.171, 0.821) (1.191, 0.818) (1.211, 0.815) (1.231, 0.812) (1.251, 0.809) (1.271, 0.806) (1.291, 0.803) (1.311, 0.800) (1.331, 0.797) (1.351, 0.794) (1.371, 0.791) (1.391, 0.788) (1.411, 0.785) (1.431, 0.783) (1.451, 0.780) (1.471, 0.777) (1.491, 0.774) (1.511, 0.771) (1.531, 0.769) (1.551, 0.766) (1.571, 0.763) (1.591, 0.760) (1.611, 0.757) (1.631, 0.755) (1.651, 0.752) (1.671, 0.749) (1.691, 0.747) (1.711, 0.744) (1.731, 0.741) (1.751, 0.739) (1.771, 0.736) (1.791, 0.733) (1.811, 0.731) (1.831, 0.728) (1.851, 0.726) (1.871, 0.723) (1.891, 0.721) (1.911, 0.718) (1.931, 0.716) (1.951, 0.713) (1.971, 0.711) (1.991, 0.708) (2.011, 0.706) (2.031, 0.703) (2.051, 0.701) (2.071, 0.698) (2.091, 0.696) (2.111, 0.694) (2.131, 0.691) (2.151, 0.689) (2.171, 0.687) (2.191, 0.684) (2.211, 0.682) (2.231, 0.680) (2.251, 0.677) (2.271, 0.675) (2.291, 0.673) (2.311, 0.671) (2.331, 0.668) (2.351, 0.666) (2.371, 0.664) (2.391, 0.662) (2.411, 0.660) (2.431, 0.657) (2.451, 0.655) (2.471, 0.653) (2.491, 0.651) (2.511, 0.649) (2.531, 0.647) (2.551, 0.645) (2.571, 0.643) (2.591, 0.640) (2.611, 0.638) (2.631, 0.636) (2.651, 0.634) (2.671, 0.632) (2.691, 0.630) (2.711, 0.628) (2.731, 0.626) (2.751, 0.624) (2.771, 0.622) (2.791, 0.620) (2.811, 0.618) (2.831, 0.616) (2.851, 0.615) (2.871, 0.613) (2.891, 0.611) (2.911, 0.609) (2.931, 0.607) (2.951, 0.605) (2.971, 0.603) (2.991, 0.601) (3.012, 0.599) (3.032, 0.598) (3.052, 0.596) (3.072, 0.594) (3.092, 0.592) (3.112, 0.590) (3.132, 0.589) (3.152, 0.587) (3.172, 0.585) (3.192, 0.583) (3.212, 0.582) (3.232, 0.580) (3.252, 0.578) (3.272, 0.576) (3.292, 0.575) (3.312, 0.573) (3.332, 0.571) (3.352, 0.570) (3.372, 0.568) (3.392, 0.566) (3.412, 0.565) (3.432, 0.563) (3.452, 0.561) (3.472, 0.560) (3.492, 0.558) (3.512, 0.557) (3.532, 0.555) (3.552, 0.553) (3.572, 0.552) (3.592, 0.550) (3.612, 0.549) (3.632, 0.547) (3.652, 0.546) (3.672, 0.544) (3.692, 0.542) (3.712, 0.541) (3.732, 0.539) (3.752, 0.538) (3.772, 0.536) (3.792, 0.535) (3.812, 0.533) (3.832, 0.532) (3.852, 0.531) (3.872, 0.529) (3.892, 0.528) (3.912, 0.526) (3.932, 0.525) (3.952, 0.523) (3.972, 0.522) (3.992, 0.520) (4.012, 0.519) (4.032, 0.518) (4.052, 0.516) (4.072, 0.515) (4.092, 0.514) (4.112, 0.512) (4.132, 0.511) (4.152, 0.509) (4.172, 0.508) (4.192, 0.507) (4.212, 0.505) (4.232, 0.504) (4.252, 0.503) (4.272, 0.501) (4.292, 0.500) (4.312, 0.499) (4.332, 0.497) (4.352, 0.496) (4.372, 0.495) (4.392, 0.494) (4.412, 0.492) (4.432, 0.491) (4.452, 0.490) (4.472, 0.489) (4.492, 0.487) (4.512, 0.486) (4.532, 0.485) (4.552, 0.484) (4.572, 0.482) (4.592, 0.481) (4.612, 0.480) (4.632, 0.479) (4.652, 0.478) (4.672, 0.476) (4.692, 0.475) (4.712, 0.474) (4.732, 0.473) (4.752, 0.472) (4.772, 0.470) (4.792, 0.469) (4.812, 0.468) (4.832, 0.467) (4.852, 0.466) (4.872, 0.465) (4.892, 0.464) (4.912, 0.462) (4.932, 0.461) (4.952, 0.460) (4.972, 0.459) (4.992, 0.458) (5.013, 0.457) (5.033, 0.456) (5.053, 0.455) (5.073, 0.454) (5.093, 0.453) (5.113, 0.451) (5.133, 0.450) (5.153, 0.449) (5.173, 0.448) (5.193, 0.447) (5.213, 0.446) (5.233, 0.445) (5.253, 0.444) (5.273, 0.443) (5.293, 0.442) (5.313, 0.441) (5.333, 0.440) (5.353, 0.439) (5.373, 0.438) (5.393, 0.437) (5.413, 0.436) (5.433, 0.435) (5.453, 0.434) (5.473, 0.433) (5.493, 0.432) (5.513, 0.431) (5.533, 0.430) (5.553, 0.429) (5.573, 0.428) (5.593, 0.427) (5.613, 0.426) (5.633, 0.425) (5.653, 0.424) (5.673, 0.423) (5.693, 0.422) (5.713, 0.421) (5.733, 0.420) (5.753, 0.419) (5.773, 0.419) (5.793, 0.418) (5.813, 0.417) (5.833, 0.416) (5.853, 0.415) (5.873, 0.414) (5.893, 0.413) (5.913, 0.412) (5.933, 0.411) (5.953, 0.410) (5.973, 0.409) (5.993, 0.409) (6.013, 0.408) (6.033, 0.407) (6.053, 0.406) (6.073, 0.405) (6.093, 0.404) (6.113, 0.403) (6.133, 0.402) (6.153, 0.402) (6.173, 0.401) (6.193, 0.400) (6.213, 0.399) (6.233, 0.398) (6.253, 0.397) (6.273, 0.397) (6.293, 0.396) (6.313, 0.395) (6.333, 0.394) (6.353, 0.393) (6.373, 0.392) (6.393, 0.392) (6.413, 0.391) (6.433, 0.390) (6.453, 0.389) (6.473, 0.388) (6.493, 0.388) (6.513, 0.387) (6.533, 0.386) (6.553, 0.385) (6.573, 0.384) (6.593, 0.384) (6.613, 0.383) (6.633, 0.382) (6.653, 0.381) (6.673, 0.380) (6.693, 0.380) (6.713, 0.379) (6.733, 0.378) (6.753, 0.377) (6.773, 0.377) (6.793, 0.376) (6.813, 0.375) (6.833, 0.374) (6.853, 0.374) (6.873, 0.373) (6.893, 0.372) (6.913, 0.371) (6.933, 0.371) (6.953, 0.370) (6.973, 0.369) (6.993, 0.369) (7.014, 0.368) (7.034, 0.367) (7.054, 0.366) (7.074, 0.366) (7.094, 0.365) (7.114, 0.364) (7.134, 0.364) (7.154, 0.363) (7.174, 0.362) (7.194, 0.361) (7.214, 0.361) (7.234, 0.360) (7.254, 0.359) (7.274, 0.359) (7.294, 0.358) (7.314, 0.357) (7.334, 0.357) (7.354, 0.356) (7.374, 0.355) (7.394, 0.355) (7.414, 0.354) (7.434, 0.353) (7.454, 0.353) (7.474, 0.352) (7.494, 0.351) (7.514, 0.351) (7.534, 0.350) (7.554, 0.349) (7.574, 0.349) (7.594, 0.348) (7.614, 0.347) (7.634, 0.347) (7.654, 0.346) (7.674, 0.345) (7.694, 0.345) (7.714, 0.344) (7.734, 0.344) (7.754, 0.343) (7.774, 0.342) (7.794, 0.342) (7.814, 0.341) (7.834, 0.340) (7.854, 0.340) (7.874, 0.339) (7.894, 0.339) (7.914, 0.338) (7.934, 0.337) (7.954, 0.337) (7.974, 0.336) (7.994, 0.336) (8.014, 0.335) (8.034, 0.334) (8.054, 0.334) (8.074, 0.333) (8.094, 0.333) (8.114, 0.332) (8.134, 0.331) (8.154, 0.331) (8.174, 0.330) (8.194, 0.330) (8.214, 0.329) (8.234, 0.329) (8.254, 0.328) (8.274, 0.327) (8.294, 0.327) (8.314, 0.326) (8.334, 0.326) (8.354, 0.325) (8.374, 0.325) (8.394, 0.324) (8.414, 0.323) (8.434, 0.323) (8.454, 0.322) (8.474, 0.322) (8.494, 0.321) (8.514, 0.321) (8.534, 0.320) (8.554, 0.320) (8.574, 0.319) (8.594, 0.318) (8.614, 0.318) (8.634, 0.317) (8.654, 0.317) (8.674, 0.316) (8.694, 0.316) (8.714, 0.315) (8.734, 0.315) (8.754, 0.314) (8.774, 0.314) (8.794, 0.313) (8.814, 0.313) (8.834, 0.312) (8.854, 0.312) (8.874, 0.311) (8.894, 0.311) (8.914, 0.310) (8.934, 0.310) (8.954, 0.309) (8.974, 0.308) (8.994, 0.308) (9.015, 0.307) (9.035, 0.307) (9.055, 0.306) (9.075, 0.306) (9.095, 0.305) (9.115, 0.305) (9.135, 0.304) (9.155, 0.304) (9.175, 0.303) (9.195, 0.303) (9.215, 0.302) (9.235, 0.302) (9.255, 0.302) (9.275, 0.301) (9.295, 0.301) (9.315, 0.300) (9.335, 0.300) (9.355, 0.299) (9.375, 0.299) (9.395, 0.298) (9.415, 0.298) (9.435, 0.297) (9.455, 0.297) (9.475, 0.296) (9.495, 0.296) (9.515, 0.295) (9.535, 0.295) (9.555, 0.294) (9.575, 0.294) (9.595, 0.293) (9.615, 0.293) (9.635, 0.293) (9.655, 0.292) (9.675, 0.292) (9.695, 0.291) (9.715, 0.291) (9.735, 0.290) (9.755, 0.290) (9.775, 0.289) (9.795, 0.289) (9.815, 0.288) (9.835, 0.288) (9.855, 0.288) (9.875, 0.287) (9.895, 0.287) (9.915, 0.286) (9.935, 0.286) (9.955, 0.285) (9.975, 0.285) (9.995, 0.285) (10.015, 0.284) (10.035, 0.284) (10.055, 0.283) (10.075, 0.283) (10.095, 0.282) (10.115, 0.282) (10.135, 0.282) (10.155, 0.281) (10.175, 0.281) (10.195, 0.280) (10.215, 0.280) (10.235, 0.279) (10.255, 0.279) (10.275, 0.279) (10.295, 0.278) (10.315, 0.278) (10.335, 0.277) (10.355, 0.277) (10.375, 0.277) (10.395, 0.276) (10.415, 0.276) (10.435, 0.275) (10.455, 0.275) (10.475, 0.274) (10.495, 0.274) (10.515, 0.274) (10.535, 0.273) (10.555, 0.273) (10.575, 0.272) (10.595, 0.272) (10.615, 0.272) (10.635, 0.271) (10.655, 0.271) (10.675, 0.271) (10.695, 0.270) (10.715, 0.270) (10.735, 0.269) (10.755, 0.269) (10.775, 0.269) (10.795, 0.268) (10.815, 0.268) (10.835, 0.267) (10.855, 0.267) (10.875, 0.267) (10.895, 0.266) (10.915, 0.266) (10.935, 0.265) (10.955, 0.265) (10.975, 0.265) (10.995, 0.264) (11.016, 0.264) (11.036, 0.264) (11.056, 0.263) (11.076, 0.263) (11.096, 0.263) (11.116, 0.262) (11.136, 0.262) (11.156, 0.261) (11.176, 0.261) (11.196, 0.261) (11.216, 0.260) (11.236, 0.260) (11.256, 0.260) (11.276, 0.259) (11.296, 0.259) (11.316, 0.259) (11.336, 0.258) (11.356, 0.258) (11.376, 0.257) (11.396, 0.257) (11.416, 0.257) (11.436, 0.256) (11.456, 0.256) (11.476, 0.256) (11.496, 0.255) (11.516, 0.255) (11.536, 0.255) (11.556, 0.254) (11.576, 0.254) (11.596, 0.254) (11.616, 0.253) (11.636, 0.253) (11.656, 0.253) (11.676, 0.252) (11.696, 0.252) (11.716, 0.252) (11.736, 0.251) (11.756, 0.251) (11.776, 0.251) (11.796, 0.250) (11.816, 0.250) (11.836, 0.250) (11.856, 0.249) (11.876, 0.249) (11.896, 0.249) (11.916, 0.248) (11.936, 0.248) (11.956, 0.248) (11.976, 0.247) (11.996, 0.247) (12.016, 0.247) (12.036, 0.246) (12.056, 0.246) (12.076, 0.246) (12.096, 0.245) (12.116, 0.245) (12.136, 0.245) (12.156, 0.244) (12.176, 0.244) (12.196, 0.244) (12.216, 0.243) (12.236, 0.243) (12.256, 0.243) (12.276, 0.242) (12.296, 0.242) (12.316, 0.242) (12.336, 0.241) (12.356, 0.241) (12.376, 0.241) (12.396, 0.240) (12.416, 0.240) (12.436, 0.240) (12.456, 0.240) (12.476, 0.239) (12.496, 0.239) (12.516, 0.239) (12.536, 0.238) (12.556, 0.238) (12.576, 0.238) (12.596, 0.237) (12.616, 0.237) (12.636, 0.237) (12.656, 0.237) (12.676, 0.236) (12.696, 0.236) (12.716, 0.236) (12.736, 0.235) (12.756, 0.235) (12.776, 0.235) (12.796, 0.234) (12.816, 0.234) (12.836, 0.234) (12.856, 0.234) (12.876, 0.233) (12.896, 0.233) (12.916, 0.233) (12.936, 0.232) (12.956, 0.232) (12.976, 0.232) (12.996, 0.232) (13.017, 0.231) (13.037, 0.231) (13.057, 0.231) (13.077, 0.230) (13.097, 0.230) (13.117, 0.230) (13.137, 0.230) (13.157, 0.229) (13.177, 0.229) (13.197, 0.229) (13.217, 0.228) (13.237, 0.228) (13.257, 0.228) (13.277, 0.228) (13.297, 0.227) (13.317, 0.227) (13.337, 0.227) (13.357, 0.226) (13.377, 0.226) (13.397, 0.226) (13.417, 0.226) (13.437, 0.225) (13.457, 0.225) (13.477, 0.225) (13.497, 0.225) (13.517, 0.224) (13.537, 0.224) (13.557, 0.224) (13.577, 0.223) (13.597, 0.223) (13.617, 0.223) (13.637, 0.223) (13.657, 0.222) (13.677, 0.222) (13.697, 0.222) (13.717, 0.222) (13.737, 0.221) (13.757, 0.221) (13.777, 0.221) (13.797, 0.221) (13.817, 0.220) (13.837, 0.220) (13.857, 0.220) (13.877, 0.219) (13.897, 0.219) (13.917, 0.219) (13.937, 0.219) (13.957, 0.218) (13.977, 0.218) (13.997, 0.218) (14.017, 0.218) (14.037, 0.217) (14.057, 0.217) (14.077, 0.217) (14.097, 0.217) (14.117, 0.216) (14.137, 0.216) (14.157, 0.216) (14.177, 0.216) (14.197, 0.215) (14.217, 0.215) (14.237, 0.215) (14.257, 0.215) (14.277, 0.214) (14.297, 0.214) (14.317, 0.214) (14.337, 0.214) (14.357, 0.213) (14.377, 0.213) (14.397, 0.213) (14.417, 0.213) (14.437, 0.212) (14.457, 0.212) (14.477, 0.212) (14.497, 0.212) (14.517, 0.212) (14.537, 0.211) (14.557, 0.211) (14.577, 0.211) (14.597, 0.211) (14.617, 0.210) (14.637, 0.210) (14.657, 0.210) (14.677, 0.210) (14.697, 0.209) (14.717, 0.209) (14.737, 0.209) (14.757, 0.209) (14.777, 0.208) (14.797, 0.208) (14.817, 0.208) (14.837, 0.208) (14.857, 0.207) (14.877, 0.207) (14.897, 0.207) (14.917, 0.207) (14.937, 0.207) (14.957, 0.206) (14.977, 0.206) (14.997, 0.206) (15.018, 0.206) (15.038, 0.205) (15.058, 0.205) (15.078, 0.205) (15.098, 0.205) (15.118, 0.205) (15.138, 0.204) (15.158, 0.204) (15.178, 0.204) (15.198, 0.204) (15.218, 0.203) (15.238, 0.203) (15.258, 0.203) (15.278, 0.203) (15.298, 0.203) (15.318, 0.202) (15.338, 0.202) (15.358, 0.202) (15.378, 0.202) (15.398, 0.201) (15.418, 0.201) (15.438, 0.201) (15.458, 0.201) (15.478, 0.201) (15.498, 0.200) (15.518, 0.200) (15.538, 0.200) (15.558, 0.200) (15.578, 0.199) (15.598, 0.199) (15.618, 0.199) (15.638, 0.199) (15.658, 0.199) (15.678, 0.198) (15.698, 0.198) (15.718, 0.198) (15.738, 0.198) (15.758, 0.198) (15.778, 0.197) (15.798, 0.197) (15.818, 0.197) (15.838, 0.197) (15.858, 0.197) (15.878, 0.196) (15.898, 0.196) (15.918, 0.196) (15.938, 0.196) (15.958, 0.195) (15.978, 0.195) (15.998, 0.195) (16.018, 0.195) (16.038, 0.195) (16.058, 0.194) (16.078, 0.194) (16.098, 0.194) (16.118, 0.194) (16.138, 0.194) (16.158, 0.193) (16.178, 0.193) (16.198, 0.193) (16.218, 0.193) (16.238, 0.193) (16.258, 0.192) (16.278, 0.192) (16.298, 0.192) (16.318, 0.192) (16.338, 0.192) (16.358, 0.191) (16.378, 0.191) (16.398, 0.191) (16.418, 0.191) (16.438, 0.191) (16.458, 0.190) (16.478, 0.190) (16.498, 0.190) (16.518, 0.190) (16.538, 0.190) (16.558, 0.190) (16.578, 0.189) (16.598, 0.189) (16.618, 0.189) (16.638, 0.189) (16.658, 0.189) (16.678, 0.188) (16.698, 0.188) (16.718, 0.188) (16.738, 0.188) (16.758, 0.188) (16.778, 0.187) (16.798, 0.187) (16.818, 0.187) (16.838, 0.187) (16.858, 0.187) (16.878, 0.186) (16.898, 0.186) (16.918, 0.186) (16.938, 0.186) (16.958, 0.186) (16.978, 0.186) (16.998, 0.185) (17.019, 0.185) (17.039, 0.185) (17.059, 0.185) (17.079, 0.185) (17.099, 0.184) (17.119, 0.184) (17.139, 0.184) (17.159, 0.184) (17.179, 0.184) (17.199, 0.184) (17.219, 0.183) (17.239, 0.183) (17.259, 0.183) (17.279, 0.183) (17.299, 0.183) (17.319, 0.182) (17.339, 0.182) (17.359, 0.182) (17.379, 0.182) (17.399, 0.182) (17.419, 0.182) (17.439, 0.181) (17.459, 0.181) (17.479, 0.181) (17.499, 0.181) (17.519, 0.181) (17.539, 0.180) (17.559, 0.180) (17.579, 0.180) (17.599, 0.180) (17.619, 0.180) (17.639, 0.180) (17.659, 0.179) (17.679, 0.179) (17.699, 0.179) (17.719, 0.179) (17.739, 0.179) (17.759, 0.179) (17.779, 0.178) (17.799, 0.178) (17.819, 0.178) (17.839, 0.178) (17.859, 0.178) (17.879, 0.178) (17.899, 0.177) (17.919, 0.177) (17.939, 0.177) (17.959, 0.177) (17.979, 0.177) (17.999, 0.177) (18.019, 0.176) (18.039, 0.176) (18.059, 0.176) (18.079, 0.176) (18.099, 0.176) (18.119, 0.176) (18.139, 0.175) (18.159, 0.175) (18.179, 0.175) (18.199, 0.175) (18.219, 0.175) (18.239, 0.175) (18.259, 0.174) (18.279, 0.174) (18.299, 0.174) (18.319, 0.174) (18.339, 0.174) (18.359, 0.174) (18.379, 0.173) (18.399, 0.173) (18.419, 0.173) (18.439, 0.173) (18.459, 0.173) (18.479, 0.173) (18.499, 0.172) (18.519, 0.172) (18.539, 0.172) (18.559, 0.172) (18.579, 0.172) (18.599, 0.172) (18.619, 0.172) (18.639, 0.171) (18.659, 0.171) (18.679, 0.171) (18.699, 0.171) (18.719, 0.171) (18.739, 0.171) (18.759, 0.170) (18.779, 0.170) (18.799, 0.170) (18.819, 0.170) (18.839, 0.170) (18.859, 0.170) (18.879, 0.169) (18.899, 0.169) (18.919, 0.169) (18.939, 0.169) (18.959, 0.169) (18.979, 0.169) (18.999, 0.169) (19.020, 0.168) (19.040, 0.168) (19.060, 0.168) (19.080, 0.168) (19.100, 0.168) (19.120, 0.168) (19.140, 0.167) (19.160, 0.167) (19.180, 0.167) (19.200, 0.167) (19.220, 0.167) (19.240, 0.167) (19.260, 0.167) (19.280, 0.166) (19.300, 0.166) (19.320, 0.166) (19.340, 0.166) (19.360, 0.166) (19.380, 0.166) (19.400, 0.166) (19.420, 0.165) (19.440, 0.165) (19.460, 0.165) (19.480, 0.165) (19.500, 0.165) (19.520, 0.165) (19.540, 0.165) (19.560, 0.164) (19.580, 0.164) (19.600, 0.164) (19.620, 0.164) (19.640, 0.164) (19.660, 0.164) (19.680, 0.163) (19.700, 0.163) (19.720, 0.163) (19.740, 0.163) (19.760, 0.163) (19.780, 0.163) (19.800, 0.163) (19.820, 0.162) (19.840, 0.162) (19.860, 0.162) (19.880, 0.162) (19.900, 0.162) (19.920, 0.162) (19.940, 0.162) (19.960, 0.162) (19.980, 0.161) (20.000, 0.161)
    };
  \end{axis}
\end{tikzpicture}

%% file: thesis-arxiv.bbl
\clearpage
\newpage

\begin{thebibliography}{99} 
	\bibitem{errwpemantle}
  Robin Pemantle.
	\textit{Phase Transition in Reinforced Random Walk and RWRE on Trees}.
	The Annals of Probability, Vol.~16, No.~3 (Jul., 1988), pp.~1229-1241. Institute of Mathematical Statistics, 1988.
	Accessible at \url{https://projecteuclid.org/euclid.aop/1176991687}.
	
	\bibitem{rrwrenlund}
	Henrik Renlund.
	\textit{Reinforced Random Walk}.
  Department of Mathematics, Uppsala University. 2005.
	Accessible at \url{http://www.diva-portal.org/smash/record.jsf?pid=diva2%3A305167&dswid=6914}

	\bibitem{rwrelyonspemantle}
	Russell Lyons and Robin Pemantle.
	\textit{Random Walk in a Random Environment and First-Passage Percolation on Trees}.
  The Annals of Probability, Vol.~20, No.~1 (Jan., 1992), pp.~125-136. Institute of Mathematical Statistics, 1992.
  Accessible at \url{https://projecteuclid.org/euclid.aop/1176989920}.
  
	\bibitem{rwrelyonspemantlecorr}
	Russell Lyons and Robin Pemantle.
	\textit{Correction: Random walk in a random environment and first-passage percolation on trees}.
  The Annals of Probability, Vol.~31, No.~1 (Jan., 2003), pp.~528-529. Institute of Mathematical Statistics, 2003.
  Accessible at \url{https://projecteuclid.org/euclid.aop/1046294319}.
  
	\bibitem{transientaidekon}
	Elie Aid\'ekon.
	\textit{Transient random walks in random environment on a Galton-Watson tree}.
  Probability Theory and Related Fields, Vol.~142, No.~3-4 (Nov., 2008), pp.~525-559. Springer, 2008.
	Accessible at \url{https://arxiv.org/abs/0710.3377}.
	
	\bibitem{probtreenet}
  Russell Lyons and Yuval Peres.
  \textit{Probability on Trees and Networks}.
  Version of 6 August 2019.
  Cambridge University Press, 2017.
	Accessible at \url{https://rdlyons.pages.iu.edu/prbtree/} (online-only corrected edition at
	\url{https://rdlyons.pages.iu.edu/prbtree/book_corr.pdf}).
	
  \bibitem{introrw}
  Pablo Lessa.
	\textit{Recurrence vs Transience: An introduction to random walks}.
  2015.
	Accessible at \url{http://www.cmat.edu.uy/~lessa/otherwork.html} under \url{http://www.cmat.edu.uy/~lessa/resource/randomwalknotes.pdf}.
	
	\bibitem{ising}
	Russell Lyons.
	\textit{The Ising Model and Percolation on Trees and Tree-Like Graphs}.
	Communications in Mathematical Physics, Vol.~125, No.~2 (Jun., 1989), pp.~337-353. Springer, 1989.
	Accessible at \url{https://projecteuclid.org/euclid.cmp/1104179469}.
	
	\bibitem{rwperctrees}
	Russell Lyons.
	\textit{Random Walks and Percolation on Trees}.
	The Annals of Probability, Vol.~18, No.~3 (Jul., 1990), pp.~931-958. Institute of Mathematical Statistics, 1990.
	Accessible at \url{https://projecteuclid.org/euclid.aop/1176990730}.
	
	\bibitem{cramerthm}
	Rapha\"el Cerf and Pierre Petit.
	\textit{A short proof of Cram\'er's theorem in $\bbR$}.
	The American Mathematical Monthly, Vol.~118, No.~10 (Dec., 2011), pp.~925-931. Mathematical Association of America, 2011.
	Accessible at \url{https://arxiv.org/abs/1002.3496}.
\end{thebibliography}
